\documentclass[11pt, a4paper]{amsart}
\usepackage{amsfonts}
\usepackage{amsmath}
\usepackage{enumerate}
\newcommand{\equ}[1]{(\ref{#1})} 
 
\newcommand{\be}{\begin{equation}} 
\newcommand{\ee}{\end{equation}} 
\newtheorem{lem}{Lemma}[section] 
\newcommand{\rr}{{\mathbb R}} 
\newcommand{\nn}{{\mathbb N}} 
\newtheorem{teo}{Theorem}[section] 
\newtheorem{hypo}{Hypothesis}
\newtheorem*{claim}{Claim} 
\newtheorem{prop}{Proposition}[section] 
\newtheorem{rem}{Remark}[section] 
\newtheorem{cor}{Corollary}[section] 
\newcommand{\calA}{{\mathcal A}} 
 
\newcommand{\calC}{{\mathcal C}}

\newcommand{\calF}{{\mathcal F}} 
\newcommand{\calG}{{\mathcal G}} 
\newcommand{\calH}{{\mathcal H}} 
\newcommand{\calI}{{\mathcal I}} 
\newcommand{\calK}{{\mathcal K}} 
\newcommand{\calL}{{\mathcal L}} 
\newcommand{\calM}{{\mathcal M}}

\newcommand{\calQ}{{\mathcal Q}}

\newcommand{\calV}{{\mathcal V}} 
\newcommand{\calJ}{{\mathcal J}} 

\begin{document} 
\bibliographystyle{plain} 
\title[Allen-Cahn equation with a  
three well potential]{Existence of a solution to a vector-valued  
Allen-Cahn equation with a  
three well potential}

\author[ M. Saez Trumper]
{Mariel Saez Trumper}

\address{Mariel Saez Trumper
\hfill\break\indent
Max Planck Institute for Gravitational Physics
\hfill\break\indent
Albert Einstein Institute
\hfill\break\indent
Am M{\"u}hlenberg 1\\
D-14476 Golm\\
Germany\\
\hfill\break\indent
Germany .}
\email{{\tt  mariel.saez@aei.mpg.de}}
\begin{abstract}
In this paper we prove existence of a vector-valued solution $v$ to
\begin{align*} -\Delta v +\frac{\nabla_v W(v)}{2}&=0, \\
 \lim_{r\to \infty}v(r \cos\theta,r\sin\theta)&= c_i \hbox{ for } \theta \in 
(\theta_{i-1}, \theta_i), \end{align*}
where $W:\rr^2\to \rr$ is non-negative function that attains its minimum 0 at
$\{c_i\}_{i=1}^3$ and the angles $\theta_i$ are determined by the function $W$.
This solution is an energy minimizer.\end{abstract}

\maketitle

\section{Introduction}

In this paper we establish existence of a vector-valued solution $v:\rr^2\to \rr^2$
to the following elliptic problem:

\begin{align} -\Delta v +\frac{\nabla_v W( v)}{2}&=0, \label{ecprinc}\\
\lim_{r\to \infty}v(r \cos\theta,r\sin\theta)&= c_i \hbox{ for }\theta \in 
(\theta_{i-1}, \theta_i), \label{condainf}\end{align}
where $W:\rr^2\to \rr$ is positive function with three local minima, given by $\{c_i\}_{i=1}^3$, and the angles $\theta_i$,   with $\theta_3=2\pi+\theta_0$, are determined by the potential $W$  (for a more precise description on how these angles are determined we refer the reader to definitions \equ{condang} and \equ{defthetai}). 
%The proof of this result is based on geometric measure theory arguments, variational arguments and parabolic equations techniques. 

In \cite{vecvalps} an analogous result was proved by P. Sternberg when  $W$ has two minima. Later on, Bronsard, Gui and Schatzman~\cite{atrhlb} 
considered potentials with three minima that were  equivariant under the symmetry group of the equilateral triangle. Under these conditions they 
proved  existence of a solution to \equ{ecprinc}-\equ{condainf}.  The system of equations given by
 \equ{ecprinc} was also
studied in \cite{hegf}, but the domains considered  were
 bounded and  Neumann boundary condition  was imposed. In that paper, under appropriate assumptions over the potential $W$, Flores, Padilla and Tonwaga established the existence of solutions that join the three minima ($c_1$, $c_2$ and $c_3$); however, no precise description of the triple junction was provided. 
 Recently,  potentials with four minima were studied in  \cite{sqcg}, establishing (under several assumptions over the potential $W$) the existence of solutions to \equ{ecprinc} that connect all the four wells.

Our interest in this problem is originated in some models of three-boundary
motion. Material scientists working on the theory of  transition layers have found that the motion of grain boundaries is governed by its local mean curvature 
(see \cite{twodimwm},\cite{twodimwm2} for example). These models naturally arise as the singular limit of the parabolic Allen-Cahn equation (see \cite{amicsa}). 
The expected relation between  grain boundaries motion and the  parabolic Allen-Cahn  equation can be described as follows:  Consider a positive potential $W:\Omega\subset \rr^n\to \rr$ with a finite number of minima $\{c_i\}_{i=i}^m$. Let $u_\epsilon:\rr^n\to\rr^n$ be a solution to 
 \be \frac{d u_\epsilon}{d t}-\Delta u_\epsilon +\frac{\nabla_v W(u_\epsilon)}{2\epsilon^2}=0 .\label{ginzlandpar}\ee
As $\epsilon \to 0$ the solutions $u_\epsilon$ will converge almost everywhere to one of the constants $c_i$ (see \cite{conofti}, \cite{fasreajr}).  For every $t$, this creates a  partition of $\Omega=\bigcup_{i=1}^m\Omega_i(t)$, where $\Omega_i(t)=\{x\in\Omega:u_\epsilon(x,t)\to c_i \hbox{ as }\epsilon \to 0\}$. The interface between these sets  correspond to the grain boundaries, which   evolve under its curvature.               
When $n=2$ and $m=3$ the solution will describe a ``three-phase'' boundary motion that might present ``triple- points'', namely points where  these 3 boundaries meet.
Bronsard and Reitich~\cite{onthrlb} predicted that at points that are away from the triple points and close to 
the interface between $c_i$ and $c_j$ the solutions to \equ{ginzlandpar} should be approximated by $\zeta_{ij}\left(\frac{d_{ij}(x,t)}{\epsilon}\right)$, where 
$d_{ij}$ is the distance function to this interface and $\zeta_{ij}$ is a solution to the equation
\be \zeta_{ij}''(\lambda)+\frac{\nabla W(\zeta_{ij}(\lambda))}{2}=0, \label{eczeta}.\ee
\be \lim_{\tau\to -\infty}\zeta_{ij}(\tau)=c_i, 
\lim_{\tau\to \infty}\zeta_{ij}(\tau)=c_j. \label{explim}\ee
On the other hand, close to the triple points it is expected that
solutions to
\equ{ginzlandpar}  will behave, after appropriate rescaling, like a solution to \equ{ecprinc}-\equ{condainf}. However, the existence of such solution has not been established  in the general case before. This is the main goal of this paper. 

Based on the previous discussion, 
in order to match the expected behavior of solutions to \equ{ginzlandpar}  near double junctions  and the one close to  triple junctions, we expect that  solutions to \equ{ecprinc}-\equ{condainf}  satisfy an extra condition at infinity. 
Namely, solutions to \equ{ecprinc}-\equ{condainf}  should  resemble solutions to   \equ{eczeta}-\equ{explim} near the half-lines of direction $\theta_i$.
%Namely,
 %solutions to \equ{ecprinc} should 
 %not only  satisfy \equ{condainf}, but also,   close to the half-lines starting at the origin  with  direction  $ \theta_i$, these solutions should  approach, in an adequate sense,  appropriate  solutions to \equ{eczeta}. 
 We will implicitly impose this condition throughout the paper. 
Therefore, we briefly discuss  the existence of solutions to  \equ{eczeta}-\equ{explim}.
% However, in order to impose the extra condition, it is necessary to study the existence of solutions to  \equ{eczeta}-\equ{explim}.
  %In order to impose the appropriate behavior, we start by   briefly  discussing the existence of  solutions
 %\equ{eczeta}-\equ{explim}:  
% Hence, the existence of such curves will be beneficial in the proof of our result.
For potentials with two wells 
 the 
 existence of such curves was proved   by P. Sternberg in
\cite{vecvalps}.
 %the conditions the following conditions:
 %\begin{enumerate}
 %\item \label{cond1W} the matrix $\frac{\partial^2 W(u)}{\partial u_i \partial u_j}$ is positive semidefinite at $\{c_i\}_{i=1}^3$, that is the minima are non-degenerate;
%\item \label{cond2W} there exist positive constants $K_1, K_2$ and $m$, and a number $p\geq 2$ such that
%$$K_1|u|^p\leq W(u)\leq K_2 |u|^p \hbox{ for } |u|\geq m;$$
%\item \label{cond3W} $V(r,\theta):=W(u+r(\cos \theta, \sin \theta))=r^2+O(r^3)$
%for $r$ sufficiently small and $u=c_i$ for some $i\in\{1,2,3\}$, where $r$ and $\theta$ are local polar coordinates.
%\end{enumerate} 
However,   the problem is more subtle  when considering arbitrary three-well potentials, even if conditions analogous to the ones imposed in  \cite{vecvalps}  hold. In
\cite{expstana} Alikakos, Betel\'u and Chen  provided some examples of potentials where  solutions to  \equ{eczeta}-\equ{explim} did not exist for certain $i,j$. 
%proved that,   analogous might not be sufficient to ensure the existence of   curve solutions to equations for every . 
On the other hand, they also established
appropriate conditions under which all these solutions in fact do exist. 
In what follows  we will  assume  we are in the latter case. Namely we assume the existence of $\zeta_{ij}$ for every $i$ and $j$. This and other technical assumptions on the potential $W$ will be discussed in detail in the following section. At the moment we  state the main theorem of this paper:
%More specifically we prove thata
\begin{teo}\label{teoprinc}
Let $W:\rr^2\to \rr$ be a proper $\calC^3$ function that satisfy
\begin{enumerate}[ (a) ]
\item\label{condW0} $W$ has only three local minima $c_1, c_2$ and $c_3$ and $W(c_i)=0$; 
\item \label{cond1W} the matrix $\frac{\partial^2 W(u)}{\partial u_i \partial u_j}$ is positive definite at $\{c_i\}_{i=1}^3$, that is the minima are non-degenerate.
 
 \item \label{condhess} The hessian of the function $W(u)$ (which we denote by $W''$) is positive semidefinite for $|u|>K$, where  $K>0$ is a fixed real number;
\item \label{cond2W} there exist positive constants $K_1, K_2$ and $m$, and a number $p\geq 2$ such that
$$K_1|u|^p\leq W(u)\leq K_2 |u|^p \hbox{ for } |u|\geq m;$$
\item \label{cond3W}  Hypothesis \ref{geod} holds (see the next section for a description on this hypothesis). In particular, there are solutions to \equ{eczeta}-\equ{explim} for every $i$ and $j$.
 \end{enumerate} 
Define 
\be \Gamma (\zeta_1,\zeta_2)=\inf\left\{\int_0^1 W^{\frac{1}{2}}(\gamma(\lambda))|\gamma'(\lambda)|d\lambda:  \gamma \in C^1([0,1],\rr^2), \gamma(0)=\zeta_1\hbox{ and } \gamma(1)=\zeta_2 \right\} \label{defdist}\ee
Consider $\{\alpha_i\}_{i=1}^3\in[0,2\pi)$ such that
\be \frac{\sin \alpha_1}{\Gamma(c_2,c_3)}= \frac{\sin \alpha_2}{\Gamma(c_1,c_3)}= \frac{\sin \alpha_3}{\Gamma(c_1,c_2)}.\label{condang}\ee
Then for $\theta_i\in [0, 2\pi)$ such that \be \alpha_i=\theta_{i}-\theta_{i-1}\label{defthetai}\ee there is a solution $v$ to \equ{ecprinc}-\equ{condainf}. Moreover, there exists  a differentiable function 
$\phi$ satisfying  \equ{condainf} such that 
for
$$\calG (w)=\int_{\rr^2}(|Dw|^2+W(w)-|D\phi|^2-W(\phi))dx,$$
we have
$$\calG (v)=\inf\{\calG (w):w\in \calV\},$$
for $\calV=\left\{w\in C^1:\int_{\rr^2}|Dw-D\phi|dx,\int_{\rr^2}|w-\phi|dx
<\infty\right\}.$ 
\end{teo}

\medskip
We  would like to remark that  the function $\phi$ in  Theorem \ref{teoprinc} will be defined explicitly in the coming section (more specifically in sub-section \ref{constphi}) and it  will capture the behavior at infinity of the solution  $u$ to \equ{ecprinc}-\equ{condainf}. 
In the construction of this function Hypothesis  \equ{cond3W} is required. 
Relaxations of this hypothesis are possible, but we will skip them in order to keep the presentation simpler.
We also want to point out that, as discussed in \cite{onthrlb}, the definitions of $\alpha_i$  and $\Gamma_i$ imply that $\alpha_1+\alpha_2+\alpha_3=2\pi$ and  $\theta_3=2\pi+\theta_0$.

Before proceeding to the coming sections, we would like to briefly  outline our proof of Theorem  \ref{teoprinc}  and its organization through the paper. The basic idea is the following:
Let $B_R$ denote the ball or radius $R$ and let $v_R$ solve equation \equ{ecprinc}  in $B_R$ with Dirichlet boundary condition $v_R |_{\partial B_r}=\phi$
 (the function $\phi$ is defined in equation \equ{deffi}  and captures the desired behaviour at infinity, as it is discussed in Remark  \ref{discfi} below).
 %We take a sequence of functions $v_R$ that solve equation \equ{ecprinc} in a ball of radius $R$  (that we will denote by $B_R$). At the boundary of the ball we impose Dirichlet boundary condition. More precisely, we construct  a continuous function $\phi$ defined in $\rr^2$, which at infinity  captures the desired behavior of the solution $v$ to \equ{ecprinc}-\equ{condainf} (for a more detailed discussion of the construction of the function $\phi$ we refer the reader to \equ{deffi} in the following section and the related discussion).  We impose $v_R(x)=\phi(x) $ for $x\in \partial B_R$. 
 The proof of Theorem \ref{teoprinc} will be equivalent to show convergence of the solutions $v_R$ in an appropriate norm.

In order to prove the convergence result we use the following key observation: In the unit ball we define the function
$$u_R(x)=v_R\left(R x\right),$$ then $u_R$ satisfies
$$-\Delta u_R +\frac{R^2 \nabla_v W(u_R)}{2}=0 \hbox{ for } x\in B_1.$$
Hence for $\epsilon =\frac{1}{R}$, the function $u_\epsilon$ satisfies
\be -\Delta u_\epsilon +\frac{ \nabla_v W(u_\epsilon)}{2\epsilon^2} =0.\label{ginzlandel}\ee
As $R\to \infty$  (or equivalently as $\epsilon\to 0$) we
expect $v_R$ to converge to the solution $v$ to \equ{ecprinc}-\equ{condainf} (this will be proved in Section \ref{pf}), and correspondingly, we expect
  the limiting solution $u_\epsilon$ to \equ{ginzlandel} to capture the behavior of $v$  at infinity.  Equation \equ{ginzlandel} has been largely studied (see for example \cite{ginlanfb} and \cite{linandpr}).
This motivates us to analyze in Section  \ref{convl1} some existing results for \equ{ginzlandel} that apply in our context and  provide useful information for our problem.
More precisely, combining results in   \cite{minintsb}, \cite{locminrk} and \cite{locminps} and using $\Gamma $-convergence techniques
 we prove that the rescaled $u_\epsilon$ converge to  a function $u_0$ in the $L^1$ norm in the unit ball. Moreover, the function $u_0$  equals $c_i$ in the the angular sectors  defined by $\theta\in (\theta_{i-1},\theta_i)$  and it  
 is minimizing for an appropriate functional (eventually, this property will imply the minimizing result in Theorem \ref{teoprinc}).  Hypotheses \equ{cond2W}  and \equ{cond3W} are essential in this section. However, we would like to point out  that it is not clear whether  they are just technical conditions (which may be removed)  or not.  On the other hand, hypotheses   \equ{condW0}  and \equ{cond1W}  (which are also used in this section)  are natural in the context of the problem.

 In order to finish the proof, in 
  Section \ref{unifconv} we show that the convergence holds in a stronger norm than $L^1$. 
The main idea in this computation  is to  use the parabolic version of equation \equ{ginzlandel} to interpolate between an approximate solutions  to \equ{ecprinc} in the ball
  (which we will denote by $U_{\vec{q}}$) and the real solution. More precisely,  we consider a function $\tilde{h}_\epsilon$ that is a solution to
  \begin{align*}
\frac{d\tilde{h}_\epsilon}{dt}-\Delta \tilde{h}_\epsilon +\frac{ \nabla_v W(\tilde{h}_\epsilon)}{2\epsilon^2} &=0 \hbox{ for } x\in B_1, t\in(0,\infty)\\
\tilde{h}(x,t)&=\phi_\epsilon(x) \hbox{ for } x\in \partial B_1, \\
\tilde{h}(x,0)&=U_{\vec{q}}.
  \end{align*}
  The "approximate solution"  $U_{\vec{q}}$ depends on $\epsilon$, satisfies $U_{\vec{q}}=\phi_\epsilon(x)$ for $x\in \partial B_1$ and
   $ \left(-\Delta U_{\vec{q}} +\frac{ \nabla_v W(U_{\vec{q}})}{2\epsilon^2}\right)(x) \to 0$ as $\epsilon\to 0$ point-wise in $B_1$. Using Theorem \ref{aprox} we prove that in fact
   $\tilde{h}_\epsilon$ and $U_{\vec{q}}$ remain appropriately close in time. We conclude by observing that, as $t\to \infty$ it holds that $\tilde{h}_\epsilon(\cdot,t)\to u_\epsilon(\cdot)$. This will imply that  
   in fact $u_\epsilon$ is close to $U_{\vec{q}}$. Also in that section, we use similar techniques to control the convergence in compact domains of the  sequence $v_\epsilon:B_{\frac{1}{\epsilon}}\to \rr^2$
   given by $v_\epsilon(x)=u\left(\epsilon x\right)$.
 The proof of Theorem \ref{teoprinc} can be easily finished by combining  the elements described above. This is achieved in 
  Section \ref{pf}.

We would like to remark that the techniques used in this paper were
already used in similar problems (see \cite{papertesis} and \cite{triodginz}). In general,
the method can be extended as long as the solutions to \equ{ecprinc} 
converge to minima of  $W$ as $\epsilon \to 0 $ and that approximate solutions with the desired characteristics (such as $U_{\vec{q}}$ in this case)  can be constructed.

The author wishes to thank the referee for the very useful comments in improving the exposition, the Max Planck Institute for Gravitational Physics for providing a great work environment and to Rafe Mazzeo and Alex Freire for very useful discussions.

\section{Definitions and preliminary lemmas}\label{defs}

We divide this section into three sub-sections.  The first one is devoted to several definitions that will be used in the analysis performed in Section \ref{convl1}.  The main objective of the second sub-section  is to construct  the function $\phi$ used in Theorem \ref{teoprinc}. In the final sub-section we summarize a collection of existing results that will be used throughout this paper.

\subsection{General definitions}
In this sub-section we will address several general definitions that will simplify the notation in the coming section.

Define the function $g_i:\rr^2\to \rr$  for any $p\in \rr^2$ as
\be g_i(p)= \Gamma(c_i,p)\ee
Where the function $\Gamma$ is defined by \equ{defdist}.
Notice that $\Gamma$ can be regarded as degenerate distance function. Hence
 $g_i(p)$ represents the distance  of a point $p$ (with respect to the distance function $\Gamma$) to the critical point $c_i$.

%Let 
%\begin{align*}g_i(u)=\inf\left \{ \int _{-1}^1\sqrt{W(\gamma(t))}|\gamma'(t)|dt: \right.&
%\gamma:[-1,1]\to \rr^2 \hbox{ that satisfies }\\
%&\left.\gamma(-1)=c_i \hbox{ and } \gamma(1)=u \right\} \end{align*}

Inspired in  \cite{vecvalps} we  consider the following assuption:

\begin{hypo} \label{geod}
Suppose that
for every $u\in \rr^2$, there exists a curve $\gamma_u^i:[-1,1]\to \rr^2$ such that $\gamma_u^i(-1)=c_i$, $\gamma_u^i(1)=u$ and 
\be  g_i(u)= \int _{-1}^1\sqrt{W(\gamma_u^i(t))}\ |\left(\gamma_u^i\right)'(t)|dt. \label{defgi}\ee
The function $g_i$ is Lipschitz continuous and satisfies
\be |D g_i(u)|=\sqrt{W(u)} \hbox{  a.e. }\label{gradgi}\ee
\end{hypo}
For potentials with two wells the 
 existence of such curves was proved  by P. Sternberg in
\cite{vecvalps}. 
He also proved that when considering a curve that joins the minima of $W$, it can be re-parametrized by a curve 
 $\beta_{ij}:(-\infty, \infty)\to (-1,1)$ such that the curves defined by
$$\zeta_{ij}(\tau)= \gamma_{c_j}^i(\beta_{ij}(\tau))$$ satisfy
\be 2g_i(c_j)=\int_{-\infty}^{\infty}W(\zeta_{ij})+|\zeta '|^2 d\tau, \ee
as well as \equ{eczeta} and
\equ{explim} 
%\be \lim_{\tau\to -\infty}\zeta_{ij}(\tau)=c_i, 
%\lim_{\tau\to \infty}\zeta_{ij}(\tau)=c_j, \label{explim}\ee
(where the limits in  \equ{explim}  are attained at an exponential rate).  In our situation, if we assume Hypothesis \ref{geod},
the previous construction can also be carried out. 
%It also holds that the curves $\zeta_{ij}$ satisfy
%\be \zeta_{ij}''(\lambda)+\frac{\nabla W(\zeta_{ij}(\lambda))}{2}=0 \label{eczeta}.\ee
%However, recently in
%\cite{expstana} Alikakos, Betel\'u and Chen proved that in general, when considering  three potentials  curve solutions to equation \equ{eczeta} that join a given  pair of minima might not exist. 
Hence, in what follows we will work under Hypothesis \ref{geod} and,  in particular,  we assume that for any pair of minima $c_i, c_j$ there is a solution to \equ{eczeta}-\equ{explim}.

%\begin{rem}
%The last assertion in Hypothesis \ref{geod} is not stated in \cite{vecvalps} as a result of the corresponding Lemma, but, as specified by Sternberg, 
%follows from the proof in this paper.
%\end{rem}

\smallskip

As mentioned in the introduction, we want to  relate equation \equ{ecprinc}-\equ{condainf}
with the following equation in the unit ball:
\begin{align}-\Delta u_\epsilon+ \frac{\nabla_vW(u_\epsilon)}{\epsilon^2}&=0 
\hbox{ for } x\in B_1\label{laeq}\\
 u_\epsilon|_{\partial B_1}(x)&=\phi_\epsilon(x).\label{bc}\end{align}
where $\phi_\epsilon$ will be properly defined in the coming sub-section.
This equation motivates us to define the following functional:
 \be \calI_\epsilon(u)=\left\{\begin{array}{cl} \int_{B_1}\epsilon|D u|^2 +\frac{1}{\epsilon}W(u) dy &\hbox{ if }u\in H^1(B_1) \hbox{ and } u|_{\partial B_1}(x)=\phi_\epsilon(x)\\ \infty & \hbox{ otherwise. }\end{array} \right. 
\label{deffeps}\ee
where
 $u:B_1 \to \rr^2$,  $\phi_\epsilon:\partial B_1\to \rr^2.$ 
It is easy to check that weak solutions to  \equ{laeq}-\equ{bc}  can be regarded as
 critical points of \equ{deffeps}.
%Notice that critical points of $F_\epsilon$ are weak solutions to the equation

We are interested in studying the limiting problem as $\epsilon\to 0$. More specifically, 
we expect the limit of the solutions $u_\epsilon$ to \equ{laeq}
will capture the behavior at infinity of the function $v$ which satisfies
 \equ{ecprinc}-\equ{condainf}. In particular, 
we 
want
to show that it is possible to obtain as the limit of the functions $u_\epsilon$ a function $u_0$ that satisfies
%Consider  . Fix $\theta_1,\theta_2\in [0, 2\pi)$, such that condition  is %satisfied for $\alpha_i=\theta_.
% We are going to consider $c_0=c_3$.
%And define the function $u_0$ almost everywhere by
\be u_0(r\cos\theta,r\sin\theta)=c_i \hbox{ for }\theta \in (\theta_{i-1},\theta_i), \label{defu0}\ee 
where $\alpha_i=\theta_i-\theta_{i-1}$ satisfy \equ{condang}. Without loss of generality we are going to assume that $\theta_0=0$ and $\theta_3=2\pi$.

%We also expect that as $\epsilon \to 0$ the solutions $u_\epsilon$  approach a function $u_0$ that satisfies almost everywhere equations \equ{laeq}-\equ{bc} with
% $\epsilon=0$. For this, it would be necessary to have $u_0\in \{c_i\}_{i=1}^3$ almost everywhere.
% and $u_\epsilon|_{\partial B_1}(x)=\phi_0(x)$ almost everywhere.  
 
 In order to study the limit of the functions $u_\epsilon$ 
%thiswe will consider functions $u:B_1\to \rr^2$and 
we define the following limit functional (that we will show corresponds to the $\Gamma$-limit of the functionals $\calI_\epsilon$):

\be \calI_0(u)=\left\{ \begin{array}{cc} \sum_{i,j=1}^3 \Gamma(c_{i}, c_j)
H_1\left(\partial_{B_1}  \Omega_{i}(u)\bigcap \partial_{B_1} \Omega_{i+1}(u)\right)&\hbox{if }g_i(u)\in BV(B_1)\\ +\sum_{i,j=1}^3 \Gamma(c_{i}, c_j)
H_1\left(\left(  \partial \Omega_{j}(u)\bigcap \partial B_1\right)\setminus
 \Phi_{i}\right) &\hbox{ and }u\in \{c_i\}_{i=0}^3 \\ 
&\\
\infty 
&\hbox{ otherwise,} \end{array} \right.\label{deff0}\ee
where $\Omega_i(u)=\{x \in B_1: u(x)=c_i\}$, $\phi_0(x)=\lim_{\epsilon \to 0}\phi_\epsilon (x)$, $\Phi_i=\{x\in \partial B_1: \phi_0(x)=c_i\}$ and $H_1$ is the one dimensional Hausdorff measure.

\subsection{The function $\phi$}\label{constphi}

As described in the introduction, the function $\phi$ should represent the boundary condition at infinity, that is, it should satisfy \equ{condainf}. In particular, we expect the sequence of functions $\phi_\epsilon$
(defined by  $\phi_\epsilon(x)=\phi\left(\frac{x}{\epsilon}\right)$) to  converge to $c_i$ as $\epsilon\to 0$
in the angular sectors of $B_1$ defined by
  $\theta \in (\theta_{i-1},\theta_i)$ (where the angles $\theta_i$ are defined by \equ{condang}-\equ{defthetai}).
 Moreover, we will construct  a function $\phi$  that away from the triple point,  approximates a solution to \equ{ginzlandel}  
(we will make this statement more precise in section \ref{unifconv}). 

More precisely,
let 
 $L_i$ be the half-lines starting at the origin, with direction $ \theta_i$. 
Away from $L_i$, 
 the
function $\phi$ is defined by one of the constants $c_j$
 (that is, one of the 
 minima of $W$). Notice that in fact $c_j$  are solutions to \equ{laeq}. 
Near the half-lines $L_i$, the function $\phi$ will be  equal to an appropriate  solution to \equ{eczeta} (that we denote $\zeta_{ij}$),  evaluated at the distance to $L_i$. These functions  are  approximate solutions in the sense to be discussed in section \ref{unifconv}.
%Hence, when $\epsilon$ approaches 0, the functions $\phi_\epsilon$ will approach $u_0$.

We summarize  the description above with  the following equations:
Consider a smooth function $\eta:\rr^2\to \rr$  such that $\eta(x)\equiv 1$ when $|x|\leq\frac{1}{2}$ 
and $\eta(x)\equiv 0$ for $|x|\geq 1$, the distance
$$d_i(x)=d(x,L_i),$$
and
 a partition of unity $\{\eta_i\}_{i=1}^6$ associated to  the family of intervals $\{\calA_j\}_{j=1}^6$, where 
%$\{(2\pi-\delta,2\pi]\bigcup[0,\delta), (\frac{\delta}{2},\theta_1-\frac{\delta}{2}),(\theta_1-\delta,\theta_1+\delta),
%(\theta_1+\frac{\delta}{2},\theta_2-\frac{\delta}{2}),(\theta_2-\delta,\theta_2+\delta),(\theta_2+\frac{\delta}{2},2\pi-\frac{\delta}{2})\}$,
$$\calA_{2i}=(\theta_i-\delta, \theta_i+\delta),$$
$$\calA_{2i+1}=\left(\theta_i+\frac{\delta}{2},\theta_{i+1}-\frac{\delta}{2}\right).$$
%aca!
 %we assume $\eta_j(x)\geq 0$, $\eta_j (x)\equiv 0$ for $x\not \in \calA_j$ and
%$\sum_j\eta_j(x)=1$ for every $x$.
Now we
define \be \phi(x)=(1-\eta(x))\left(\eta_5(\theta)c_3+\eta_6(\theta) \zeta_{31}(d_0(x))+\sum_{i=1}^2\left( \eta_{2i}(\theta)\zeta_{ii+1}(d_i(x))+\eta_{2i-1}(\theta)c_i \right) \right) \label{deffi}\ee
and \be \phi_\epsilon(x)=\phi\left(\frac{x}{\epsilon}\right).\label{deffieps}\ee
Notice that since $L_i$ is a half-line we have that $d_i\left(\frac{x}{\epsilon}\right)=\frac{d_i(x)}{\epsilon}$.

\begin{rem}\label{discfi}
The functions $\phi_\epsilon$ are not only well defined on the boundary of $B_1$, but also in the interior. Moreover,  under these definitions we have that 
%$$\phi_\epsilon \to  \phi_0 \hbox{ uniformly in }B_1$$ and
$$\phi_0(x):=\lim_{\epsilon\to 0}\phi_\epsilon(x)=u_0(x) \hbox{ a.e.}$$

Furthermore, in Section \ref{unifconv} will be shown that near the boundary (more precisely for $|x|>\epsilon^\alpha$) the function $\phi_\epsilon$ is an ''approximate solution" to the equation \equ{ecprinc} (in the sense that
for every $x$ holds $\left(-\Delta \phi_\epsilon+\frac{\nabla_v W(u)}{2\epsilon^2}\right)(x)\to 0 $ as $\epsilon\to 0$.)  We will prove  that in fact for every $\alpha<1$ holds $\sup_{\epsilon^\alpha<|x|<1}|u_\epsilon-\phi_\epsilon|\to 0$ as $\epsilon \to 0$. Correspondingly, for $v_\epsilon:B_{\frac{1}{\epsilon}}\to \rr^2$ defined by $v_\epsilon(x)=u_\epsilon(\epsilon x)$ holds  
 $\sup_{\epsilon^{\alpha-1}<|x|<\frac{1}{\epsilon}}\left|v_\epsilon-\phi\right|\to 0$ as $\epsilon \to 0$.
%for the rescaled equation). 

On the other hand, it is not expected  that the functions $\phi_\epsilon$ are  good approximations to the solution inside the ball of radius $\epsilon^\alpha$ (or correspondingly, $\phi$ is not a good approximation of $v_\epsilon$ in the ball of radius $\epsilon^{\alpha-1}$). This can be illustrated as follows: 
The choice of the functions $\phi_\epsilon$ in \equ{deffieps} flexible as long as the features described above are preserved (namely, for $|x|>\epsilon^\alpha$ they approach $u_0$ and they are an approximated solution to the equation). 
For example, it is possible to consider $\tilde{\phi}_\epsilon(x)=\phi\left( 
\frac{x}{\epsilon}+\ln(\epsilon) x_0\right)$.  In fact, for every $k\in\nn$ holds $\sup_{|x|>\epsilon^\alpha}|D^k\phi_\epsilon-D^k\tilde{\phi}_\epsilon|\to 0
$ as $\epsilon \to 0$. 
However, for every $\sigma<\epsilon^{1-\alpha}$ we have
$\min{c_i}<\sup_{|x|<\epsilon^{\alpha}}\left|\tilde{\phi}_\epsilon(x)-\tilde{\phi}_\sigma\left(\frac{\sigma x}{\epsilon}\right)\right|$, which contrasts with the second inequality in Theorem \ref{aprox}. In particular, it is clear that $\tilde{\phi}$ cannot not be a good approximation of the solution inside the ball of radius $\epsilon^\alpha$. Similarly, it is not expected that $\phi_\epsilon$ approximates the solution $u_\epsilon$ inside the ball of radius $\epsilon^\alpha$ (or that the corresponding  function $v_\epsilon$ would be
approximated by $\phi$ inside the ball of radius $\epsilon^{\alpha-1}$).

\end{rem}

%We are also going to restate some  lemmas that had been proven before in the literature.
\medskip

\subsection{Techinal Lemmas}
Now we  state some technical lemmas.  
The first one was originally proved  in \cite{tesis}:

\begin{lem} \label{cota} 
Let $u_\epsilon(x)\in \calC^2$ 
satisfy \equ{laeq}-\equ{bc}, where $W:\rr^2\to \rr$ is a proper function in
$ \calC^2$ 
 bounded below, with a finite number of critical points (that we label as $\{c_i\}_{i=1}^m$),   
%such that $W(v)\to \infty$ as  $|v|\to \infty$
 and such that the Hessian
of $W(u)$ is positive semidefinite for $|u|\geq K$ for some real number $K$. Suppose that the functions $\phi_\epsilon$ are uniformly bounded. Then 
there is a constant $C$ depending only on  uniform bounds over $\phi_\epsilon$ and  $W$, but not on $\epsilon$, such that
$$\sup|u_\epsilon|\leq C. $$
%where $C$ only depends on uniform bounds over $\phi_\epsilon$ and  $W$.
\end{lem}
\begin{proof}

Consider $\omega_\epsilon(x)=W(u_\epsilon)(x)$; then 
\begin{align*} -\Delta \omega_\epsilon &=
-\sum_i (\nabla_v W(u_\epsilon) \cdot (u_\epsilon)_{x_i})_{x_i}\\
&= - 
(W''(u_\epsilon) D u_\epsilon)\cdot D u_\epsilon-\nabla_v W(u_\epsilon) \cdot \Delta u_\epsilon,
\end{align*}
where $W''$ denotes the Hessian matrix of $W$ and the dot product between two $2\times 2$ matrices is the standard dot product in $\rr^4$.
Since $u_\epsilon$ satisfies \equ{laeq}, this becomes
\be -\Delta \omega_\epsilon +\frac{|W'(u_\epsilon)|^2}{2\epsilon ^2} 
+(W''(u_\epsilon)D u)\cdot D u_\epsilon=0. \label{eqv0} \ee

If the maximum of $\omega_\epsilon$ is attained at the boundary, then it is bounded by  the maximum of $W(\phi_\epsilon(x))$.

Suppose that $\omega_\epsilon$ has an interior maximum at $x_0$ and 
$|u_\epsilon(x_0)|\geq K$. 
Since $x_0$ is a maximum for $\omega_\epsilon$, it holds that
 $\Delta \omega_\epsilon(x_0)\leq 0$. We also have by hypothesis that
$W''(u)$ is positive semidefinite for $|u|\geq K$, hence
$$-\Delta \omega_\epsilon +\frac{|D _u W(u_\epsilon)|^2}{\epsilon ^2} +(W''(u_\epsilon)D u_\epsilon)\cdot
D u_\epsilon \geq 0.$$
The inequality is strict (which contradicts \equ{eqv0}) 
unless 
$$\frac{|D _u W(u_\epsilon)|^2}{\epsilon ^2} =(W''(u_\epsilon)D u_\epsilon)\cdot
D u_\epsilon =0.$$
If $\nabla_v W (u_\epsilon(x_0))=0$, we would have  $u_\epsilon(x_0) =c_i$ for some $i$ and this implies (since the maximum is attained at this point) that 
$W(u_\epsilon(x,t))\leq W(c_i)$. Hence we have $\omega_\epsilon\leq \max\{\sup_{|u|\leq K} W(u_\epsilon), W(\phi_\epsilon), \max_{i=1\ldots m } W(c_i)\}$.

Since $W$ is a proper function, we conclude the result of the Lemma. \end{proof}

\medskip

We will also  use Lemma A.1 and Lemma A.2   in \cite{asyforfb0}.  We restate 
them here without proof:

\begin{lem}\label{lema2}[Lemma A.1 in  \cite{asyforfb0}]
Assume that $u$ satisfies 
$$-\Delta u=f \hbox{ on }\Omega \subset \rr^n$$
%$$u=0 \hbox{ on } \partial \Omega $$
%where $\Omega$ is a smooth bounded domain. 
Then 
\be |D u(x)|^2 \leq C \left(\| f \|_{L^\infty(\Omega)}\|u\|_{L^\infty(\Omega)}+\frac{1}{dist^2(x,\partial \Omega)}\|u\|_{L^\infty(\Omega)}^2\right)\quad \forall x\in \Omega, \ee
where $C$ is a constant depending only on $n$.
\end{lem}

%Assume that $u$ satisfies 
%$$-\Delta u=f \hbox{ on }\Omega \subset \rr^n$$
%$$u=0 \hbox{ on } \partial \Omega $$
%where $\Omega$ is a smooth bounded domain. Then it holds
%\be \|D u\|^2_{L^\infty(\Omega)} \leq C \| f \|_{L^\infty(\Omega)}\|u\|_{L^\infty(\Omega)}\ee
%where $C$ is a constant depending only on $\Omega$.

\begin{lem}\label{lema1} [Lemma A.2 in  \cite{asyforfb0}]

Assume that $u$ satisfies 
$$-\Delta u=f \hbox{ on }\Omega \subset \rr^n$$
$$u=0 \hbox{ on } \partial \Omega $$
where $\Omega$ is a smooth bounded domain. Then it holds
\be \|D u\|^2_{L^\infty(\Omega)} \leq C \| f \|_{L^\infty(\Omega)}\|u\|_{L^\infty(\Omega)}\ee
where $C$ is a constant depending only on $\Omega$.

\end{lem}
\section{Convergence in $L^1$}\label{convl1}
 
%Let us fist define our candidate to limit solution:

In this section we show that  solutions $u_\epsilon$ to equation \equ{laeq}-\equ{bc} converge in $L^1$. More precisely, we prove
the following result 

\begin{prop} \label{stern} Let  $u_0$ be defined by \equ{defu0}. Consider $\calI_\epsilon$ and $\calI_0$ defined by \equ{deffeps} and \equ{deff0} respectively.
For $\phi_\epsilon$ defined by  \equ{deffi}-\equ{deffieps} there exists a sequence of minimizers $u_\epsilon$
 of $\calI_\epsilon$,  such that $\calI_\epsilon(u_\epsilon)\to \calI_0(u_0)$ and $u_\epsilon \to u_0$ in $L^1$. 
\end{prop}

As stated in \cite{locminps},  when considering the Neumman boundary condition problem, Proposition \ref{stern} follows from  results in  \cite{minintsb}, \cite{locminrk} and \cite{locminps}.
In what follows we are going to state these  results and point out the necessary modifications in our setting.

%An analogous the result prove by P. Sternberg in \cite{vecvalps}. We will sketch the proof following this paper. We will restate Theorem 1 used by Sternberg, that was original proven in \cite{locminrk}.

\begin{teo}\label{minimizer}(\cite{locminps})
Let $u_0$  be defined by \equ{defu0} and  $C_i=\{x\in\Omega:u_0(x)=c_i\}$. 
Consider a domain $\Omega$ and partition $(E,F,G)$ of $\Omega$. Define
\begin{align*} \calF(E,F,G)=\Gamma(c_1,c_2) &H_1(\partial_\Omega E\bigcap \partial_\Omega G) +\Gamma(c_1,c_3)H_1(\partial_\Omega E\bigcap \partial_\Omega F)\\&+\Gamma(c_3,c_2)H_1(\partial_\Omega F\bigcap \partial_\Omega G).\end{align*}
Then
the partition formed by $C_1, C_2$ and $C_3$ is an isolated local minimizer 
of $\calF$, that is
 \be \calF(C_1,C_2,C_3)=\min  \calF(E,F,G)\ee
%\be \min\{\Gamma(c_1,c_2)H_1(\partial_\Omega E\bigcup \partial_\Omega G) +\Gamma(c_1,c_3)H_1(\partial_\Omega E\bigcup \partial_\Omega F)+\Gamma(c_3,c_2)H_1(\partial_\Omega F\bigcup \partial_\Omega G)\}\ee
where the minimum is taken over all the partitions $(E,F,G)$ of $\Omega$ satisfying the condition
\be|C_1\Delta E|+|C_2\Delta F|+|C_3\Delta G|\leq\delta,\ee
where $\delta$ is some small positive number.
\end{teo}
\medskip
\begin{rem}
The proof of Lemma 3.1 in \cite{locminps} implies that this $\delta$ can be uniformly chosen for balls of all radii.
\end{rem}
\medskip
\begin{teo}(Theorem 2.5 in \cite{minintsb}) \label{teo1b}
Let 
\be \tilde{\calI}_{\epsilon,\Omega}(u)
=\left\{\begin{array}{cl} \int_{\Omega}\epsilon|D u|^2 +\frac{1}{\epsilon}W(u) dy &\hbox{ if }u\in H^1(\Omega) \hbox{ and } \int_{\Omega}u(x)dx=m
\\ \infty & \hbox{ otherwise. }\end{array} \right. \label{defitildeps}\ee
and 
\be \tilde{\calI}_{0,\Omega}(u)=\left\{ \begin{array}{cl} \sum_{i,j=1}^3 \Gamma(c_{i}, c_j)
H_1\left(\partial_{B_1}  \Omega_{i}(u)\bigcap \partial_{B_1} \Omega_{j}(u)\right)&\hbox{ if } g_i(u)\in BV(\Omega) \hbox{ for }i\in\{1,2,3\}, \\
& W(u(x))=0\hbox { a.e.} \hbox{ and }\int_\Omega u(x)dx=m
\\ \infty 
&\hbox{ otherwise} \end{array} \right. \label{defitild0}\ee

It holds for every $\epsilon_h\to 0$ that
\begin{itemize} 

\item For every $u_{\epsilon_h}\to u$ in $L^1(\Omega)$ we have that $\tilde{\calI}_0(u)\leq \liminf_{h\to\infty}\tilde{\calI}_{\epsilon_h}(u_{\epsilon_h})$
\item There is  $u_{\epsilon_h}\to u$ in $L^1(\Omega)$ such that $\tilde{\calI}_0(u)\geq \limsup_{h\to\infty}\tilde{\calI}_{\epsilon_h}(u_{\epsilon_h})$
\end{itemize}
\end{teo}
\medskip

\begin{prop}\label{teo2b}( Proposition 2.2 in \cite{minintsb}) The function $g_i$ is locally Lipschitz-continuous. Moreover, if $u\in H^1(\Omega)\bigcup L^\infty(\Omega)$, then $g_i(u)\in W^{1,1}(\Omega)$ and the following inequality holds:
\be \int_\Omega |D(g_i(u))|dx\leq\int_\Omega \sqrt{W(u)}|Du|dx. \ee
\end{prop}
\medskip
\begin{rem}\label{remb}
Following the proof of Theorem \ref{teo1b} in \cite{minintsb} it is easy to see that the restriction $\int_\Omega u(x)dx=m$, imposed by Baldo in his work, can be removed  from Theorem \ref{teo1b} without modifying 
the proof.
\end{rem}
\medskip
\begin{teo}\label{teostern} \cite{locminrk}
Suppose that a sequence of functionals $\{\calI_\epsilon\}$ and a functional $\calI_0$ satisfying the following conditions:
\begin{enumerate}
\item \label{lowsemcont} if $w_\epsilon \to w_0$ in $L^1(\Omega)$ as $\epsilon \to 0$, then $\liminf \calI_\epsilon(w_\epsilon)\geq \calI_0(w_0)$;
\item\label{minatt}  for any $w_0 \in L^1(\Omega)$ there is a family $\{\rho_\epsilon\}_{\epsilon>0}$ with $\rho_\epsilon \to w_0$ in $L^1(\Omega)$ and $\calI_\epsilon (\rho_\epsilon)\to \calI_0(w_0)$;
\item \label{comp} any family $\{w_\epsilon\}_{\epsilon>0}$ such that $\calI_\epsilon(w_\epsilon)\leq C<\infty$ for all $\epsilon>0$ is compact in $L^1(\Omega)$;
\item \label{minex} there exits an isolated $L^1$-local minimizer $u_0$ of $\calI_0$; that is, $\calI_0(u_0)<\calI_0(w)$ whenever $0<\|u_0-w\|_{L^1(\Omega)}\leq \delta$ for some $\delta >0$.

\end{enumerate}

Then there exits an $\epsilon_0>0$ and a family $\{u_\epsilon\}$ for $\epsilon<\epsilon_0$ such that $u_\epsilon$ is an $L^1$-local minimizer of $\calI_\epsilon$ and $u_\epsilon \to u_0$ in $L^1(\Omega)$
\end{teo}

Theorem \ref{teo1b}  establish conditions
\ref{lowsemcont} and \ref{minatt} of Theorem \ref{teostern} for $\tilde{\calI}_{\epsilon,\Omega}$  (defined by \equ{defitildeps}) and $\tilde{\calI}_0$(defined by \equ{defitild0}). 
Theorem \ref{minimizer} establishes that $u_0$ is a local minimizer for $\tilde{\calI}_{0,\Omega}$ (condition \ref{minex} of Theorem \ref{teostern}). 
We need to show that these theorems imply that the conditions  of Theorem    \ref{teostern}  also hold for $\calI_\epsilon$ and $\calI_0$  (defined by \equ{deffeps} and \equ{deff0}, respectively). In addition, we need to prove  that condition
\ref{comp} holds for these functionals.

\begin{lem}\label{locminf0}
Theorem \ref{minimizer} implies that $u_0$ is a local minimizer for $\calI_0$.

\end{lem}
\begin{proof}

Let $C_i=\{x \in B_1 :u_0(x)=c_i\}$ and for any $w$ let $\Omega_i(w)=\{x\in B_1:w(x)=c_i\}$.
Consider $\delta$ for $B_1$ as is Theorem \ref{minimizer}. We are going to show by contradiction
that for every $w$ such that $w(x)\in\{c_i\}_{i=1}^3$ almost everywhere and
$$|C_1\Delta \Omega_1(w)|+|C_2\Delta \Omega_2(w)|+|C_3\Delta \Omega_3(w)|\leq
\delta $$
holds that
$$\calI_0(u_0)\leq \calI_0(w).$$

Suppose that there is a $w$ such that 
\be |C_1\Delta \Omega_1(w)|+|C_2\Delta \Omega_2(w)|+|C_3\Delta \Omega_3(w)|\leq
\delta  \label{difsim}\ee
and
\be \calI_0(u_0)> \calI_0(w). \label{f0mal}\ee

Consider $\sigma>0$ and $B_{1+\sigma}$. 
Define $$\calI_\epsilon^\sigma(u)=\tilde{\calI}_{\epsilon,B_{1+\sigma}}(u).$$
Notice first that $u_0$ (given by \equ{defu0})  is well defined for every $x \in \rr^2$. In particular is well defined for every $x \in B_{1+\sigma}$ for any $\sigma >0$. Hence, we can
define
\be w^\sigma(x)=\left\{\begin{array}{cc}w(x)&\hbox{ if } x\in \bar{B_1}\\
u_0(x)&\hbox{ if } x \in B_{1+\sigma}\setminus B_1
\end{array} \right..\label{defvsigma}\ee

Let
$$\tilde{C}_i=\{x\in B_{1+\sigma}:u_0(x)=c_i\}$$
$$\tilde{\Omega}_i(w)=\{x\in B_{1+\sigma}:w^\sigma(x)=c_i\}.$$
Using definition \equ{defvsigma} and equation \equ{difsim} we also have
\be |\tilde{C}_1\Delta \tilde{\Omega}_1(w)|+|\tilde{C}_2\Delta \tilde{\Omega}_2(w)|+|\tilde{C}_3\Delta \tilde{\Omega}_3(w)|\leq
\delta. \label{difsim2}\ee

Notice that every subset on the boundary where $w$  does not agree with $u_0$  becomes an interior boundary term for $w^\sigma$ in $B_{1+\sigma}$. By the definition of definition $\calI_0^\sigma$ we have that
$$\calI^\sigma_0(w^\sigma)=\calI_0(w^\sigma)+\sigma \sum_{i,j=1}^3\Gamma(c_i,c_j)$$
and 
$$\calI^\sigma_0(u_0)=\calI_0(u_0)+\sigma \sum_{i,j=1}^3\Gamma(c_i,c_j).$$
 
Inequality \equ{f0mal} implies that
\be \calI^\sigma_0(w^\sigma)<\calI^\sigma_0(u_0), \label{ftildemal}\ee

which together with \equ{difsim2} contradicts the local minimality of $u_0$ given by Theorem \ref{minimizer}. \end{proof}
\medskip

\begin{proof}[Proof of Proposition \ref{stern}]

In what follows, we are going to show that  Theorem \ref{teo1b} and Proposition \ref{teo2b} imply conditions \ref{lowsemcont} and \ref{minatt}  of Theorem \ref{teostern}  for the functionals defined by \equ{deffeps} and  \equ{deff0}. 

Recall that $\phi_\epsilon$ is given by \equ{deffieps}, $\phi_0=\lim_{\epsilon\to 0}\phi_\epsilon$ and
$\phi_0=u_0$ a.e.
%We are also going to assume that the function $\phi_0(x)\in\{c_i\}_{i=1}^3$. 
%\begin{lem}\label{condw}
%For the functionals $F_\epsilon$ and $F_0$ defined by \equ{deffeps} -\equ{deff0} and the function $\phi_0=u_0|_{\partial B_1}$  (\equ{defu0}), conditions \ref{lowsemcont}, \ref{minatt} and \ref{comp} hold.
%\end{lem}
\medskip

\begin{proof}[Proof of condition \ref{lowsemcont}]

Let \be w_\epsilon \to w_0 \hbox{ in } L^1. \label{conl1vls}\ee

As in the proof of Lemma \ref{locminf0}, consider $\sigma>0$ and define
\be \calI_\epsilon^\sigma(u)=\tilde{\calI}_{\epsilon,B_{1+\sigma}}(u),\label{deffsigls}\ee
\be w^\sigma_\epsilon(x)=\left\{\begin{array}{cc}w_\epsilon(x)&\hbox{ if } x\in \bar{B_1}\\
\phi_\epsilon(x)&\hbox{ if } x \in B_{1+\sigma}\setminus B_1
\end{array} \right..\label{defvepssigmals}\ee
and
\be w^\sigma_0(x)=\left\{\begin{array}{cc}w_0(x)&\hbox{ if } x\in \bar{B_1}\\
\phi_0(x)&\hbox{ if } x \in B_{1+\sigma}\setminus B_1
\end{array} \right..\label{defv0sigmals}\ee

Notice that again the boundary portions of $w_0$ that do not agree with $\phi_0$ become interior boundaries of $w_\sigma^0$. Hence, as before, if $\calI_0^\sigma(w_0)\ne \infty$ we have that 
\be \calI_0^\sigma(w_0)=\calI_0(w_0)+\sigma \sum_{i,j=1}^3\Gamma(c_i,c_j). \label{f0v0ls}\ee

Using  \equ{conl1vls} and definitions \equ{defvepssigmals} and 
\equ{defv0sigmals} we have that 
$$ w^\sigma_\epsilon \to w^\sigma_0 \hbox{ in } L^1.$$ Theorem \ref{teo1b} and Remark \ref{remb} imply that 
\be \calI_0^\sigma(w_0^\sigma)\leq \liminf_{\epsilon\to 0}\calI_\epsilon^\sigma(w_\epsilon^\sigma). \label{lsfsigma}\ee

We can explicitly compute that  
\be \tilde{\calI}_{\epsilon,B_{1+\sigma}\setminus B_1}(\phi_\epsilon)\to \sigma \sum_{i,j=1}^3\Gamma(c_i,c_j).\label{restols}\ee

It is also easy to check that
\be \calI_\epsilon^\sigma(w_\epsilon)=\calI_\epsilon(w_\epsilon)+ \tilde{\calI}_{\epsilon,B_{1+\sigma}\setminus B_1}(\phi_\epsilon).
\label{fepsvls}\ee
Equations \equ{restols} and \equ{fepsvls} imply that
%we have that 
$$\calI_\epsilon^\sigma(w_\epsilon)\to \infty \hbox{ if and only if }
\calI_\epsilon (w_\epsilon)\to \infty.$$

We can assume that
$\liminf_{\epsilon \to 0}\calI_\epsilon(w_\epsilon)<\infty$ (otherwise the result is trivial). Equations
\equ{f0v0ls}, \equ{lsfsigma}, \equ{fepsvls} and \equ{restols} imply that
\begin{align*}\calI_0(w_0)+\sigma \sum_{i,j=1}^3\Gamma(c_i,c_j) =&\calI_0^\sigma(w_0)\\
\leq &\liminf_{\epsilon\to 0}\calI_\epsilon^\sigma(w_\epsilon^\sigma)\\
=& \liminf_{\epsilon\to 0} \calI_\epsilon(w_\epsilon)+\sigma \sum_{i,j=1}^3\Gamma(c_i,c_j).\end{align*}
This implies
$$\calI_0(w_0)
\leq
\liminf_{\epsilon\to 0} \calI_\epsilon(w_\epsilon),$$
which proves the result. \end{proof}

%\medskip+\sigma \sum_{i,j=1}^3\Gamma(c_i,c_j)\\
%\leq &\liminf_{\epsilon\to 0}F_\epsilon^\sigma(v_\epsilon^\sigma)\\
%=& 
\medskip

\begin{proof}[ Proof of condition \ref{minatt}]
                                                                              
The proof of condition \ref{minatt} follows directly from the proof in
\cite{minintsb} of the equivalent statement. Hence, we are going to follow Baldo's proof, use some of his  constructions   and point out the necessary modifications
in our setting. For more details, we refer the reader to \cite{minintsb}.

As in the proof of condition \ref{lowsemcont}, let $\calI^\sigma_\epsilon$ be defined by 
\equ{deffsigls}, that is
$$ \calI_\epsilon^\sigma(u)=\tilde{\calI}_{\epsilon,B_{1+\sigma}}(u).$$

Consider $w_0\in \{c_i\}_{i=1}^3$, such that $\calI_0(w_0)<\infty$ 
(otherwise the result is trivial).                                                                                
As before, we extend the domain to $B_{1+\sigma}$, for some 
$\sigma>0$, and we extend $w_0$ by $\phi_0$ outside the unit ball. We label this extension as $w_0^\sigma$.

Let $\rho^\sigma_\epsilon$ be the sequence of functions given by Theorem \ref{teo1b} that satisfy $\rho_\epsilon^\sigma \to w_0^\sigma$  in $L^1$ and $\calI^\sigma_\epsilon(\rho_\epsilon)\to \calI_0^\sigma(w_0^\sigma)$.

We can write $w_0=\sum_{i=1}^3 c_i 1_{\Omega_i}$. The functions $\rho_\epsilon^\sigma$ constructed by Baldo in \cite{minintsb} are uniformly bounded functions, that  $\epsilon$- near the boundaries $\partial \Omega_i\bigcap \partial \Omega_j \bigcap 
B_{1+\sigma}$ are equal to the geodesic $\zeta_{ij}$. In the interior of 
$\Omega_i$, $\rho^\sigma_\epsilon$ approaches $c_i$ uniformly. In particular, we have that $\rho_\epsilon\to w_0$ almost everywhere and it is uniformly bounded. By dominated convergence theorem we have that the restriction of $\rho^\sigma_\epsilon$ to $B_1$, that we will label as $\rho_\epsilon$, converges  to $w_0$ in the $L^1$ norm. 

As in the proof of \ref{lowsemcont}, we have
\be \calI_0^\sigma(w_0^\sigma)=\calI_0(w_0)+\sigma \sum_{i,j=1}^3\Gamma(c_i,c_j).\label{ig1minatt}\ee

By the definitions of $\calI_\epsilon^\sigma$, $\calI_\epsilon$ , $\rho^\sigma_\epsilon$ and $\rho_\epsilon$,  for every $\sigma>0$ holds that
\be \calI_\epsilon^\sigma(\rho^\sigma_\epsilon)\geq \calI_\epsilon(\rho_\epsilon).\label{des1minatt}\ee

Combining \equ{ig1minatt}, \equ{des1minatt} and Theorem \ref{teo1b} we have 
\begin{align*} \calI_0(w_0)+\sigma \sum_{i,j=1}^3\Gamma(c_i,c_j)&= \calI_0^\sigma(w_0^\sigma)\\ &=\lim_{\epsilon\to 0}\calI_\epsilon^\sigma(\rho_\epsilon^\sigma)\\&\geq 
\lim_{\epsilon\to 0}\calI_\epsilon(\rho_\epsilon).\end{align*}

Taking $\sigma \to 0$ follows that 
$$ \calI_0(w_0)\geq 
\lim_{\epsilon\to 0}\calI_\epsilon(\rho_\epsilon).$$
Combining this equation and Condition \ref{lowsemcont} (that we proved above)
we conclude that
$$ \calI_0(w_0)=
\lim_{\epsilon\to 0}\calI_\epsilon(\rho_\epsilon),$$
which finishes the proof. \end{proof}

\medskip

\begin{proof}[ Proof of condition \ref{comp}]

We will follow the proof in \cite{vecvalps}.
Suppose that $\calI_\epsilon(w_\epsilon)\leq C<\infty$ for some family $\{w_\epsilon\}_{\epsilon>0}$.

Define $$G_\epsilon(x)=g_1(w_\epsilon(x)).$$
Proposition \ref{teo2b} implies that 
\begin{align*}\int_{B_1}|DG_\epsilon(x)|dx\leq&\int_{B_1}\sqrt{W(w_\epsilon)}|Dw_\epsilon|dx \\
\leq & \epsilon\int_{B_1}|Dw_\epsilon|^2 dx+\frac{1}{\epsilon}\int_{B_1}W(w_\epsilon)dx\\
\leq& C.
\end{align*}

 Hypothesis (\ref{cond2W}) of Theorem \ref{teoprinc} implies that $w_\epsilon$ are uniformly bounded in $L^p(B_1)$ for some $p$.
Hence, $G_\epsilon$ are uniformly bounded in $L^1(B_1)$ and
$$\|G_\epsilon\|_{BV(B_1)} \leq C.$$

Since bounded sequences in BV are compact in $L^1$(\cite{minsureg}), there is a subsequence $G_\epsilon$ convergent to $G_0$ in $L^1$. This function $G_0$ takes the form
$$G_0(x)=\left\{\begin{array}{cc}0 &\hbox{ if }x\in C_1 \\
g_1(c_2)& \hbox{ if }x\in C_2 \\g_1(c_3)& \hbox{ if }x\in C_3. 
\end{array}\right.$$
 
Since $c_1$ is the only value $x$ such that $g_1(x)=0$ and $g_1$ is continuous, we have that there is a subsequence $\{w_{\epsilon_j}\}$ that converges in measure to $c_1$ on $C_1$. The uniform bounds in $L^p$ (provided by hypothesis \equ{cond2W}) imply that $\{w_{\epsilon_j}\}$  converge on $C_1$ also in the $L^1$ norm. The proof can be finished  by repeating the same argument for $g_2$ and $g_3$.\end{proof}

Since Lemma \ref{locminf0} implies condition \equ{minex} of Theorem \ref{teostern}, using that theorem we conclude the result of Proposition \ref{stern}. \end{proof}

From Theorem \ref{teostern} we conclude the following corollary:
\begin{cor}\label{pwconv}
Let $u_0$ be defined as in Theorem \ref{teostern}. Then
there is a subsequence of the family $\{u_\epsilon\}$ that converges point-wise almost everywhere to $u_0$.
\end{cor}

\section{Uniform Convergence} \label{unifconv}

In this section we focus on improving the convergence bounds proved in the previous section. Namely,  we prove

\begin{teo}\label{aprox}

  Fix $0<\alpha< 1$. 
Let $0<\sigma\leq \epsilon^{1-\alpha}$ then for every $m>0$ there is a constant $C$ (that might depend on $\alpha$ and $m$) such that

\begin{itemize}
\item $$\sup_{|x|\geq \epsilon^{\alpha}}|u_\epsilon-\phi_\epsilon|\leq C\epsilon^m.$$

\item $$\sup_{|x|\leq \frac{\epsilon^{\alpha}}{2}}\left|u_\epsilon(x)-
u_\sigma \left(\frac{\sigma x}{\epsilon}\right)\right|\leq C\epsilon^m
. $$

\end{itemize}
\end{teo}

\smallskip

There are two main ingredients in the proof of this theorem. The first is  the construction of a function $U_{\vec{q}}$ that satisfies $U_{\vec{q}}(x)=\phi_\epsilon(x)$ for $x\in B_1\setminus B_{\epsilon^\alpha}$,  $U_{\vec{q}}(x)=u_\epsilon\left(\frac{\sigma x}{\epsilon}\right)$ for $x\in B_{\frac{\epsilon^\alpha}{2}}$ and 
$\left|-\Delta U_{\vec{q}}+\frac{\nabla_v W(U_{\vec{q}})}{2\epsilon^2}\right|(x)\to 0$ point-wise; the second one  is 
Theorem \ref{cotaprinc0}. 
The idea is the following: We consider $U_{\vec{q}}$ as the initial condition for the parabolic equation \equ{ginzlandpar} in the unit ball. Since $U_{\vec{q}}$ is almost a solution to this equation, we expect that the actual solution to \equ{ginzlandpar} will stay close $U_{\vec{q}}$. This assertion it is ensured by Theorem \ref{cotaprinc0}.  However, in order to apply that theorem is necessary to consider solutions to an equation with $0$ boundary condition. For this reason, instead of considering equation \equ{ginzlandpar}  we take \equ{ecpar'}-\equ{cbpar'}-\equ{cipar'} (which correspond to subtract  the function $U_{\vec{q}}$ from the solution to \equ{ginzlandpar} ).
We finally conclude Theorem \ref{aprox} by observing that 
 our solution to \equ{ginzlandpar} converge to $u_\epsilon$ as $t\to \infty$.

 %\begin{lem} \label{aproxu<eps}
 
%Let $M$ be the constant in Lemma \ref{aprox>eps}.
%and  $|x|\geq \epsilon^\alpha$ there  are constants $K$ and $c$, independent of $\alpha, \epsilon$ and $\alpha$ such that 
%Let $0 < \alpha< 1$ and $\sigma \leq \epsilon$.
%   Then there is a constant $C$, independent of $\epsilon$ and $\sigma$ such 
%that
%for every $|x| \leq \epsilon^\alpha$ and  $0\leq t \leq \frac{M\epsilon^2}{3}$  holds

%\leq 
%C \left(\beta \epsilon^{2-2\alpha} +\epsilon^{1-\alpha}+\beta
%\int_{B_1}|u_\epsilon(y)-\phi_{\sqrt{\epsilon}}(y)|dy \right.$$
%$$\left. \hspace{5cm}+\beta \int_{B_1}|u_\sigma(y)-\phi_\sigma(y)|dy \right).$$
 %\end{lem}
We would also like to remark that the minimizing property of solutions $u_\epsilon$ will not be used in this section. In fact, the construction presented here would work for any type of critical point of the functional $\calI_\epsilon$ with the appropriate boundary values. 
However, the minimizing property will be used again in section \ref{pf} in order to show the minimizing statement of Theorem \ref{teoprinc}.

%\end{rem}

Now we proceed with the construction of the function $U_{\vec{q}}$. Since this functions depends also from other parameters besides $\epsilon$ (such as $\alpha$ above and $\sigma$, 
which will be shortly introduced) the subindex $\vec{q}$ stands for
$\vec{q}=(\epsilon,\sigma, \alpha)$.

Let
$$v_\epsilon(x)=u_\epsilon(\epsilon x) \hbox{ and, }$$
$$u_\sigma^\epsilon (x)=u_\sigma \left(\frac{\sigma x}{\epsilon}\right).$$
Consider a positive function $\eta:\rr \to \rr$ such that $\eta(x)=0$ for $|x|\leq \frac{1}{2}$ and  $\eta(x)=1$ for $|x|\geq 1$. Fix $0<\alpha <1$ and 
\be E=2 \epsilon^\alpha-\epsilon^{2m+4-\alpha}.\label{defE}\ee
%>1\geq \alpha$  then define
Then define for $y\in \rr^2$
 the function
$$\eta_{\alpha}(y) =\eta\left(\frac{\epsilon}{2E}|y|+1- \frac{\epsilon^{\alpha}}{2E}\right).$$
Notice that the function $\eta_{\alpha}(y) $ satisfies  $\eta_{\alpha}(y)=0$ for $|y|\leq \epsilon^{\alpha-1}-\frac{E}{\epsilon}$ and  $\eta_{\alpha}(x)=1$ for $|y|\geq \epsilon^{\alpha-1}$. Moreover, defining
  $$\eta_{\alpha}^\epsilon(y)=\eta_{\alpha}\left(\frac{y}{\epsilon}\right) $$ it satisfies  $\eta^\epsilon_{\alpha}\left(y\right) =0$ for $|y|\leq \epsilon^{\alpha}-E$ (where $E$ is defined by \equ{defE}) and  $\eta^\epsilon_{\alpha}\left(y\right) =1$ for $|y|\geq \epsilon^{\alpha}$.

We will denote by $\calH_{\Omega}$ the heat Kernel in $\Omega\subset \rr^2$.
A more detailed description and some properties of the Heat Kernel can be found in the Appendix.

Let \be \calQ=\{( \epsilon,\sigma, \alpha)\in(0,1]\times(0,1]\times [0,1]: \quad \sigma\leq \epsilon^{1-\alpha}\}. \label{calq}\ee
Define for  ${\vec{q}}=
(\epsilon,\sigma, \alpha)\in \calQ$  the function
  $$V_{\vec{q}}(y)=\eta_{\alpha}(y) \phi (y)+(1-\eta_{\alpha}(y)) v_\sigma (y).$$ 
  Now we take 
    $$U_{\vec{q}}(y)= V_{\vec{q}}
    \left(\frac{y}{\epsilon}\right).$$

    Let us denote by $\calC_{S}$ the set of continuous functions from $S$ to $\rr^2$.
For $\vec{q}$ as above consider the function $F_{\vec{q}}:\calC_{B_1\times[0,T]} \times \calC_{B_1}\to \calC_{B_1\times[0,T]}$ 
%(where $\calC_{B_1\times[0,T]}$ is the set of uniformly bounded continuous functions) 
defined by
 $$F_ {\vec{q}}(h,\psi)(x,t)=\int_0^t \int_{B_{\frac{1}{\epsilon}}}\calH_{B_{\frac{1}{\epsilon}}}(x,y,t-s)\left(-\nabla_v W(h+V_{\vec{q}})(y,s)+\Delta V_{\vec{q}}\right)dyds$$
$$+ \int_{B_{\frac{1}{\epsilon}}}\calH_{B_{\frac{1}{\epsilon}}}(x,y,t)\psi (y)dy.$$
 
Notice that, for a given $\psi$,  Duhamel's formula implies that, if there is a  fixed point $h_{\vec{q},\psi}$  of $F_ {\vec{q}}(\cdot,\psi)$, it would satisfy
\begin{align}
\frac{d h_{\vec{q},\psi}}{d t}-\Delta h_{\vec{q},\psi}
+\frac{\nabla_v W(h_{\vec{q},\psi}+V_{\vec{q}})}{2}&=\Delta V_{\vec{q}} \hbox{ in }B_{\frac{1}{\epsilon}}\label{ecpar}\\
h(x,t)&=0 \hbox{ on }\partial B_{\frac{1}{\epsilon}}\label{cbpar}\\
h(x,0)&=\psi(x).\label{cipar}
\end{align}
%\begin{rem}: Notice that in particular the function $u_\epsilon(\epsilon x)$ is a solution to this equation when $\psi(x)=u_\epsilon(\epsilon x)$.
%For example, when $\psi(x)=v_\epsilon(x)-V_{\vec{q}}(x)$ we have $h=v_\epsilon(x)-V_{\vec{q}}(x)$.

The next lemma shows the existence of such a fixed point:

\begin{lem}\label{existprobpar} Fix a uniformly bounded continuous function $\psi_\epsilon$ and ${\vec{q}}\in \calQ$, where 
%${\vec{q}}=(\epsilon,\sigma, \alpha)$ and 
$\calQ$ is defined by \equ{calq}. The function $F_ {\vec{q}}(\cdot,\psi):\calC_{B_1\times[0,T]}\to \calC_{B_1\times[0,T]}$ has a unique fixed point that we label $h_{{\vec{q}},\psi}$. Moreover, for $K>0$ and functions $w_{\vec{q}}$ satisfying $|w_{\vec{q}}|\leq K$, there are  constants $M$ and $\beta$  (that might depend on $K$),  such that for every $T\geq 0$ holds 
\begin{align}\sup_{B_{\frac{1}{\epsilon}}\times \left[T, T+\frac{2 \beta}{M}\right]}|w_{\vec{q}}-h_{{\vec{q}},\psi}|
 \leq \frac{1}{1-\beta} &\left(2\sup_{B_{\frac{1}{\epsilon}}\times \left[T, T+\frac{2 \beta}{M}\right]} |F_{{\vec{q}}}(w_{\vec{q}},\psi)-w_{\vec{q}}|\right.
 \notag \\
& \quad \left. +
\sup_{x\in B_{\frac{1}{\epsilon}}}|w_{\vec{q}}-h_{{\vec{q}},\psi}|(x, T)\right).\label{ineq0}\end{align}

\end{lem}
%\noindent{\bf Proof:}
We postpone the proof of this Lemma to the Appendix. 

\medskip

% By noticing that 
%$v_\epsilon-V_{\vec{q}}$ is a solution to \equ{ecpar}-\equ{cbpar}-\equ{cipar}
% we conclude
%\begin{cor} \label{corueps}
%For every fixed ${\vec{q}}\in \calQ$ such that ${\vec{q}}=(\epsilon,\sigma,\beta)$ it holds $v_\epsilon-v_{{\vec{q}}} \equiv v_{{\vec{q}},v_\epsilon-v_{{\vec{q}}} }. $ 
%is a fixed point of the function $F_ q(\cdot,v_\epsilon):\calC \to \calC$.
%\end{cor}
\medskip
From Lemma \ref{existprobpar} we can prove the following theorem (which provides one of the essential tools in the proof of Theorem
 \ref{aprox}):

%in this section is the following
\begin{teo}\label{cotaprinc0}  Under the hypothesis of Lemma \ref{existprobpar}, one of the two following alternatives hold:
\begin{enumerate}
\item \label{mejorcaso}$\lim_{n \to \infty}\sup_{B_{\frac{1}{\epsilon}}\times [0,T_n]}|w_n-h_{{\vec{q}}_n,\psi_n}|\to 0$,
or
\item \label{peorcaso} there is a constant $C$, independent of ${\vec{q}}_n$ and $T_n$
such that 
\be \sup_{B_{\frac{1}{\epsilon}}\times [0,T_n]}|w_n-h_{{\vec{q}}_n,\psi_n}|\leq C \sup_{B_{\frac{1}{\epsilon}}\times [0,T_n]}
|F_{{\vec{q}}_n}(w_n,\psi_n)-w_n|.\label{eqpeorcaso}\ee

\end{enumerate}
\end{teo}

\begin{rem}
Notice that in Theorem \ref{cotaprinc0}  it is possible to choose $T_n=\infty$ for every $n$.
\end{rem}

\begin{proof}[Proof of Theorem \ref{cotaprinc0}]
Consider sequences   of continuous functions  $\psi_n,
w_n $ satisfying \\ $\sup_{B_1} |\psi_n|,\ \sup_{B_1\times [0,T_n)}|w_n|\leq K$ and
 ${\vec{q}}_n\in \calQ$. Suppose that neither \equ{mejorcaso} nor \equ{peorcaso} hold. Then there are subsequences such that
\be \lim_{n \to \infty}\sup_{B_{\frac{1}{\epsilon_n}}\times [0,T_n]}|w_n-h_{{\vec{q}}_n,\psi_n}|
\not \to 0 \hbox{ and, } \label{nomejorcaso}\ee
\be \sup_{B_{\frac{1}{\epsilon_n}}\times [0,T_n]}|w_n-h_{{\vec{q}}_n,\psi_n}|= n \sup_{B_{\frac{1}{\epsilon_n}}\times 
[0,T_n]}
|F_{{\vec{q}}_n}(w_n,\psi_n)-w_n|.\label{nopeorcaso}\ee

 The a priori bounds shown in Theorem \ref{vvc} and the boundedness hypothesis imply that there is a constant independent of $n$ such that $|w_n-h_{{\vec{q}}_n,\psi_n}|\leq C$. Then, \equ{nopeorcaso} implies
\be \sup_{B_{\frac{1}{\epsilon_n}}\times 
[0,T_n]}
|F_{{\vec{q}}_n}(w_n,\psi_n)-w_n|\to 0.\label{to0dem}\ee

Applying inequality \equ{ineq0} recursively
we have that for every $0\leq T<\infty$ there is a constant that depends on $T$ (but independent of ${\vec{q}}_n$) such that
\be \sup_{B_{\frac{1}{\epsilon_n}}\times [0,T]}|w_n-h_{{\vec{q}}_n,\psi_n}|\leq C(T) \sup_{B_{\frac{1}{\epsilon_n}}\times 
[0,T]}
|F_{{\vec{q}}_n}(w_n,\psi_n)-w_n|.\label{peorcasoacotado}\ee
Therefore if the $T_n$ are unifromly bounded, case \equ{peorcaso} holds trivially, which contradicts
\equ{nopeorcaso}.
Hence we may assume $T_n\to \infty$. We will show that in this case 
$$\lim_{n \to \infty}\sup_{B_{\frac{1}{\epsilon_n}}\times [0,T_n]}|w_n-h_{{\vec{q}}_n,\psi_n}|
 \to 0, $$
contradicting a \equ{nomejorcaso}.

Let 
$$\tau=\{(S_n)_{n\in \nn}: 0\leq \ S_n\leq T_n, \ \lim_{n \to \infty}\sup_{B_{\frac{1}{\epsilon_n}}\times [0,S_n]}|w_n-h_{{\vec{q}}_n,\psi_n}|
 \to 0 \}.$$
For the set of sequences in $\rr_+$ we consider the topology defined by the basis of open sets given by $B_{\sigma}((S_n)_{n\in \nn})=\{(\tilde{S}_n)_{n\in \nn} :\ \tilde{S}_n \geq 0 \hbox{ and }  \sup_{n\in \nn}|S_n-\tilde{S}_n|\leq \sigma \}$ for any $\sigma >0$.
Notice that in particular inequality \equ{peorcasoacotado} implies that $\tau$ is a non-empty set, since at least $S_n=\inf_n T_n \in \tau$.

%every finete $T$ belongs to $\tau$.

\noindent{\bf Claim: $\tau$ is open}

Consider $(S_n)_n \in \tau.$ Let $\tilde{S}_n=\min\{S_n+\frac{2\beta}{M}, T_n\}$. Using inequality \equ{ineq0} we have
$$\sup_{B_{\frac{1}{\epsilon_n}}\times [S_n, \tilde{S}_n]}|w_n-h_{{\vec{q}}_n,\psi_n}|
 \leq \frac{1}{1-\beta}\left(2\sup_{B_{\frac{1}{\epsilon_n}}\times [S_n,\tilde{S}_n]} |F_{{\vec{q}}_n}(w_n,\psi_n)-w_n|+
\sup_{x\in B_{\frac{1}{\epsilon_n}}}|w_n-h_{{\vec{q}}_n,\psi_n}|(x, S_n)\right).$$
Since $\tilde{S}_n\leq T_n$  and $S_n\in \tau$, taking $n\to \infty$ we have that
$$\lim_{n\to \infty}\sup_{B_{\frac{1}{\epsilon_n}}\times [S_n, \tilde{S}_n]}|w_n-h_{{\vec{q}}_n,\psi_n}|=0,$$
and $  B_{\frac{2\beta}{M}}\bigcap \tau\subset \tau$. Hence $\tau$ is open.

\medskip
\noindent{\bf Claim: $\tau$ is closed}

Suppose that $S^k=(S_n^k)_n\in \tau$ satisfy 
$S^k\to \tilde{S}=(\tilde{S}_n)_n$ as $k\to \infty$.
By the definition of the topology we have that there is a $k_0$ such that for every $n\in \nn$ and $k\geq k_0$
holds $|S_n^k-\tilde{S}_n|\leq \frac{2\beta}{M}$. Using inequality 
\equ{ineq0} 
$$\sup_{B_{\frac{1}{\epsilon_n}}\times [S_n^{k_0},\tilde{S}_n]}|w_n-h_{{\vec{q}}_n,\psi_n}|
 \leq \frac{1}{1-\beta}\left( 2\sup_{B_{\frac{1}{\epsilon_n}}\times [S_n^{k_0},\tilde{S}_n]}|F_{{\vec{q}}_n}(w_n,\psi_n)-w_n|+
\sup_{x\in B_{\frac{1}{\epsilon_n}}}|w_n-h_{{\vec{q}}_n,\psi_n}|(x, S_n^{k_0})\right).$$
Using that  $(S_n^{k_0})_n\in \tau$ and \equ{to0dem}, when $n\to \infty$ we have 
$$\lim_{n \to \infty}\sup_{B_{\frac{1}{\epsilon_n}}\times [0,\tilde{S}_n]}|w_n-h_{{\vec{q}}_n,\psi_n}|=
\max\left\{\sup_{B_{\frac{1}{\epsilon_n}}\times [0,S_n^{k_0}]}|w_n-h_{{\vec{q}}_n,\psi_n}|, \sup_{B_{\frac{1}{\epsilon_n}}\times [S_n^{k_0},\tilde{S}_n]}|w_n-h_{{\vec{q}}_n,\psi_n}|\right\}\to 0.$$ Therefore $\tilde{S}\in \tau$ and $\tau$ is closed.

\medskip
Since $\tau$ is open, closed and non-empty we conclude that $\tau=\{(S_n)_{n\in \nn}: \ 0\leq S_n\leq T_n\}$. In particular $(T_n)_n\in \tau$, which contradicts \equ{nomejorcaso} and proves the Theorem.\end{proof}
 
 Following the proof of Theorem \ref{cotaprinc0}  we obtain:

 \begin{cor}\label{cotaprinc} 
Consider the sequences $\psi_n,
w_n $,
 ${\vec{q}}_n\in \calQ$ and $T_n>0$  as in Theorem \ref{cotaprinc0}. Assume in addition that there  are constants $C,m$ such that
$  \sup_{B_{\frac{1}{\epsilon_n}}\times [0,T_n]}
|F_{{\vec{q}}_n}(w_n,\psi_n)-w_n|\leq C\epsilon_n^m$. Then for every $\tilde{m}<m$
 holds either
\begin{enumerate}
\item \label{mejorcasoc}$\lim_{n \to \infty}\sup_{B_{\frac{1}{\epsilon_n}}\times [0,T_n]}\frac{|w_n-h_{{\vec{q}}_n,\psi_n}|}{\epsilon_n^{\tilde{m}}}\to 0$,
or
\item \label{peorcasoc} there is a constant $C$, independent of ${\vec{q}}_n$ and $T_n$
such that 
\be \sup_{B_{\frac{1}{\epsilon_n}}\times [0,T_n]}\frac{|w_n-h_{{\vec{q}}_n,\psi_n}|}{{\epsilon^{\tilde{m}}
}}\leq C \sup_{B_{\frac{1}{\epsilon_n}}\times [0,T_n]}
\frac{|F_{{\vec{q}}_n}(w_n,\psi_n)-w_n|}{\epsilon_n^{\tilde{m}}}.\label{eqpeorcasoc}\ee

\end{enumerate}

 In particular, there is a constant $C$ such that
 $$\sup_{B_{\frac{1}{\epsilon_n}}\times [0,T_n]}|w_n-h_{{\vec{q}}_n,\psi_n}|\leq C \epsilon_n^{\tilde{m}}.$$

 \end{cor}

\medskip

Now we would like to rescale the estimates of the previous Theorem and Corollary to the unit ball. Namely, instead of considering the function
$h_{\vec{q},\psi^\epsilon _\epsilon}: B_{\frac{1}{\epsilon}}\times [0,\frac{T}{\epsilon^2}]\to \rr^2$ we define 
the function $k_{\vec{q},\psi^\epsilon _\epsilon}:B_1\times [0,T]\to\rr^2 $  by 
\be  k_{\vec{q},\psi^\epsilon _\epsilon}(x,t)=h_{\vec{q}}\left(\frac{x}{\epsilon},\frac{y}{\epsilon^2}\right).\label{defk}\ee
Notice that under this definition for every $\epsilon>0$ we can write the left hand side of equation \equ{eqpeorcaso} as
 $$\sup_{B_{\frac{1}{\epsilon}}\times [0,\frac{T}{\epsilon^2}]}
\left| h_{\vec{q},\psi^\epsilon _\epsilon}
\left(x,t\right)
-w_\epsilon(x,t)\right|= \sup_{B_1\times [0,T]}
\left| k_{\vec{q}}(x,t)
-w_\epsilon^\epsilon(x,t)\right|,$$
where $w_\epsilon^\epsilon(x,t)=w_\epsilon\left(\frac{x}{\epsilon}, \frac{t}{\epsilon^2}\right)$.

Now we would like to rescale the right hand side of inequality  \equ{eqpeorcaso}.
Notice that by applying the function $F_{{\vec{q}}}$ to any function pair of continuous functions $w_\epsilon, \phi_\epsilon$ we obtain a continuous function
 $F_{{\vec{q}}}(w_\epsilon^\epsilon ,\psi^\epsilon_\epsilon):B_{\frac{1}{\epsilon}}\times[0,\frac{T}{\epsilon^2}]\to \rr^2$, which satisfies (via Duhamel's formula) the following eqaution:
\begin{align*}
\frac{d F_ {\vec{q}}(w_\epsilon,\psi_\epsilon) }{d t}-\Delta F_ {\vec{q}}(w_\epsilon,\psi_\epsilon)
+\frac{\nabla_v W(w_\epsilon+V_{\vec{q}})}{2}&=\Delta V_{\vec{q}}
 \hbox{ in }B_{\frac{1}{\epsilon}}\times [0,\frac{T}{\epsilon^2}]\\
F_{\vec{q}}(w_\epsilon, \psi_\epsilon)(x,t)&=0 \hbox{ on }\partial B_{\frac{1}{\epsilon}}\\
F_{\vec{q}}(w_\epsilon, \psi_\epsilon)(x,0)&=\psi_\epsilon(x).\end{align*}
  
Let us define the function $\calL_{\vec{q}}:\calC_{B_1\times[0,T]}\times \calC_{B_1}\to \calC_{B_1\times[0,T]}$ as
$\calL_{\vec{q}}(w_\epsilon^\epsilon, \psi_\epsilon^\epsilon)(x,t)=F_{\vec{q}}(w_\epsilon,\psi_\epsilon)\left(\frac{x}{\epsilon},\frac{y}{\epsilon^2}\right)$, where as before
$w_\epsilon^\epsilon(x,t)=w_\epsilon\left(\frac{x}{\epsilon}, \frac{t}{\epsilon^2}\right)$ and similarly 
$\psi_\epsilon^\epsilon(x,t)=\psi_\epsilon\left(\frac{x}{\epsilon}, \frac{t}{\epsilon^2}\right)$. A simple computation shows that for any $w_\epsilon^\epsilon,\psi_\epsilon^\epsilon$ the function
obtained by evaluating  $\calL_{\vec{q}}$ at $(w_\epsilon^\epsilon, \psi_\epsilon^\epsilon)$, denoted by $\calL_{\vec{q}}(w_\epsilon^\epsilon, \psi_\epsilon^\epsilon)$,  satisfies

\begin{align*}
\frac{d \calL_ {\vec{q}}(w_\epsilon^\epsilon,\psi_\epsilon^\epsilon) }{d t}-\Delta \calL_ {\vec{q}}(w_\epsilon^\epsilon,\psi_\epsilon^\epsilon)
+\frac{\nabla_v W(w_\epsilon^\epsilon+U_{\vec{q}})}{2\epsilon^2}&=\Delta U_{\vec{q}}
 \hbox{ in }B_{1}\times [0,T]\\
\calL_{\vec{q}}(w_\epsilon^\epsilon, \psi_\epsilon^\epsilon)(x,t)&=0 \hbox{ on }\partial B_{1}\\
\calL_{\vec{q}}(w_\epsilon^\epsilon, \psi_\epsilon^\epsilon)(x,0)&=\psi_\epsilon^\epsilon(x).\end{align*}

Using again Duhamel's formula we conclude that 
\begin{align}\calL_ {\vec{q}}(w_\epsilon^\epsilon,\psi_\epsilon^\epsilon)(x,t)
=\int_0^{t} \int_{B_{1}} &\calH_{B_1}\left(x,y,t-s\right)\left(-\frac{\nabla_v W(w_\epsilon^\epsilon+
U_{\vec{q}})(y,s)}{\epsilon^2}\right.
\notag \\
&  \left.\frac{}{}+\Delta U_{\vec{q}} (y)
\right)dyds + \int_{B_{1}}\calH_{B_1}\left(x,y,t\right)\psi_\epsilon^\epsilon (y)dy.\label{defgfunc}\end{align}
In particular, we have that $k_{{\vec{q}},\psi}$ defined by \equ{defk} is a fixed point of $\calL_{\vec{q}}(\cdot,\psi)$.

Hence, the right hand side of equation \equ{eqpeorcaso} reads
$$\sup_{B_{\frac{1}{\epsilon}}\times[0,\frac{T}{\epsilon^2}]}|F_{\vec{q}}(w_\epsilon,\psi_\epsilon)-w_\epsilon|=\sup_{B_1\times[0,T]} |\calL_{\vec{q}}(w^\epsilon_\epsilon,\psi_\epsilon)-w_\epsilon^\epsilon|.$$
%On the other hand, if the function $w_\epsilon^\epsilon$ is regular enough  and $w_\epsilon(x,t)=0 $ for $|x|=1$, we can express it as (using again Duhamel's formula)

%$$w^\epsilon_\epsilon(x,t)=\int_0^{t} \int_{B_{1}}\calH_{B_1}\left(x,y,t-s\right) P w^\epsilon_\epsilon(y,s) dyds + \int_{B_{1}}\calH_{B_1}\left(x,y,t\right)w^\epsilon_\epsilon (y,0)dy,$$
%where $P$ is the heat operator, namely  $Pw=\frac{dw}{dt}-\Delta w$.

In this context we  can re-formulate Theorem \ref{cotaprinc0} (dropping the super-indeces to simplify the notation) as
 
\begin{teo}\label{cotaprinc1}
Let $k_{{\vec{q}},\psi}$ be  defined by \equ{defk}. Then is the unique fixed point of $\calL_{{\vec{q}}}(\cdot,\psi)$. Moreover, for any fixed
 $K>0$ and sequences of continuous functions $\psi_n,
w_n $ satisfying $\sup |\psi_n|,\ \sup|w_n|\leq K$ and vectors
 ${\vec{q}}_n\in \calQ$ and $T_n>0$ holds either

\begin{enumerate}
\item $\sup_{B_1\times [0,T]}
\left| k_{\vec{q},\psi _\epsilon}
\left(x,t\right)
-w_\epsilon(x,t)\right|\to 0$ or

\item there is a constant $C$, independent of $\epsilon, \sigma$ and $T$
such that 
$$\sup_{B_1\times [0,T]}\left|k_{\vec{q},\psi_\epsilon}(x,t)
-w_\epsilon(x,t)\right|\leq C 
\sup_{B_1\times [0,T]}
|\calL_{{\vec{q}}}(w_\epsilon ,\psi_\epsilon)-w_\epsilon|,$$
where ${\vec{q}}=(\epsilon,\sigma, \alpha)$.
\end{enumerate}
\end{teo}

Now we can devote ourselves to prove Theorem \ref{aprox}. We divide the proof into two steps: Lemma \ref{kgoesto0} and Lemma \ref{tainf}.

 Notice  first 
that the function $k_{\vec{q}, \psi}$ defined by \equ{defk} is a solution to the following equation:  
\begin{align}
Pk_{\vec{q}, \psi}
+\frac{\nabla_v W(k_{\vec{q}, \psi}+U_{\vec{q}})}{2\epsilon^2}&=\Delta U_{\vec{q}} \hbox{ in }B_{1}\label{ecpar'}\\
k_{\vec{q}, \psi}(x,t)&=0 \hbox{ on }\partial B_{1}\label{cbpar'}\\
k_{\vec{q}, \psi}(x,0)&=\psi.\label{cipar'}
\end{align}
where $Pk_{\vec{q}, \psi}=\frac{d k_{\vec{q}, \psi}}{d t}-\Delta k_{\vec{q}, \psi}$. In order to simplify the 
notation, when $\psi\equiv 0$ we will simply denote this solution by $k_{\vec{q}}$
 (instead of $k_{\vec{q},0}$). In Lemma \ref{kgoesto0}  we show that 
$$\lim_{\epsilon\to 0}\sup_{B_1\times[0,\infty]}|k_{\vec{q}}(x,t)|=0.$$
In order to do this computation we will use several estimates from the Appendix.  Thereafter
we will conclude the proof of  Theorem \ref{aprox}  by showing  in Lemma \ref{tainf} that for every fixed $\epsilon$ there is a sequence $0<t_n \nearrow \infty$ satisfying
$$\lim_{n\to \infty}\sup_{B_1}|k_{\vec{q}}(x,t_n)-u_\epsilon+U_{\vec{q}}|=0.$$

\begin{lem} \label{kgoesto0}
Let $k_{\vec{q}}$ be the solution to \equ{ecpar'}-\equ{cbpar'}-\equ{cipar'}, for $\psi=0$. Then
 $$\lim_{\epsilon\to 0}
\sup_{B_1\times[0,\infty]}|k_{\vec{q}}(x,t)|=0.$$
\end{lem}

\begin{proof}
Suppose that 
$$\sup_{B_1\times[0,\infty)}|k_{\vec{q}}|\not \to 0.$$
Theorem  \ref{cotaprinc1}  implies that (by choosing $w_\epsilon=\psi_\epsilon=0$)
\be \sup_{B_1\times [0,\infty)}|k_{\vec{q}}|\leq C\sup_{B_1\times[0,\infty)}|\calL_{\vec{q}}(0,0)|. \label{boundk}\ee

Set $S_\epsilon=\sup_{B_1\times[0,\infty)}|\calL_{\vec{q}}(0,0)|$ (possibly infinity). Fix $\delta>0$ and notice that, by definition of supremum, there is a $t_\epsilon$ such that
$\sup_{x\in B_1}|\calL_{\vec{q}}(0,0)(x,t_\epsilon)-S_\epsilon|\leq \delta$ (or when
$S_\epsilon=\infty$  pick $t_\epsilon$ such that $\sup_{x\in B_1}|\calL_{\vec{q}}(0,0)(x,t_\epsilon)|\geq \delta^{-1}$).

We will show that, independently of $\delta$, holds $\sup_{x\in B_1}|\calL_{\vec{q}}(0,0)|(x,t_\epsilon)\to 0$ as $\epsilon \to 0$ (notice that this immediately 
contradicts $S_\epsilon=\infty$). 
%Given a sequence of $t_\epsilon$ we need to consider two cases: 
%\begin{enumerate}
%\item There is a $T>0$ such that $t_\epsilon \leq T$ for every $\epsilon$, or
%\item there is a subsequence such that $t_\epsilon\nearrow \infty$.
%\end{enumerate}
%Hence we will 
Recall first that 
\be \calL_ {\vec{q}}(0,0)(x,t)
=\int_0^{t} \int_{B_{1}}\calH_{B_1}\left(x,y,t-s\right)\left(-\frac{\nabla_v W(
U_{\vec{q}})(y,s)}{\epsilon^2}
+\Delta U_{\vec{q}} (y)
\right)dyds.\label{calgat00}\ee

%and

%$$I_3(x,t)=  \int_{B_1} |u_\epsilon-U_{\vec{q}}| dy. \hspace{6cm}$$

Notice that for $|x|\leq \epsilon^\alpha-E$ we have
$$\frac{-\nabla_v W(U_{\vec{q}})}{\epsilon^2}+\Delta U_{\vec{q}}=
\frac{-\nabla_v W(u^\epsilon_\sigma)}{\epsilon^2}+\Delta u^\epsilon_\sigma=0.$$
Hence, \equ{calgat00} implies
\be |\calL_{\vec{q}}(0,0)|(x,t)\leq I_1(x,t)+I_2(x,t), \label{cotagfunc}\ee
where
$$I_1(x,t)= \int_0^{t}\int_{\{|y|\geq \epsilon^\alpha \}}\calH_{B_1}(x,y,t-s)\left|\frac{-\nabla_v W(\phi_\epsilon)}{\epsilon^2}+\Delta \phi_\epsilon\right|(y,s)
dyds  \hspace{1cm}$$
\begin{align*}I_2(x,t)=& \int_0^t\int_{\{\epsilon^\alpha-E\leq |y|\leq \epsilon^\alpha \}}\calH_{B_1}(x,y,t-s) \left|
  \frac{-\nabla_v W(U_{\vec{q}})}{\epsilon^2}+\eta^\epsilon_\alpha \Delta \phi_\epsilon\right.
\\ & \hspace{1cm}\frac{}{}+\Delta (\eta^\epsilon_\alpha )\left(h_\sigma^\epsilon-\phi_\epsilon\right)+\left.\nabla(\eta^\epsilon_\alpha )\cdot D\left( u^\epsilon_\sigma- \phi_\epsilon\right)
  \right| (y,s)dy ds.  \frac{}{}\end{align*} 

Now   we find bounds for  $I_1$ and $I_2$. For each of these integrals we will consider two ranges for the variable $t$, namely
$t\leq T$ and $t\geq T$, where  $T>0$ is fixed positive constant.%Fix $\delta>0$.
\begin{itemize}
\item{Bounds over $I_1$:}

%{\bf Claim:} There is a $0<\delta<1$ such that 

%$$\sup_{(x,t)\in B_1\times}$$

%Suppose that such $\delta$ does not exist. Then, there is a sequence

%By definition of $\phi_\epsilon$ we have that 
Since $\epsilon<\epsilon^\alpha$ (when $\epsilon<1$) we have that for every $|x|\geq \epsilon^\alpha$ the function $\eta(x)\equiv 0$ and for such $x$ we have

\begin{align*} \Delta \phi_\epsilon(x)=\frac{1}{\epsilon^2}\left( \epsilon^2
\Delta \eta_{6}\frac{}{} \right.& \zeta_{31}\left(d_0  \left(\frac{x}{\epsilon}\right)\right)+  \eta_{6} \zeta''_{31}\left(d_0\left(\frac{x}{\epsilon}\right)\right)\\
& + 2\epsilon \nabla \eta_{6}\cdot \nabla d_0\left(\frac{x}{\epsilon}\right) \zeta'_{6}\left(d_0\left(\frac{x}{\epsilon}\right)\right)+\epsilon ^2 \Delta \eta_{5}c_i \hspace{1cm}\\
+ \sum_{i=1}^3 &\epsilon^2 \Delta \eta_{2i}\zeta_{ii+1}\left(d_i\left(\frac{x}{\epsilon}\right)\right)+  \eta_{2i} \zeta''_{ii+1}\left(d_i\left(\frac{x}{\epsilon}\right)\right) \\
& \left.+2 \epsilon \nabla \eta_{2i}\cdot \nabla d_i\left(\frac{x}{\epsilon}\right) \zeta'_{ii+1}\left(d_i\left(\frac{x}{\epsilon}\right)\right)+\epsilon^2\Delta \eta_{2i-1}c_i \right). \end{align*}
%and $$\phi_\epsilon(x,t)=\phi\left(\frac{x}{\epsilon}\right)$$
  %$\{[2\pi-\delta,2\pi]\bigcup[0,\delta], [\frac{\delta}{2},\theta_1-\frac{\delta}{2}],[\theta_1-\delta,\theta_1+\delta],
%[\theta_1+\frac{\delta}{2},\theta_2-\frac{\delta}{2}],[\theta_2-\delta,\theta_2+\delta],[\theta_2+\frac{\delta}{2},2\pi-\frac{\delta}{2}]\}$
%Let $$\phi(x,t)=c_0\eta_0(x,t)+(1-\eta_0(x,t))\sum_{i=1}^3\eta_{2i-1}(\theta)g_i(d_i(x,t))+\eta_{2i}(\theta)c_i $$
%and $$\phi_\epsilon(x)=\phi\left(\frac{x}{\epsilon}\right)$$
Since the functions $\eta_j$ depend only on the angle $\theta$ we have that
\begin{align*}\Delta \eta_j=&\frac{\eta_j''}{r^2} \hbox{ and }\\
|\nabla \eta_j|\leq &|\eta_j'|.\end{align*}
In particular for $|x|\geq \epsilon^\alpha$
 \begin{align*}|\Delta \eta_j|\leq 
 & \frac{4 |\eta''_j| }{\epsilon^{2\alpha}} \hbox{ and }\\
|\nabla \eta_j|\leq &|\eta_j'|.\end{align*}
Recall that for $\theta\in \left[\theta_i-\frac{\delta}{2},\theta_i+\frac{\delta}{2}\right]$ we have 
$\eta_{2i}\equiv 1 \hbox{ and }  \eta_{j}\equiv 0 \hbox{ for every } j\ne 2i.$
 Then  
 \be  \frac{\nabla_v W(\phi_\epsilon)}{\epsilon^2}+\Delta \phi_\epsilon=0 
\hbox{ for }\theta\in \left[\theta_i-\frac{\delta}{2},\theta_i+\frac{\delta}{2}\right].\label{cotacerca}\ee

Now we need to find bounds for $\theta \in \left[\theta_i+\frac{\delta}{2},\theta_{i+1}-\delta\right]$.
Notice first that  
 \begin{align}|\Delta \eta(\theta)|=&\left|\frac{\eta''(\theta)}{r^2}\right|
 \leq \frac{K}{r^2}  \leq\frac{K}{\epsilon^{2\alpha}} \hbox{ for }|x|\geq \epsilon^\alpha  \label{cotalapeta}\\
 |\nabla \eta|=&\left|\frac{\eta'}{r}\right|\leq \frac{K}{r} \leq \frac{K}{\epsilon^\alpha}  \hbox{ for }|x|\geq \epsilon^\alpha .\label{cotagraeta} \end{align}
  Notice also that  $\eta_{j}\ne 0$ only for $j=21,  2i-1$ and $$\eta_{2i}+ \eta_{2i-1}=1.$$
Hence 
 \begin{align*} \Delta \phi_\epsilon &=\frac{1}{\epsilon^2}\left( \epsilon^2 \Delta \eta_{2i}(\theta)\zeta_{ii+1}\left(d_i\left(\frac{x}{\epsilon}\right)\right)+  \eta_{2i}(\theta) \zeta''_{ii+1}\left(d_i\left(\frac{x}{\epsilon}\right)\right) 
\right.
\\   &\hspace{1cm} \left.+2 \epsilon \nabla \eta_{2i}(\theta)\cdot \nabla d_i\left(\frac{x}{\epsilon}\right) 
\zeta'_{ii+1}\left(d_i\left(\frac{x}{\epsilon}\right)\right)+
\epsilon^2\Delta \eta_{2i-1}(\theta)c_i \right) \\
&=\Delta \eta_{2i}(\theta)\left(\zeta_{ii+1}\left(d_i\left(\frac{x}{\epsilon}\right)\right)-c_i\right)+  \eta_{2i}(\theta)\frac{-\nabla_v W( \zeta_{ii+1})}{\epsilon^2}\left(d_i\left(\frac{x}{\epsilon}\right)\right) 
\\   &\hspace{1cm} \left.+2\frac{1}{ \epsilon} \nabla \eta_{2i}(\theta)\cdot \nabla d_i\left(\frac{x}{\epsilon}\right) 
\zeta'_{ii+1}\left(d_i\left(\frac{x}{\epsilon}\right)\right). \right.\end{align*}
Using Hypothesis \ref{geod} we have that there  are constants $K,c>0$ such that
 \be \left| \frac{\nabla_v W(\phi_\epsilon)}{\epsilon^2}+\Delta \phi_\epsilon\right|\leq K\frac{e^{-c\frac{d_i}{\epsilon}}}{\epsilon^2} \hbox{ for }|x|\geq \epsilon^\alpha \hbox{ and }\theta \in \left[\theta_i+\frac{\delta}{2},\theta_{i+1}-\frac{\delta}{2}\right]. \label{cotalejos00}\ee

Furthermore, for $|x|>\epsilon^\alpha$ and $\theta\in \left[\theta_i+\frac{\delta}{2},\theta_{i+1}-\frac{\delta}{2}\right]$ we have 
$$|d_i|\geq \epsilon^\alpha \sin \delta.$$
Hence, 
 \be \left| \frac{\nabla_v W(\phi_\epsilon)}{\epsilon^2}+\Delta \phi_\epsilon\right|\leq K\frac{e^{-c\frac{\epsilon^\alpha \sin \delta}{\epsilon}}}{\epsilon^2} \hbox{ for }|x|>\epsilon^\alpha \hbox{ and }\theta\in \left[\theta_i+\frac{\delta}{2},\theta_{i+1}-\frac{\delta}{2}\right].\label{cotalejos}\ee

Now we proceed to find bounds in two different cases:

\begin{enumerate}

\item Suppose that $t\leq T$.
Equations \equ{cotacerca} and \equ{cotalejos} imply
$$ I_1(x,t)\leq K\frac{e^{-c\frac{\epsilon^\alpha \sin \delta}{\epsilon}}}{\epsilon^2} \int_0^{t}\int_{\{|x|\geq \epsilon^\alpha \}}\calH_{B_1}(x,y,t-s) dyds.$$

Using Lemma \ref{cotH} we have
\be I_1(x,t) \leq K\frac{e^{-c\frac{\epsilon^\alpha \sin \delta}{\epsilon}}}{\epsilon^2} \hbox{ for every }x\in B_1 \hbox{ and } 0\leq t\leq T.\label{i1} \ee

\item Suppose that $t\geq T$
%---------------------------------------------------------------------------
%Standard heat kernel estimates show that
%$\calH_{B_1}(x,y,t)=\frac{e^{-\frac{|x-y|^2}{4t}}}{4\pi t}+\calH^{bd}(x,y,t)$, where $ \calH^{bd}(x,y,t)=O\left(\left[\frac{d(x)}{\sqrt{t}}+\frac{d(y)}{\sqrt{t}}\right]^{-\infty}\right)$, where $d(x)$ the distance of $x$ to the boundary.

%Moreover,  $ \calH^{bd}(x,y,t)=t^{-1}\tilde{A}\left(t, \frac{x'-y'}{\sqrt{t}},x_n, y_n, y'\right)$ where $\tilde{A}(t, X',\psi_n, \eta_n, y')=O\left((\psi_n+ \eta_n+|X'|)
%^{-\infty}\right)$
%-----------------------------------------------------------------------------

Let $$f_\epsilon= \left| \frac{\nabla_v W(\phi_\epsilon)}{\epsilon^2}+\Delta \phi_\epsilon\right|$$ and fix $\delta>0$.
Now we divide $I_1$ in the three following integrals:
\begin{align*}I_{11}(x,t)=&\int_0^{t-\delta}\int_{\{|y|\geq \epsilon^\alpha\}
\bigcap\{|x-y|\leq \frac{\sqrt{t-s}}{t} \}} \calH_{B_1}(x,y,t-s)f_\epsilon(y,s)dyds, \\
I_{12}(x,t)=&\int_0^{t-\delta}\int_{\{|y|\geq \epsilon^\alpha\}
\bigcap\{|x-y|\geq \frac{\sqrt{t-s}}{t} \}} \calH_{B_1}(x,y,t-s)f_\epsilon(y,s)dyds,  \\
I_{13}(x,t)=&\int_{t-\delta}^{\delta}\int_{\{|y|\geq \epsilon^\alpha\}} \calH_{B_1}(x,y,t-s)f_\epsilon(y,s)dyds. \end{align*}
Then $$I_1=I_{11}+I_{12}+I_{13}.$$

By Theorem \ref{H2} we have that $|\calH_{B_1}(x,y, t-s)|\leq \frac{C}{(t-s)}$, then
\begin{align*} I_{11}(x,t)\leq & C \int_0^{t-\delta}\frac{\sup_{|y|\geq \alpha}f_\epsilon}{(t-s)}\frac{(t-s)}{t^2}\pi ds\\ =& \frac{\sup_{|y|\geq \alpha}f_\epsilon}{t^2}(t-\delta) \\ \leq&
C \frac{e^{-c\frac{d_i}{\epsilon}}}{t \epsilon^2}.\end{align*}

Using again Theorem \ref{H2} , for $|x-y|\geq \frac{\sqrt{t-s}}{t} $ we have  $|\calH_{B_1}(x,y, t-s)|=O\left(\left[\frac{1}{t}\right]^{-\infty}\right)$. In particular there is a constant $C$ such that  $|\calH_{B_1}(x,y, t-s)|\leq \frac{C}{t},$ then

\begin{align*}I_{12}(x,t) \leq & \int_0^{t-\delta} \frac{C}{t} \int_{B_1}f_\epsilon(y)dy\\ \leq  & t \frac{C}{t}\int_{B_1}f_\epsilon(y)dy\\ \leq &  C  \frac{e^{-c\frac{d_i}{\epsilon}}}{ \epsilon^2}.\end{align*}

Finally, using Lemma \ref{cotH} we have
$$I_{13}(x,t)\leq \delta \sup f_\epsilon \leq  C  \frac{e^{-c\frac{d_i}{\epsilon}}}{ \epsilon^2}.$$
%Now suppose that such $\delta>0$ does not exist.
%Then for every $n\in \nn$ there is a $t_n$ and a $x_n \in B_1\setminus B_{1-\frac{1}{n}}$ such that $$I_1(x,t_n)\leq I_1(x_n,t_n).$$

%Consider any convergent subsequences of $x_n\to x$ and $t_n\to t$. Since
  %$x_n \in B_1\setminus B_{1-\frac{1}{n}}$, it must hold that $x\in \partial B_1$. By definition of $I_1$ we have that 
%$$I_1(x,t)=0 \hbox{ for } x\in \partial B_1 \hbox{ and every }t.$$
%Hence $$I_1(x_n,t_n)\to 0.$$
%for every $t\in[0, T_\epsilon]$, the maximum of $I_1(\cdot,t)$ is attained at a point $x\in B_{1-\delta}$.

Combining the previous estimates we obtain
\be I_1(x,t)\leq C  \frac{e^{-c\frac{d_i}{\epsilon}}}{ \epsilon^2}  \hbox{ for every }x\in B_1 \hbox{ and }  t\geq T. \label{i1'}\ee
\end{enumerate}
\item{ Bounds over $I_2$:}

Using the definitions of $U_{\vec{q}}$, $\phi_\epsilon$, Theorem \ref{vvc} and Lemma \ref{lema2} 
 we have
$$
\left|  \frac{-\nabla_v W(U_{\vec{q}})}{\epsilon^2}+\eta^\epsilon_\alpha \Delta \phi_\epsilon\right|\leq \frac{C}{\epsilon^2}$$
$$
|\Delta (\eta^\epsilon_\alpha )\left(h_\sigma^\epsilon-\phi_\epsilon\right)|\leq \frac{C}{E^2}$$
$$\left|\nabla(\eta^\epsilon_\alpha )\cdot D\left( u^\epsilon_\sigma- \phi_\epsilon\right)
   \right| \leq \frac{C}{E \epsilon}.$$ 
 Hence:
\begin{enumerate} 
\item For $t\leq T$
$$I_2(x,t) \leq C\int_0^t\int_{\epsilon^\alpha-E\leq |x|\leq \epsilon^\alpha}\calH(x,y,t-s)\left(\frac{1}{\epsilon^2}+\frac{1}{E^2}+\frac{1}{E \epsilon}\right)dyds.$$

%Standard estimates for the heat kernel 
Theorem \ref{H2} implies that for $t-s\geq \epsilon^{m+2}$ there is a constant $C$ independent of
 $x,y$ such that $|\calH(x,y,t-s)|\leq \frac{C}{\epsilon^{m+2}} $. Moreover, by definition 
$\epsilon^\alpha \leq E=\epsilon^{\alpha}(2-\epsilon^{2m+4})\leq 2\epsilon^\alpha$. Hence 
\begin{align*} I_2\leq &
 \int_0^{t-\epsilon^{m+2}}\int_{\epsilon^\alpha-E\leq |x|\leq \epsilon^\alpha}\frac{C}{\epsilon^{m+2}}\left(\frac{1}{\epsilon^2}+\frac{1}{\epsilon^{2\alpha}}+\frac{1}{ \epsilon^{1+\alpha}}\right)
dyds\\ &+\int_{t-\epsilon^{m+2}}^t\int_{\epsilon^\alpha-E\leq |x|\leq \epsilon^\alpha}\calH_{B_1}(x,y,t-s) \frac{1}{\epsilon^2}
\left(1+\epsilon^{2-2\alpha}+\epsilon^{1-\alpha}
\right)dyds\\ 
\leq &
\frac{C}{\epsilon^{m+4}}\int_0^{t-\epsilon^{m+2}}
\left(1+\epsilon^{2-2\alpha}+\epsilon^{1-\alpha}
\right)\pi(\epsilon^{2\alpha}-(\epsilon^\alpha-E)^2)ds
\\ &
+\frac{C}{\epsilon^2}\int_{t-\epsilon^{m+2}}^t\int_{B_1}\calH(x,y,t-s) 
dyds.
\end{align*}

Using that $t\leq T$, Lemma \ref{cotH} and the definition of $E$ we conclude
\begin{align}I_2(x,t)\leq &\frac{C}{\epsilon^{m+4}}E(2\epsilon^\alpha-E)
+\frac{C}{\epsilon^2}\epsilon^{m+2} \notag \\ 
 \leq& \frac{C}{\epsilon^{m+4}}\epsilon^\alpha
\epsilon^{2m+4-\alpha}
+C \epsilon^m \notag\\
\leq & C\epsilon^{m} \hbox{ for } x\in B_1 \hbox{ and }0\leq t \leq T.\label{i2}\end{align}

\item For $t\geq T$

The previous estimates show that the integrand of $I_2$ can be bounded by $\frac{C}{\epsilon^2}$. Dividing up the integral as we did for $I_1$ we obtain

\begin{align*}I_2\leq & \int_0^{t-\epsilon^{m+2}}
\int_{\{\epsilon^\alpha-E\leq |y|\leq \epsilon^\alpha\}\bigcap\{|x-y|\leq \frac{\sqrt{t-s}}{t} \}} \calH_{B_1}(x,y,t-s)\frac{C}{\epsilon^2}dyds \\ 
& +
\int_0^{t-\epsilon^{m+2}}\int_{\{\epsilon^\alpha-E\leq |y|\leq \epsilon^\alpha\
\}
\bigcap\{|x-y|\geq \frac{\sqrt{t-s}}{t} \}} \calH_{B_1}(x,y,t-s)\frac{C}{\epsilon^2}dyds \\ &
+\int_{t-\epsilon^{m+2}}^{t}\int_{\{\epsilon^\alpha-E\leq |y|\leq \epsilon^\alpha \}} \calH_{B_1}(x,y,t-s)\frac{C}{\epsilon^2}dyds. \end{align*}

Using H\"older's inequality in the first integral for $p<2$ we get

$$I_2\leq   \int_0^{t-\epsilon^{m+2}}
\left(\int_{\{\epsilon^\alpha-E\leq |y|\leq \epsilon^\alpha \}}\frac{C}{\epsilon^{2p}}dy\right)^{\frac{1}{p}}
\left(
\int_{\{|x-y|\leq \frac{\sqrt{t-s}}{t} \}} \calH^q_{B_1}(x,y,t-s)
dy\right)^{\frac{1}{q}}ds $$
$$+
\int_0^{t-\epsilon^{m+2}}
%O((\frac{1}{t})^{-\infty})
\int_{\{\epsilon^\alpha-E\leq |y|\leq \epsilon^\alpha\
\}\bigcap\{|x-y|\geq \frac{\sqrt{t-s}}{t} \}}\calH_{B_1}(x,y,t-s) \frac{1}{\epsilon^2}dy ds+C
\frac{\epsilon^{m+2}}{\epsilon^2}.$$

As before, Theorem \ref{H2} implies
$|\calH_{B_1}|\leq \frac{C}{t-s}$ and that  for $|x-y|\geq\frac{t-s}{t}$ holds
$\calH_{B_1}(x,y,t-s)=O((\frac{1}{t})^{-\infty})$, therefore
\begin{align*}I_2 \leq & C \left[\int_0^{t-\epsilon^{m+2}}
\left(\frac{C\epsilon^{2m+4}}{\epsilon^{2p}}\right)^{\frac{1}{p}}
\frac{1}{t-s} 
\left(
 \frac{t-s}{t^2} 
\right)^{\frac{1}{q}}ds \right. \\
&\quad + \left.
\int_0^{t-\epsilon^{m+2}}
\int_{\{\epsilon^\alpha-E\leq |y|\leq \epsilon^\alpha\
\}}
\frac{C}{t}
\frac{\epsilon^{2m+4}}{\epsilon^2} +
C\frac{\epsilon^{m+2}}
{\epsilon^2}\right]
 \\ \leq &
C\left(\epsilon^{2m+4-2p}\right)^{\frac{1}{p}}
 \frac{t^{\frac{1}{q}}-\epsilon^{\frac{m+2}{q}}}{t^\frac{2}{q}} 
+
t \frac{C}{t}
\epsilon^{2m+2} +
C\epsilon^m
.\end{align*}
Therefore, for $t\geq T$ and  $p<2$ holds
\be I_2(x,t)\leq C\left(\frac{\epsilon^{\frac{2m+4-2p}{p}}}{T^{\frac{1}{q}}}+\epsilon^2+\epsilon^2 \right)\leq C\epsilon^2\label{i2'}\ee
\end{enumerate}
\medskip
\end{itemize}

Now we can conclude the result of Lemma by combining \equ{i1}, \equ{i1'}, 
\equ{i2} in \equ{i2'} in \equ{cotagfunc} and \equ{boundk}. More precisely:
$$ \sup_{B_1\times[0,\infty)} |k_{\vec{q}}|\leq C \frac{e^{-c\frac{\epsilon^\alpha \sin \delta}{\epsilon}}}{\epsilon^2} +C\epsilon^m\leq C\epsilon^{m},$$
where $C$ depends on $\alpha$ and $m$.
This implies the desired Lemma. \end{proof}

To finish the proof of Theorem \ref{aprox} we need the following Lemma
\begin{lem}\label{tainf}

Fix $\epsilon>0$ and let $k_{\vec{q}}$ be the solution \equ{ecpar'}-
\equ{cbpar'}
-\equ{cipar'}. Then, there is a sequence of times $t_n\nearrow \infty$ such that $$\lim_{n\to \infty}\sup_{B_1}|k_{\vec{q}}(x,t_n)-u_\epsilon(x)+U_{\vec{q}}(x)|=0.$$ 
\end{lem}
\begin{proof}

Corollary \ref{cotkder}  in the appendix shows that for every $t>0$ there is a constant $C$
 such that $|Dk_{\vec{q}}(x,t)|\leq \frac{C}{\epsilon}$. Similarly, by taking derivatives on the equation, we can find bounds over the second and third space
 derivatives (these bounds will depend on $\epsilon$). Since $\epsilon$ is fixed, using Arzela-Ascoli's Theorem we conclude for every sequence
$t_n\nearrow \infty$ there is a subsequence $k_{\vec{q}}(x,t_n)$ that converges in $\calC^2$.
Let us denote this limit by $k^\infty_{\vec{q}}(x)$ and the convergent subsequence $\{t_n\}_{n\in \nn}$ as well. 

We will show that
$k^\infty_{\vec{q}}(x)$ satisfies
\be \Delta k^\infty_{\vec{q}}(x)=\frac{\nabla_v W(k_{\vec{q}}^\infty+U_{\vec{q}})}{\epsilon^2}-\Delta U_{\vec{q}} \hbox{ for } x\in B_1\label{kinfec}\ee
\be k^\infty_{\vec{q}}|_{\partial B_1}=0.\label{cikinf}\ee

First we need to show that for every $\tau>0$ the sequence $k_{\vec{q}}(x,t_n+\tau)$ also converges in $\calC^2$ to $k^\infty_{\vec{q}}(x)$.
Define
$$\calJ(t)=\int_{B_1}\left(\frac{|\nabla k_{\vec{q}}|^2}{2}+\frac{ W(k_{\vec{q}}+U_{\vec{q}})}{\epsilon^2}- k_{\vec{q}} \cdot \Delta U_{\vec{q}}\right)(x,t) dx.$$

Using Theorem \ref{vvc} and the definition of $U_{\vec{q}}$ it is easy to see that $\calJ(t)$ is bounded below for every $t$. Moreover,
taking time derivative we have
 \begin{align*}\frac{d\calJ}{dt}
=& \int_{B_1}\left( \nabla k_{\vec{q}}\cdot
\nabla (k_{\vec{q}})_t+\frac{ \nabla W(k_{\vec{q}}+U_{\vec{q}})}{\epsilon^2}
\cdot(k_{\vec{q}})_t-\Delta U_{\vec{q}}\cdot (k_{\vec{q}})_t\right) (x,t)dx\\
=& \int_{B_1}\left[\left(-\Delta k_{\vec{q}}
+\frac{ \nabla W(k_{\vec{q}}+U_{\vec{q}})}{\epsilon^2}
-\Delta U_{\vec{q}}\right)\cdot (k_{\vec{q}})_t \right](x,t)dx\\
=& -\int_{B_1}|(k_{\vec{q}})_t|^2 (x,t) dx.
\end{align*}
Therefore $\calJ$ is bounded below and decreasing, hence it converges.  Moreover
for every fixed $\tau>0$
$$\int_{t_n}^{t_n+\tau}\int_{B_1}|k_{\vec{q}}|^2_t (x,s)dxds=\calJ(t_n)-\calJ(t_n+\tau)\to 0.$$
Since for every fixed $x$ we can write
$k_{\vec{q}}(x,t_n+\tau)-k_{\vec{q}}(x,t_n)=\int_{t_n}^{t_n+\tau}(k_{\vec{q}})_t (x,s)ds$,  we have that
\begin{align*}\int_{B_1}|k_{\vec{q}}(x,t_n+\tau)-k_{\vec{q}}(x,t_n)|dx\leq& \int_{t_n}^{t_n+\tau}\int_{B_1}|(k_{\vec{q}})_t|(x,s)dx ds
\\ 
\leq& C \left(\int_{t_n}^{t_n+\tau}\int_{B_1}|(k_{\vec{q}})_t|^2(x,s)dx ds\right)^{\frac{1}{2}}\to 0\hbox{ as }n\to \infty.\end{align*}
Hence $k_{\vec{q}}(x,t_n+\tau)-k_{\vec{q}}(x,t_n)$ converges to 0 almost 
everywhere. 
Let us show that this convergence is also uniform.
Suppose that $\sup_{x\in B_1}|k_{\vec{q}}(x,t_n+\tau)-k_{\vec{q}}(x,t_n)|\not \to 0$ as $n\to \infty$. Then there is a $\delta>0$ and a subsequence of times 
such that  \be \sup_{x\in B_1}|k_{\vec{q}}(x,t_n+\tau)-k_{\vec{q}}(x,t_n)|\geq 
\delta.\label{pasoindel}\ee As before, there is subsequence of these $\{t_n\}$ that converges uniformly. Since it converges almost everywhere to 0, the uniform limit must be 0 contradicting \equ{pasoindel}.

Since $\calJ(t_n)-\calJ(t_n+\tau)\to 0$, from the definition for $\calJ$ and the previous estimate we can see that $$\int_{B_1} (| \nabla k_{\vec{q}}|^2(x,t_n)-| \nabla k_{\vec{q}}|^2(x,t_n+\tau))dx\to 0 \hbox{ as }n\to \infty.$$
As above we can conclude that this convergence is almost everywhere and  uniform.
 Standard parabolic estimates imply that also  $k_{\vec{q}}(x,t_n+\tau)-k_{\vec{q}}(x,t_n)$ in the $\calC^2$ norm.

Now we can prove that $k^\infty_{\vec{q}}$ is a solution to the elliptic equation \equ{kinfec}.
Since $k_{\vec{q}}$ solves equation \equ{ecpar'}-\equ{cbpar'}-\equ{cipar'}, we have that
for any $\varphi \in C^\infty (B_1)$ 
$$\int_{B_1} 
(k_{\vec{q}}(y,t_n+1)-k_{\vec{q}} (y,t_n))\varphi(y) dy=
$$
$$
\int_{t_n}^{t_n+1 }\int_{B_1}\left(\Delta  k_{\vec{q}}(y, t_n+\tau) -
\frac{\nabla_v W (k_{\vec{q}}^\infty)}{\epsilon^2}(y, t_n+\tau )- \Delta U_{\vec{q}}\right)\varphi (y)dy d\tau . $$
Letting $n \to \infty$ we get
$$ \int_{B_1}\left(\Delta  k^\infty_{\vec{q}} - \frac{\nabla_v W(k_{\vec{q}}^\infty)}{\epsilon^2}-\Delta U_{\vec{q}}\right)\varphi (y)dy 
=0.$$
Moreover, since for every $t$ holds $k_{\vec{q}}(x,t)|_{\partial B_1}=0$
it must hold $k^\infty_{\vec{q}}|_{\partial B_1}=0$. 
Uniqueness  of solution implies that necessarily $k^\infty_{\vec{q}}\equiv u_\epsilon-U_{\vec{q}}$, which proves the Lemma.
\end{proof}

Now the proof of Theorem \ref{aprox} is direct

\medskip

\begin{proof}[ Proof of Theorem \ref{aprox}]

Fix $\epsilon>0$ and $m>0$. Consider $t_n$ as in Lemma \ref{tainf}, then
\begin{align*}\sup_{B_1}|u_\epsilon-U_{\vec{q}}|\leq &
\sup_{B_1}|u_\epsilon(x)-U_{\vec{q}}(x)-k_{\vec{q}}(x,t_n)|+  
\sup_{B_1\times[0,\infty)}|k_{\vec{q}}(x,t)|\\ \leq & \sup_{B_1}|u_\epsilon(x)-U_{\vec{q}}(x)-k_{\vec{q}}(x,t_n)|+C\epsilon^{m}.\end{align*}
Taking $t_n\to \infty$ we have
\begin{align*}\sup_{B_1}|u_\epsilon-U_{\vec{q}}|\leq C\epsilon^{m}.\end{align*}
Recalling the definition of $U_{\vec{q}}$ 
we have the result. \end{proof}
%\end{proof}

\medskip 

It is easy to see that the size of the radius of the inner ball in Theorem \ref{aprox} 
(that is the ball where $u_\epsilon(x)-u_\sigma\left(\frac{\sigma x}{\epsilon} \right)$ converges to $0$)
can be increased to $\epsilon^\alpha$. Namely,  we 
let
$$\tilde{U}_{\vec{q}}(y)=\tilde{\eta}^\epsilon_{\alpha}(y) \phi_\epsilon (y)+(1-
\tilde{\eta}^\epsilon_{\alpha}(y)) u_\sigma (y),$$
where $\tilde{\eta}:\rr \to \rr$ is 
a positive function such that $\tilde{\eta}(x)=0$ for $|x|\leq 1$ and  $\tilde{\eta}(x)=2$ for $|x|\geq 1$ and
$$\tilde{\eta}^\epsilon_{\alpha}(y) =\tilde{\eta}\left(\frac{1}{2\tilde{E}}|y|+2- \frac{2\epsilon^{\alpha}}{2\tilde{E}}\right),$$   with $\tilde{E}=4 \epsilon^\alpha-\epsilon^{2m+4-\alpha}$.  As before, $\alpha >0$.

Notice that 
$$\tilde{U}_{\vec{q}}(x)=\left\{\begin{array}{cl}\phi_\epsilon(x) & \hbox{ for }|x|\geq 2\epsilon^\alpha \\
u_\sigma\left(\frac{\sigma x}{\epsilon}\right)&\hbox{ for }|x|\leq 2\epsilon^\alpha-E \end{array}\right..
$$
Hence, following the  proof of Theorem \ref{aprox}, but changing $U_{\vec{q}}$ for 
$\tilde{U}_{\vec{q}}$ 
we have
\begin{cor}\label{aprox'}
  Fix $0<\alpha< 1$. 
Let $0<\sigma \leq 2 \epsilon^{1-\alpha}$. 
Then for every $m>0$ there is a constant $C$ (that might depend on $\alpha$ and $m$) such that

\begin{itemize}
\item $$\sup_{|x|\geq 2\epsilon^{\alpha}}|u_\epsilon-\phi_\epsilon|\leq C\epsilon^m.$$

\item $$
\sup_{|x|\leq \epsilon^{\alpha}}\left|u_\epsilon(x)-
u_\sigma \left(\frac{\sigma x}{\epsilon}\right)\right|\leq C \epsilon^m
. $$
\end{itemize}
\end{cor}

Using Lemma  \ref{lema2} and \ref{lema1} we can also prove
\begin{cor}\label{deru}
  Fix $0<\alpha< 1$. 
Let $0<\sigma \leq 2 \epsilon^{1-\alpha}$. 
Then for every $m>0$ there is a constant $C$ (that might depend on $\alpha$ and $m$) such that

\begin{itemize}
\item $$\sup_{|x|\geq \epsilon^{\alpha}}|Du_\epsilon-D\phi_\epsilon|\leq C\epsilon^m.$$

\item $$
\sup_{|x|\leq \epsilon^{\alpha}}\left|Du_\epsilon(x)-
\frac{\sigma }{\epsilon} Du_\sigma \left(\frac{\sigma x}{\epsilon}\right)\right|\leq C \epsilon^m
. $$
\end{itemize}
\end{cor}
\begin{proof} We start by  proving the first inequality of the corollary. To prove this inequality we estimate separately in two different sets, namely we first prove the inequality for $x \in
 B_{1-\frac{\epsilon^\alpha}{2}}\setminus B_{\epsilon^\alpha}$ and then for  $x
\in B_1\setminus B_{1-\frac{\epsilon^\alpha}{2}}$ (in fact, in the second step we find a bound in a larger set:  $B_1\setminus B_{\frac{3}{4}}$) .

We consider the function $u_\epsilon-\phi_\epsilon$ in the domain $B_1\setminus
B_{\frac{\epsilon}{2}}$.
Then
$$\Delta (u_\epsilon-\phi_\epsilon)=\frac{\nabla W(u_\epsilon)-
\nabla W(\phi_\epsilon)}{\epsilon^2}-\Delta \phi_\epsilon+\frac{
\nabla W(\phi_\epsilon)}{\epsilon^2}.$$
Using Lemma \ref{lema2} we have for every 
$x\in B_{1-\frac{\epsilon^\alpha}{2}}\setminus B_{\epsilon^\alpha}$

\begin{align*}|D(u_\epsilon-\phi_\epsilon)|^2 (x)\leq & C\left(\sup_{|x|\geq \frac{\epsilon^\alpha}{2}}\left|u_\epsilon-\phi_\epsilon\right|
\sup_{|x|\geq \frac{\epsilon^\alpha}{2}}
\left|\frac{\nabla W(u_\epsilon)- \nabla W(\phi_\epsilon)}{\epsilon^2}-\Delta \phi_\epsilon+\frac{
\nabla W(\phi_\epsilon)}{\epsilon^2}\right| \right.
\\ &\left.\qquad+\frac{1}{\epsilon^\alpha}\sup_{|x|\geq \frac{\epsilon^\alpha}{2}}|u_\epsilon-\phi_\epsilon|^2
\right)\\
\leq & C\left(\frac{M}{\epsilon^2}\sup_{|x|\geq \frac{\epsilon^\alpha}{2}}\left|u_\epsilon-\phi_\epsilon\right|^2+\left|u_\epsilon-\phi_\epsilon\right|
\sup_{|x|\geq \frac{\epsilon^\alpha}{2}}\left|-\Delta \phi_\epsilon+\frac{
\nabla W(\phi_\epsilon)}{\epsilon^2}\right|\right. 
\\ &
\left.\qquad+\frac{1}{\epsilon^\alpha}\sup_{|x|\geq\frac{\epsilon^\alpha}{2}}|u_\epsilon-\phi_\epsilon|^2
\right).
\end{align*}

Using Theorem \ref{aprox} and the estimates  for $\left|-\Delta \phi_\epsilon+\frac{
\nabla W(\phi_\epsilon)}{\epsilon^2}\right|$ in its proof we have for $m>0$ a constant $C$ (that depends on $m$ and $\alpha$) such that
$$|D(u_\epsilon-\phi_\epsilon)|^2 (x)\leq C\epsilon^m,$$
for  
$x\in B_{1-\frac{\epsilon^\alpha}{2}}\setminus B_{\epsilon^\alpha}$. 

In order to find bounds for $x
\in B_1\setminus B_{1-\frac{\epsilon^\alpha}{2}}$ 
we consider a smooth function $\eta$ such that $\eta(x)\equiv 1$ for $x\geq \frac{3}{4}$ and $\eta\equiv 0$ for $x\leq \frac{1}{2}$ and we consider the function 
$\eta(u_\epsilon-\phi)$ (notice that in fact  this will provide bounds in a larger set, namely $B_1\setminus B_{\frac{3}{4}}$) . Then
$\eta(u_\epsilon-\phi)$ satisfies
$$\Delta(\eta(u_\epsilon-\phi))=\Delta \eta (u_\epsilon-\phi)+\nabla\eta
\nabla(u_\epsilon-\phi)+\eta \left(\frac{\nabla W(u_\epsilon)- \nabla W(\phi_\epsilon)}{\epsilon^2}-\Delta \phi_\epsilon+\frac{
\nabla W(\phi_\epsilon)}{\epsilon^2}\right).$$
Lemma \ref{lema1}, Theorem \ref{aprox} and the previous estimates imply that
$$|D(\eta(u_\epsilon-\phi))|^2(x)=|D(u_\epsilon-\phi)|^2(x)\leq C\epsilon^m\hbox{ for }\frac{3}{4}\leq |x|\leq 1,$$
finishing the proof of the first inequality.

Now we need to prove the second inequality.
Let $u_\sigma^\epsilon(x)=u_\sigma\left(\frac{\sigma x}{\epsilon}\right)$.
To prove the second estimate we consider $u_\epsilon(x)-u_\sigma^\epsilon(x)$ in $B_{\frac{3\epsilon^\alpha}{2}}$.
Since 
$$\Delta(u_\epsilon-u_\sigma^\epsilon)=\frac{\nabla W(u_\epsilon)-
\nabla W(u_\sigma^\epsilon)}{\epsilon^2},$$
Lemma \ref{lema2} implies for every $x\in B_{\epsilon^\alpha}$
\begin{align*}|D(u_\epsilon-u_\sigma^\epsilon)|^2\leq &C\left(\sup_{|x|\leq \frac{\epsilon^\alpha}{2}}\left|u_\epsilon-u_\sigma^\epsilon\right|
\sup_{|x|\leq \frac{\epsilon^\alpha}{2}}
\left| \frac{\nabla W(u_\epsilon)-
\nabla W(u_\sigma^\epsilon)}{\epsilon^2}\right|\right.
\\ &\left. \qquad+\frac{1}{\epsilon^\alpha}
\sup_{|x|\leq \frac{\epsilon^\alpha}{2}}|u_\epsilon-u_\sigma^\epsilon|^2
\right)\\
\leq &C\left(\frac{1}{\epsilon^2}+\frac{1}{\epsilon^\alpha}\right)
\sup_{|x|\leq \frac{\epsilon^\alpha}{2}}
\left|u_\epsilon-u_\sigma^\epsilon\right|^2
\end{align*}

Corollary \ref{aprox'} implies that for every $m>0$ there is a constant $C$ such that
$$|D(u_\epsilon-u_\sigma^\epsilon)|^2\leq C\epsilon^m,$$
which finishes the proof.
\end{proof}
\section{Proof of Theorem \ref{teoprinc}}\label{pf}
%In this section we prove a Theorem equivalent to Theorem \ref{teoprinc}. More specifically, 
Let 
\be v_\epsilon(x)=u_\epsilon(\epsilon x).\label{defveps} \ee
It holds
\begin{align} -\Delta v_\epsilon +\nabla_vW(v_\epsilon)&=0 \hbox{ for  } x\in B_{\frac{1}{\epsilon}}, \\
v_\epsilon(x)&=\phi(x)\hbox{ for } x\in\partial B_{\frac{1}{\epsilon}}.\end{align}

We define the following sequence of continuous function $\tilde{v}_\epsilon :\rr^2\to \rr^2$
\be \tilde{v}_\epsilon(x)=\left\{ \begin{array}{cc}v_\epsilon( x)&\hbox{ for } |x|\leq \frac{1}{\epsilon}\\
\phi(x)& \hbox{ if } |x|\geq \frac{1}{\epsilon}\end{array}.\right.  \label{defvtilde}\ee

We will divide the proof of Theorem \ref{teoprinc} into two different theorems: Theorem \ref{teoprincv1.1} and Theorem \ref{teoprinc2}.
First
we prove

\begin{teo}\label{teoprincv1.1}
There is a subsequence of $\tilde{v}_\epsilon$ such that   $\tilde{v}_\epsilon \to v$ uniformly on compact sets as $\epsilon \to 0$ and $v$ satisfies
\be -\Delta v+\nabla_vW(v)=0 \hbox{ for  } x\in \rr^2, \label{laec}\ee
\be\lim_{|x| \to \infty}|v(x)-\phi(x)|=0. \label{cinf}\ee
\end{teo}
\begin{proof}

Recall first that $\tilde{v}_\epsilon$ is given by \equ{defvtilde}.
We will use the following strategy to prove  Theorem \ref{teoprincv1.1}:

\begin{enumerate}
\item \label{convtildk} Using the results of Section \ref{unifconv}, we show that  $\tilde{v}_\epsilon$ is a Cauchy sequence in the sup norm. Therefore, $\tilde{v}_\epsilon$ has a uniform limit $v$.

\item \label{cikbet} 
Using the definition of $\tilde{v}_\epsilon$ and the first step we show that the limit $v$ satisfies \equ{cinf}.

\item \label{convbto0} Finally, we represent $v_\epsilon$ via Green's formula in compact sets. Taking limits, we conclude that $v$ satisfies \equ{laec}.

%\item \label{comainftildk} Prove that for every $t\in [0,T]$ holds that $$\lim_{|x|\to \infty}|k^\beta_0(x,t)-\phi(x)|=0.$$
%\item\label{convatt0} We observe that pointwise $$\lim_{\beta\to \infty}\tilde{k}_0^\beta(x,0)=v(x).$$
%Moreover, 
%$ \tilde{k}_0^\beta(x,0)\to v(x)$ as $\beta \to \infty$ uniformly on compact sets.
%\item\label{ecsk0} We show that the function $k_0$ satisfies:
%\be Pk^\beta_0+\nabla_vW(k^\beta_0) =0.\label{eck0}\ee
%\ k_0(x,0)&=v(x)\label{cik0}\end{align}
%\item\label{concl} We  observe that 
%\ref{convatt0} and \ref{ecsk0} imply that for every $t>0$
%\be \lim_{\beta \to \infty}\sup_{\rr^2}|k^\beta_0(x,t)-v(x)|=0.\label{limepsbeta}\ee
 
%Finally we conclude the result by using \equ{limepsbeta} and \ref{comainftildk}.
% \ref{convtildk} together with
%\ref{convatt0} imply that for every $t$   holds $$k_0(x,t)=v(x).$$ Which, combined with \ref{comainftildk}, finishes the proof.
\end{enumerate}

Now we prove  these steps:

%\begin{itemize}
\noindent{\bf  Proof of Step \ref{convtildk} :{\em $ \{ \tilde{v}_\epsilon\}$ is a Cauchy sequence in the sup norm.}} 

Consider $\delta>0$ and take $0<\sigma<\epsilon<1$.
We will show that there is an $\epsilon_0$ such that for every
$0<\sigma<\epsilon<\epsilon_0$
$$|\tilde{v}_\epsilon(x)-\tilde{v}_\sigma (x)|\leq \delta
\hbox{ for every }x\in \rr^2.$$
We will mainly  use Theorem \ref{aprox} and Corollary \ref{aprox'}  with $\alpha=\frac{1}{2}$.
\begin{itemize}
\item{\bf If $|x|\leq  \epsilon^{-\frac{1}{2}}$:}

Notice first, that also holds $|x|\leq  \sigma^{-\frac{1}{2}}$ (since $\sigma<\epsilon$).
By the definitions of $\tilde{v}_\epsilon$ and $v_\epsilon$ we have that 
\begin{align*} \tilde{v}_\epsilon(x)-\tilde{v}_\sigma(x)=&u_\epsilon(\epsilon x)-
u_\sigma(\sigma x)\\=& u_\epsilon(y)- u_\sigma\left(\frac{ \sigma y}{\epsilon}\right), \end{align*}
where $y=\epsilon x$. Notice that $|y|=\epsilon|x|\leq \epsilon^{\frac{1}{2}}.$
%Since $\sigma<\epsilon<\epsilon^{\frac{1}{2}}$, 
Corollary \ref{aprox'}  implies that there is a $\epsilon_0$ such that
%constant  $C_{\frac{\delta}{2}}$ such that
%$$\left| h^\alpha_\epsilon(\epsilon x)-\tilde{v}_\sigma\left(\frac{\epsilon x}{\epsilon}\right)\right|\leq C \left(\epsilon^{2-2\alpha} +\epsilon^{1-\alpha}\right)+ C_{\frac{\delta}{2}} \epsilon^{2\alpha}+\frac{\delta}{2}.$$
%Therefore, there is an $\epsilon_0$ such that 
for every $\epsilon<\epsilon_0$
\be |\tilde{v}_\epsilon(x)-\tilde{v}_\sigma (x)|\leq \delta \hbox { for }|x|\leq\epsilon^{-\frac{1}{2}}. \label{cot1teoprinc}\ee

\item{\bf If  $|x|\geq  \epsilon^{-\frac{1}{2}}$ and $|x|\geq  \sigma^{-\frac{1}{2}}$:}

By the definition of $\phi$ and $\phi_\epsilon$ we have that 
$$\phi(x)=\phi_\epsilon(\epsilon x)=\phi_\sigma(\sigma x).$$

This implies 
\be  |\tilde{v}_\epsilon(x)-\tilde{v}_\sigma (x)|\leq  |\tilde{v}_\epsilon(x)-\phi_\epsilon(\epsilon x) |+ |\phi_\sigma(\sigma x)-\tilde{v}_\sigma (x)|.\label{sepor}\ee

If $|x|\geq\epsilon^{-1}$, by definition $\tilde{v}_\epsilon(x)=\phi_\epsilon(\epsilon x) $, hence
\be|\tilde{v}_\epsilon(x)-\phi_\epsilon(\epsilon x) |=0.\label{cota2int0}\ee
For $ \epsilon^{-\frac{1}{2}}\leq|x |\leq\epsilon^{-1}$, by definition 
$\tilde{v}_\epsilon(x)=u_\epsilon(\epsilon x)$. It also holds that $|\epsilon x|\geq \epsilon \epsilon^{-\frac{1}{2}}$. Therefore,
Theorem \ref{aprox}  implies that there is an $\epsilon_1$ such that for every $\epsilon<\epsilon_1$  
\be|\tilde{v}_\epsilon(x)-\phi_\epsilon(\epsilon x) |\leq \frac{\delta}{2} \hbox{ for }x\in B_1 . \label{cota2int1}\ee
Combining  \equ{cota2int0} and \equ{cota2int1} we have that for $\epsilon<\epsilon_1$
\be|\tilde{v}_\epsilon(x)-\phi_\epsilon(\epsilon x) |\leq \frac{\delta}{2} \hbox{ for } |x|\geq\epsilon^{-\frac{1}{2}}. \label{cota2int3}\ee
Since $\sigma<\epsilon<\epsilon_1$ it also holds that
\be|\tilde{v}_\sigma(x)-\phi_\sigma(\sigma x) |\leq \frac{\delta}{2} \hbox{ for } |x|\geq \sigma^{-\frac{1}{2}}. \label{cota2int4}\ee

Equations \equ{sepor}, \equ{cota2int3}and \equ{cota2int4} imply that
\be  |\tilde{v}_\epsilon(x)-\tilde{v}_\sigma (x)|\leq  \delta \hbox { for }|x|\geq\epsilon^{-\frac{1}{2}}, |x|\geq  \sigma^{-\frac{1}{2}}
.\label{cot2teoprinc}\ee

\item{\bf If $\epsilon^{-\frac{1}{2}}\leq |x|\leq  \sigma^{-\frac{1}{2}}$:}

Let us fix any $x$ in this range and define $\tilde{\sigma}=\frac{1}{|x|^2}$. As before, $$\phi(x)=\phi_\epsilon(\epsilon x)=\phi_{\tilde{\sigma}}(\tilde{\sigma} x).$$
Then, we have
$$  |\tilde{v}_\epsilon (x)-\tilde{v}_\sigma (x)|\leq  |\tilde{v}_\epsilon(x)-\phi_\epsilon(\epsilon x) |+ |\phi_{\tilde{\sigma}}(\tilde{\sigma} x)-u_{\tilde{\sigma}} (\tilde{\sigma}x)|$$
\be\hspace{2cm}
+\left|u_{\tilde{\sigma}}(y)-u_\sigma\left(\frac{\sigma y}{\tilde{\sigma}}\right)\right|,\label{sepor2}\ee
where $y=\tilde{\sigma}x$.
As before if $|x|\geq \epsilon^{-1}$, by definition 
\be|\tilde{v}_\epsilon(x)-\phi_\epsilon(\epsilon x) |=0\label{cota113}.\ee
 If $\epsilon^{-\frac{1}{2}}\leq |x|\leq \epsilon^{-1}$, by definition
 $\tilde{v}_\epsilon(x)=u_\epsilon(\epsilon x)$. Hence, 
Theorem \ref{aprox} implies that there is a $\epsilon_2$ such that for every $\epsilon<\epsilon_2$
\be |\tilde{v}_\epsilon(x)-\phi_\epsilon(\epsilon x) |\leq\frac{\delta}{3} .\label{cota123}\ee

Combining \equ{cota113} and \equ{cota123} we have for 
$\epsilon<\epsilon_2$
\be |\tilde{v}_\epsilon(x)-\phi_\epsilon(\epsilon x) |\leq\frac{\delta}{3} \hbox{ for }\epsilon^{-\frac{1}{2}}\leq |x|\leq  \sigma^{-\frac{1}{2}}.\label{cota13}\ee

By the definition of $\tilde{\sigma}$ we have that $|\tilde{\sigma}x|=\frac{1}{|x|}=\tilde{\sigma}^{\frac{1}{2}}$ and $\tilde{\sigma}\leq \epsilon$. Hence, 
using Theorem \ref{aprox}
for $\tilde{\sigma}\leq \epsilon<\epsilon_2$ we have
\be|\phi_{\tilde{\sigma}}(\tilde{\sigma} x)-u_{\tilde{\sigma}} (\tilde{\sigma}x)|\leq\frac{\delta}{3} \hbox{ for }\epsilon^{-\frac{1}{2}}\leq |x|\leq  \sigma^{-\frac{1}{2}}.\label{cota23}\ee

Finally, as $|\tilde{\sigma}x|=\tilde{\sigma}^{\frac{1}{2}}$ and $\sigma\leq \tilde{\sigma}\leq\tilde{\sigma}^{\frac{1}{2}}$, Corollary \ref{aprox'} implies that there is a $\epsilon_3$ such that 
\be \left|u_{\tilde{\sigma}}(\tilde{\sigma} x)-\tilde{v}_\sigma\left(\frac{\tilde{\sigma} x}{\tilde{\sigma}}\right)\right|\leq \frac{\delta}{3} \hbox{ for }\epsilon^{-\frac{1}{2}}\leq |x|\leq  \sigma^{-\frac{1}{2}}.\label{cota33}\ee

Equations \equ{cota13}, \equ{cota23} and \equ{cota33} in \equ{sepor2} imply that 
\be |\tilde{v}_\epsilon(x)-\tilde{v}_\sigma (x)|\leq \delta  \hbox{ for }\epsilon^{-\frac{1}{2}}\leq |x|\leq  \sigma^{-\frac{1}{2}}.
\label{cot3teoprinc}\ee
\end{itemize}

\smallskip

Combining equations \equ{cot1teoprinc}, \equ{cot2teoprinc} and \equ{cot3teoprinc} we conclude that $\tilde{v}_\epsilon:\rr^2\to \rr^2$ is a Cauchy sequence in the sup norm, hence there is a continuous function $v(x)$ such that $\tilde{v}_\epsilon \to v$ uniformly in  $\rr^2$ as $\epsilon\to 0$.

\medskip
\noindent{\bf Proof of Step \ref{cikbet}:  {\em  v satisfies \equ{cinf}}} 
Consider any sequence of points $x_n$ such that $|x_n|\to \infty$. Showing that
$\lim_{n\to \infty}|v(x_n)-\phi(x_n)|=0$ is equivalent to \equ{cinf}.
Let $\epsilon_n=\frac{1}{|x_n|}$. Then for any $\beta>0$ the definition of 
$\tilde{v}_{\epsilon_n}$ implies:
\begin{align*}|v(x_n)-\phi(x_n)|=&|v(x_n)-\tilde{v}_{\epsilon_n}(x_n)|
\\ \leq &\sup_{\rr^2}|v(x)-\tilde{v}_{\epsilon_n}(x)|
\end{align*}
Taking $n\to \infty$, step \equ{convtildk} implies that
$$\lim_{n\to\infty}|v(x_n)-\phi(x_n)|\to 0,$$
which finishes the proof.

\medskip

\noindent {\bf Proof of Step \ref{convbto0}:  {\em  $v$ satisfies \equ{laec}}}

Let us fix a ball of radius $\rho$ in $\rr^2$.
%By Lemma \ref{lema2} we have for every  $v_\epsilon $ with $\epsilon < (\rho+1)^{-1}$ it holds that
% $$|D v_\epsilon(x)|^2\leq C \sup|-\nabla_vW(v_\epsilon)|\sup|v_\epsilon|. $$
% Lemma \ref{cota} implies that there is a $K$ independent of $\epsilon$ such that
%\be|D v_\epsilon(x)|\leq K \hbox{ for }|x|\leq \rho.\label{cotaderveps}\ee
% Hence, by Arzela-Ascoli's theorem there is a subsequence of  $\{v_\epsilon\}$ that converges uniformly in $B_\rho$.
%Since this can be done for any fixed ball, by diagonal argument we can construct  a function $v$ such that there is a subsequence of $\{v_\epsilon\}$ that converges to $v$ 
%uniformly in every compact subset of $\rr^2$. Let us show that $v$ satisfies \equ{laec}.
%Notice that if  $$\epsilon _1, \epsilon_2< \left(\rho +\frac{1}{2}\right)$ then $w=v_{\epsilon_1}-v_{\epsilon_2}$ satisfies
%$$-Delta w=\nabla_vW(v_{\epsilon_1})-\nabla_vW(v_{\epsilon_2})$$
%Hence  Lemma A.1 in \cite{ implies that
%$$|\nabla w|^2\leq C\left( \sup|\nabla_vW(v_{\epsilon_1})-\nabla_vW(v_{\epsilon_2})|\sup|w| + \frac{1}{2}\sup|w|^2\right)$$

In every fixed ball $B_\rho$ we can use Green's formula to represent $v_\epsilon$. We have for 
$\epsilon\leq\frac{1}{\rho}$ that
$$v_\epsilon(x)=-\int_{\partial B_\rho} v_\epsilon(y)\frac{\partial \calK}{\partial \nu}(x,y)dS(y)+
\int_{B_\rho}\nabla_v W(v_\epsilon)(y)\calK(x,y)dy,$$
where $\calK$ is the Green's function in the ball.
Since in $B_\rho$  we have $v_\epsilon \to v$ uniformly as $\epsilon \to 0$, the function $v$ satisfies
$$v(x)=-\int_{\partial B_\rho} v(y)\frac{\partial \calK}{\partial \nu}(x,y)dS(y)+
\int_{B_\rho}\nabla_v W(v)(y)\calK(x,y)dy.$$

Hence,
$$-\Delta v+\nabla_v W(v)=0 \hbox{ for }x\in B_\rho .$$

Since this is true for arbitrary $x$ and $\rho$ we have 
that $v$ satisfies \equ{laec} for every $x\in\rr^2$, which concludes the proof of the Theorem.
%Moreover, since Lemma \ref{cota} implies that $v_\epsilon$ is uniformly 
%bounded,independent of $\epsilon$ we have that there is a constant $C$ such that $|v_\epsilon|\leq C$. This implies, using Theorem \ref{lema2}  that there there is a constant $C$ such that 
%\be |\nabla v|\leq C. \label{cotagradv}\ee
%Now we need to show that \equ{cinf} holds. 
%In order to prove this property we are going to use the sequence $k_{\vec{q}}^\beta$ defined by \equ{defk} in Section \ref{stconv}. Let 
%\be\tilde{v}_\epsilon(xy)=\left\{ \begin{array}{cc}k^\beta_\epsilon( x, t)&\hbox{ for } |x|\leq \frac{1}{\epsilon}\\
%\phi(x)& \hbox{ if } |x|\geq \frac{1}{\epsilon}.\end{array}\right.  \label{defktilde}\ee 
%To show \equ{cinf} we will use the followig strategy:

\end{proof}

Now we finish the Proof of Theorem \ref{teoprinc} by showing
\begin{teo} \label{teoprinc2}
Let
 $$\calV=\left\{w\in C^1:\int_{\rr^2}|Dw-D\phi|dx,\int_{\rr^2}|w-\phi|dx
<\infty\right\}.$$
Define the energy functional
 \be 
\calG(w)=\left\{\begin{array}{cl} \int_{\rr^2}\left(|D w|^2 +W(w) 
-|D\phi|^2-W(\phi)\right)dy &\hbox{ if }w
\in \calV \\
\infty & \hbox{ otherwise. }\end{array} \right. 
\label{deff}\ee

  The energy $\calG$ is bounded below and
the solution $v$ described by Theorem \ref{teoprincv1.1} minimizes $\calG$. That is
$$\calG (v)=\inf_{w\in C^1} \calG (w).$$
\end{teo}

\smallskip
 
\begin{proof}
Define 
\be \tilde{\calG}_\epsilon(w)=\left\{\begin{array}{cl} 
\int_{B_{\epsilon^{-1}}}|D w|^2 +W(w) dy & \hbox{ if }w\in H^1(B_{\epsilon^{-1}}) \hbox{ and } w|_{\partial B_{\epsilon^{-1}}}(x)=\phi_\epsilon(x)
\\ \infty & \hbox{ otherwise. }\end{array} \right. 
\label{defgeps}\ee
and consider $v_\epsilon$ as in the previous Theorem. We will divide the proof of Theorem \ref{teoprinc2} into the following steps:
\begin{enumerate}

\item \label{min}
$v_\epsilon$ is a minimizer for $\tilde{\calG}_\epsilon$ among 
$w_\epsilon\in H^1\left(B_{\epsilon^{-1}}\right)$. 
%satisfying $w_\epsilon|_{\partial B_{\epsilon^{-1}}}=\phi(x)$.
This implies that $v_\epsilon$ minimizes 
$\calG\epsilon(w)=\tilde{\calG}_\epsilon(w)-\tilde{\calG}_\epsilon(\phi)$ in the same class of functions.

\item \label{converg} The sequence
$\calG\epsilon(v_\epsilon)$ is 
%decreasing and bounded below, hence 
convergent.

%\item $|\nabla v_\epsilon|\to |\nabla v|$ pointwise.
\item \label{venv} $v\in \calV$.

\item \label{minv}For every $w$ in $\calV$ there is a sequence $w_\epsilon$ such that
$w_\epsilon\in H^1(B_{\epsilon^{-1}}), 
w_\epsilon|_{\partial B_{\epsilon^{-1}}}(x)=\phi_\epsilon(x)$ and
$\calG\epsilon(w_\epsilon)\to \calG (w)$.

\item \label{convgepsgv}$ \calG_\epsilon(v_\epsilon)\to \calG(v)$.

\item \label{concl} Conclude the result using the previous steps.

\end{enumerate}

\smallskip

\noindent{\bf Proof of Step \equ{min}:}
Notice first that for every $w_\epsilon\in H^1\left(B_{\epsilon^{-1}}\right)$ 
satisfying $w_\epsilon|_{\partial B_{\epsilon^{-1}}}=\phi(x)$ holds that
 $w^\epsilon_\epsilon(x)=w_\epsilon\left(\frac{ x}{\epsilon}\right)\in H^1(B_{1}
)$ and $w^\epsilon_\epsilon|_{\partial B_1}=\phi_{\epsilon}(x).$
Recall that $u_\epsilon$ is a minimizer for $\calI_\epsilon$(defined by \equ{deffeps}), that is
for every $w^\epsilon_\epsilon\in H^1(B_{1}
)$ satisfying $w^\epsilon_\epsilon|_{\partial B_1
}=\phi_{\epsilon}(x)$
holds
$$\calI_\epsilon(u_\epsilon)\leq \calI_\epsilon(w_\epsilon^\epsilon).$$
Dividing by $\epsilon$ and changing variables  holds  

$$\frac{1}{\epsilon}\calI_\epsilon(u_\epsilon)=
\int_{B_{\frac{1}{\epsilon}}}\left(|D v_\epsilon|^2 +W(v_\epsilon) \right)dy\leq
\frac{1}{\epsilon}\calI_\epsilon(w_\epsilon^\epsilon)= \int_{B_{\frac{1}{\epsilon}}}\left(|D w_\epsilon|^2 +W(w_\epsilon)\right) dy,$$
or equivalently
$$\tilde{\calG}_\epsilon(v_\epsilon)\leq \tilde{\calG}_\epsilon(w_\epsilon), \hbox{ for every }w_\epsilon\in H^1\left(B_{\frac{1}{\epsilon}}\right).$$
By subtracting $\tilde{\calG}_\epsilon(\phi)$ we get
$$\calG_\epsilon(v_\epsilon)\leq \calG_\epsilon(w_\epsilon), \hbox{ for every }w_\epsilon\in H^1\left(B_{\frac{1}{\epsilon}}\right).$$

\smallskip

\noindent{\bf Proof of Step \equ{converg}:}
Fix $0<\epsilon< \sigma$.  We need to show that $\calG_\epsilon(v_\epsilon)$ is a Cauchy sequence. Namely, we prove that for every $\delta>0$ there is an $\epsilon_0$ such that $$|\calG_\epsilon(v_\epsilon)-\calG_\sigma(v_\sigma)|\leq \delta \hbox { for } \epsilon,\sigma\leq \epsilon_0.$$
We will  study separately two cases:  $\sigma\geq \sqrt{\epsilon}$ and $\sigma<\sqrt{\epsilon}$.
\begin{itemize}
\item{\bf $\sigma\geq\sqrt{\epsilon}$}
\begin{align*}
|\calG_\epsilon(v_\epsilon)-\calG_\sigma(v_\sigma)|\leq &\left|\int_{B_{\frac{1}{\epsilon}}\setminus B_{\frac{1}{\sigma}}}\left(|Dv_\epsilon|^2-|D\phi|^2+W(v_\epsilon)-W(\phi)\right)dx
\right. \\  &
\left.+\int_{ B_{\frac{1}{\sigma}}}\left(|Dv_\epsilon|^2-|Dv_\sigma|^2
+ W(v_\epsilon)-W(v_\sigma) \right)dx \right| \\ \leq &
\int_{B_{\frac{1}{\epsilon}}\setminus B_{\frac{1}{\sqrt{\epsilon}}}}\left(\left| |Dv_\epsilon|^2-|D\phi|^2\right|+\left|W(v_\epsilon)-W(\phi)\right| \right)dx\\ &+ \int_{ B_{\frac{1}{\sqrt{\epsilon}}}\setminus  B_{\frac{1}{\sigma}}}\left(\left||Dv_\epsilon|^2-|Dv_{\sqrt{\epsilon}}|^2\right|+\left|W(v_\epsilon)-W(v_{\sqrt{\epsilon}})\right|\right)dx\\ &+ \int_{ B_{\frac{1}{\sqrt{\epsilon}}}\setminus  B_{\frac{1}{\sigma}}}\left(\left| |Dv_{\sqrt{\epsilon}}|^2-|D\phi|^2\right|+\left|W(v_{\sqrt{\epsilon}})-W(\phi)\right|\right)dx\\ 
&+\int_{ B_{\frac{1}{\sigma}}}\left(\left| |Dv_\epsilon|^2-|Dv_\sigma|^2\right|+\left|W(v_\epsilon)-W(v_\sigma)\right|\right)dx.
\end{align*}

Let $u_\sigma^\epsilon(x)=v_\sigma\left(\frac{x}{\epsilon}\right)$. Changing variables we have
\begin{align*}
|\calG_\epsilon(v_\epsilon)-\calG_\sigma(v_\sigma)|
\leq &
\int_{B_1\setminus B_{\sqrt{\epsilon}}}\left(\left||Du_\epsilon|^2-|D\phi_\epsilon|^2\right|+\left|\frac{W(u_\epsilon)-W(\phi_\epsilon)}{\epsilon^2}\right|\right)dx
\\ &+ \int_{ B_{\sqrt{\epsilon}}\setminus  B_{\frac{\epsilon}{\sigma}}}\left(\left||Du_\epsilon|^2-\left|Du^\epsilon_{\sqrt{\epsilon}}\right|^2\right|+\left|\frac{W(u_\epsilon)-W(u^\epsilon_{\sqrt{\epsilon}})}{\epsilon^2}\right|\right)dx\\ &+ \int_{ B_{1}\setminus  B_{\frac{\sqrt{\epsilon}}{\sigma}}}\left(\left| |Du_{\sqrt{\epsilon}}|^2-|D\phi_{\sqrt{\epsilon}}|^2\right|+\left|\frac{W(u_{\sqrt{\epsilon}})-W(\phi_{\sqrt{\epsilon}})}{\epsilon}\right|\right)dx
\\ 
&+\int_{ B_{\frac{\epsilon}{\sigma}}}\left(\left||Du_\epsilon|^2-|Du^\epsilon_\sigma|^2\right|+\left|\frac{W(u_\epsilon)-W(u_\sigma)}{\epsilon^2}\right|\right)dx
\end{align*}

Notice that since $\sigma\geq \sqrt{ \epsilon}$ we have that $\frac{\epsilon}{\sigma}\leq \sqrt{\epsilon}$. Then using Theorem \ref{aprox} and Corollaries \ref{aprox'} and \ref{deru} we have that for every $m$ there is a constant, that depends on $m$, such that
\be |\calG_\epsilon(v_\epsilon)-\calG_\sigma(v_\sigma)|
\leq C\epsilon^m.\label{convg1}\ee

\item{\bf $\sigma\leq \sqrt{\epsilon}$}
\begin{align*}
|\calG_\epsilon(v_\epsilon)-\calG_\sigma(v_\sigma)|\leq &\int_{B_{\frac{1}{\epsilon}}\setminus B_{\frac{1}{\sigma}}}\left(\left||Dv_\epsilon|^2-|D\phi|^2\right|+|W(v_\epsilon)-W(\phi)|\right)dx
\\  &
+\int_{ B_{\frac{1}{\sigma}}}\left(\left||Dv_\epsilon|^2-|Dv_\sigma|^2\right|+|W(v_\epsilon)-W(v_\sigma)|\right)dx\\ \leq &
\int_{B_{\frac{1}{\epsilon}}\setminus B_{\frac{1}{\sigma}}}\left(\left| |Dv_\epsilon|^2-|D\phi|^2\right|+|W(v_\epsilon)-W(\phi)|\right)dx\\ &+ \int_{ B_{\frac{1}{\sigma}}\setminus B_{\frac{1}{\sqrt{\epsilon}}}}\left(\left||Dv_\epsilon|^2-|D\phi|^2\right|+|W(v_\epsilon)-W(\phi)|\right)dx\\ &+ \int_{  B_{\frac{1}{\sigma}}\setminus B_{\frac{1}{\sqrt{\epsilon}}} }\left(\left||D\phi|^2-|Dv_\sigma|^2\right|+|W(v_{\sigma})-W(\phi)|\right)dx\\ 
&+\int_{ B_{\frac{1}{\sqrt{\epsilon}}}}\left(\left||Dv_\epsilon|^2-|Dv_\sigma|^2\right|+|W(v_\epsilon)-W(v_\sigma)|\right)dx.
\end{align*}

Let $u_\sigma^\epsilon(x)=v_\sigma\left(\frac{x}{\epsilon}\right)$. Changing variables we have
\begin{align*}
|\calG_\epsilon(v_\epsilon)-\calG_\sigma(v_\sigma)|
\leq &
\int_{B_1\setminus B_{\sqrt{\epsilon}}}\left(\left| |Du_\epsilon|^2-|D\phi_\epsilon|^2\right|+\left|\frac{W(u_\epsilon)-W(\phi_\epsilon)}{\epsilon^2}\right|\right)dx
\\ &+ \int_{ B_{1}\setminus  B_{\frac{\sigma}{\sqrt{\epsilon}}}}\left( \left| |Du_{\sigma}|^2-|D\phi_{\sigma}|^2\right|+\left|\frac{W(u_{\sigma})-W(\phi_{\sigma})}{\sigma^2}\right|\right)dx
\\ 
&+\int_{ B_{\sqrt{\epsilon}}}\left(\left| |Du_\epsilon|^2-|Du^\epsilon_\sigma|^2\right|+\left|\frac{W(u_\epsilon)-W(u^\epsilon_\sigma)}{\epsilon^2}\right|\right)dx.
\end{align*}
Since $\sigma>\epsilon$, we have that $\frac{\sigma}{\sqrt{\epsilon}}\geq\sqrt{\sigma}$. Then,
 Theorem \ref{aprox} and Corollaries \ref{aprox'} and \ref{deru} imply that  for every $m$ there is a constant, that depend on $m$, such that
\be |\calG_\epsilon(v_\epsilon)-\calG_\sigma(v_\sigma)|
\leq C(\epsilon^m +\sigma^m).\label{convg2}\ee

%Consider $\tilde{v}_\epsilon$ as before, defined by \equ{defvtilde}. For every $\sigma\geq \epsilon$ we have  that $\tilde{v}_\sigma$ satisfies
%$\tilde{v}_\sigma\in H^1\left(B_{\epsilon^{-1}}\right)$ and $
%\tilde{v}_\sigma|_{\partial B_{\epsilon^{-1}}}$. Therefore
%$$calG\epsilon(v_\epsilon)\leq \calG_\epsilon(\tilde{v}_\sigma).$$
%By the definitions of $\calG_\epsilon$ and $\tilde{v}_\sigma$ we have
% \begin{align*}\calG_\epsilon(\tilde{v}_\sigma)= &\int_{B_{\epsilon^{-1}}}
%\left(|D \tilde{v}_\sigma|^2 +W(\tilde{v}_\sigma)-|D \phi
%|^2 -W(\phi)\right)dx\\
%= &\int_{B_{\frac{1}{\sigma}}}\left(|D v_\sigma|^2 +W(v_\sigma)-|D \phi
%|^2 -W(\phi)\right)dx=\calG_{\sigma}(v_\sigma).
%\end{align*}
%Hence $\calG_\epsilon(v_\epsilon)$ is increasing in $\epsilon$. 

%Now we show that $\calG_\epsilon(v_\epsilon) $ is bounded below
%acotar las diferencias directamente!

\end{itemize}
We conclude from 
\equ{convg1} and \equ{convg2} that for every $m>0$ there is a constant $C$ such that
$$|\calG_\epsilon(v_\epsilon)-\calG_\sigma(v_\sigma)|\leq C(\epsilon^m+\sigma^m).$$
Therefore $\calG_\epsilon(v_\epsilon)$ is a Cauchy sequence of real numbers, thus convergent.

\medskip

\noindent{\bf Proof of Step \equ{venv}:}
Following the same method of the previous step we can prove the the sequences $\int_{B_{\frac{1}{\epsilon}}}|Dv_\epsilon-D\phi|$ and $\int_{B_{\frac{1}{\epsilon}}}|v_\epsilon-\phi|$ are Cauchy sequences and therefore uniformly bounded.  Fatou's Lemma  implies that
 $$\int_{\rr^2}|Dv-D\phi|dx\leq \int_{B_{\frac{1}{\epsilon}}}|Dv_\epsilon-D\phi|dx<\infty, $$
 $$\int_{\rr^2}|v-\phi|\leq \int_{B_{\frac{1}{\epsilon}}}|v_\epsilon-\phi|dx <\infty.$$
That is $v\in \calV$.

\medskip

\noindent{\bf Proof of Step \equ{minv}:}
Consider a smooth function $\eta$ satisfying $\eta(x)=1$ for $|x|\leq \frac{1}{2}$ and $\eta(x)=0$ for $|x|\geq 1$. Define
$$w_\epsilon(x)=\eta(\epsilon x)w(x)+(1-\eta(\epsilon x))\phi.$$
Then
\begin{align*}|\calG_\epsilon(w_\epsilon)-\calG(w)|=&\left|\int_{\rr^2 \setminus B_{\frac{1}{2\epsilon}}}\left(|Dw|^2-|D\phi|^2+W(w)-W(\phi)\right)dx\right.\\
 &-\int_{B_{\frac{1}{\epsilon}}\setminus B_{\frac{1}{2\epsilon}}}\left(
|\eta Dw+(1-\eta)D\phi+ D\eta(w-\phi)|^2-|D\phi|^2\right)dx\\
& \left.\left. \qquad \qquad+W(\eta(\epsilon x)w(x)+(1-\eta(\epsilon x))\phi)-W(\phi)\right)dx\right|\\
\leq & C\left|\int_{\rr^2 \setminus B_{\frac{1}{2\epsilon}}}\left(|Dw-D\phi|+|w-\phi|\right)dx\right|\\
&+\left|\int_{B_{\frac{1}{\epsilon}}\setminus B_{\frac{1}{2\epsilon}}}\left(
|\eta Dw+(1-\eta)D\phi+ D\eta(w-\phi)-D\phi|\right.\right. \\
&\left.\left.\qquad \qquad +C|\eta(\epsilon x)|w-\phi|\right)dx \frac{}{}\right|  \\
\leq& C\left|\int_{\rr^2 \setminus B_{\frac{1}{2\epsilon}}}\left(|Dw-D\phi|+|w-\phi|\right)dx\right|.\end{align*}
Since $w\in \calV$ we have
$$\lim_{\epsilon\to 0}|\calG_\epsilon(w_\epsilon)-\calG(w)|=0.$$

\medskip

\noindent{\bf Proof of Step \equ{convgepsgv}}

The previous step implies there is a $\tilde{v}_\epsilon$ such that
$$\calG_\epsilon(\tilde{v}_\epsilon)\to \calG(v).$$
Since $v_\epsilon$ is a minimizer of $\calG_\epsilon$ we have that
$$\calG_\epsilon(v_\epsilon)\leq \calG_\epsilon(\tilde{v}_\epsilon).$$
Taking limits when $\epsilon\to 0$ we have
$$\lim_{\epsilon\to 0}\calG(v_\epsilon)\leq \calG(v).$$
In particular, $\calG(v)$ is bounded below.
 Fatou's Lemma allow us to conclude the other inequality:
$$\calG(v)\leq \lim_{\epsilon\to 0}\calG(v_\epsilon).$$

\noindent{\bf Proof of Step \equ{concl}}
Consider $w\in \calV$, then take $w_\epsilon$ as in step \equ{minv}.
Then the minimality of $v_\epsilon$ implies
$$\calG_\epsilon(v_\epsilon)\leq \calG(w_\epsilon).$$
Taking limits as $\epsilon\to 0$ we conclude that
$$\calG(v)
\leq \calG(w),$$
which finishes the proof.
\end{proof}
%\item{ \bf $v$ satisfies \equ{cinf}}

\renewcommand{\theteo}{A-\arabic{teo}}
 \setcounter{teo}{0}  
\renewcommand{\thecor}{A-\arabic{cor}}
 \setcounter{cor}{0}  
\renewcommand{\thelem}{A-\arabic{lem}}
 \setcounter{lem}{0}  
 \renewcommand{\theequation}{A-\arabic{equation}}
  % redefine the command that creates the equation no.
  \setcounter{equation}{0}  % reset counter 

  \section*{APPENDIX}  % use *-form to suppress numbering

In this appendix we present a collection of technical results used in the previous sections.

 We start by stating some results about the Heat Kernel, used mainly in Section \ref{unifconv}. Consider a ball $B\subset\rr^2$. Then $\calH_{B}$
 can be
described as follows:
\begin{align}
(\frac{d}{dt}-\Delta_x)\calH_{B}(x,y,t)=&0,\label{ech}\\
\calH_{B}(x,y,t)=&0 \hbox{ whenever }x\in \partial B,\label{cbh}\\
\lim_{t\to 0^+}\calH_{B}(x,y,t)=&\delta_y(x).\label{cih}
\end{align}
Hence, the solution to the equation
\begin{align*}(\frac{d}{dt}-\Delta_x)u(x,t)=&f(x,t) ,\\
 u(x,t) =&0 \hbox{ whenever }x\in \partial B,\\
u(x,0)=&g(x),
\end{align*}
can be represented as
\be u(x,t)=\int_0^t \int_{B}\calH_{B}(x,y,t-s)f(y,s)dyds +  \int_{B}\calH_{B}(x,y,t)g(y)dy. \label{repform}\ee
We will use this representation to prove the following lemmmas. 
Let us define $P$ to be the heat operator, that is
 \be Pu=\frac{d}{dt} u-\Delta u. \label{defp}\ee
First we  prove some bounds over $\calH_{B}$:

\begin{lem}\label{cotH}
It holds that
\begin{itemize}
\item \be 0\leq \int_{B}\calH_{B}(x,y,t-s)dy ds\leq 1,\label{intcot0}\ee
\item \be 0\leq \int_{s}^t\int_{B}\calH_{B}(x,y,t-s)dy ds\leq (t-s).\label{intcot}\ee
\end{itemize}
\end{lem}
\begin{proof}
The proof follows by maximum principle.
Notice that the single-valued function
$$v(x,t)= \int_{B}\calH_{B}(x,y,t-s)dy ds$$
satisfies the equation
\begin{align}P v(x,t)=&0 ,\label{ecs}\\
 v(x,t) =&0 \hbox{ whenever }x\in \partial B,\label{cis}\\
v(x,s)=&1.\label{cbs}
\end{align}

Since the function 0 is a sub-solution to \equ{ecs}-\equ{cis}-\equ{cbs} we have that 
$$0\leq v(x,t).$$ Similarly,
the function
$1$ is a
 super-solution. Hence,
$$v(x,t)\leq 1,$$
which proves \equ{intcot0}. Equation \equ{intcot} follows by integrating inequality \equ{intcot0}. 

\medskip
 
We also include without proof the following theorem (see \cite{notondg}, 
\cite{theatirm} for example).

\begin{teo}(Theorem 3.1 in \cite{notondg}) \label{H2}
Let $\calM$ be a $n$ dimensional compact Riemannian manifold with boundary. Then there is a Dirichlet heat kernel, that is a function
$$\calH \in C^{\infty}(\calM\times \calM  \times(0,\infty)).$$
satisfying \equ{ech}-\equ{cih}-\equ{cbh}.

%\begin{align*}(\partial_t-\Delta)\calH(x,y,t)=& 0\\
%\calH(x,y,t)=& 0\hbox{ whenever }x\in \partial \calM \\
%lim_{t\to 0^+}\calH(x,y,t)=&\delta_y(x)
%\end{align*}
The smoothness of $\calH(x,x,t)$ may be described as follows
$$\calH(x,x,t)=t^{-\frac{n}{2}}(A(x,t)+B(x,t))$$
with $A\in C^{\infty}(\calM \times [0,\infty))$ and $B$ is supported near the boundary, where in local coordinates $(x',x_n)\in U'\times[0,\tilde{\delta})\subset \calM, \quad U'\subset\rr^{n-1} $ open, one has
$$B(x,t)=b\left(x', \frac{x_n}{\sqrt{t}},t\right), b\in  C^{\infty}(u'\times \rr_+ \times [0,\infty)_{\sqrt{t}})$$ with  $b(x',\psi_n,t)$ rapidly decaying as $\psi_n\to \infty$.

\end{teo}

Now we devote ourselves to prove Lemma \ref{existprobpar}. We start with the following a priori bound: 

\begin{teo}\label{vvc}
Let $\tilde{h}_\epsilon(x,t):\rr^2\to \rr^2$ 
satisfy 
\begin{align} P\tilde{h}_\epsilon+\frac{\nabla_v W(\tilde{h}_\epsilon)}{2}&=0
\hbox{ for } x\in B_{\frac{1}{\epsilon}}\label{ecpar2}\\
\tilde{h}_\epsilon(x,t)|_{\partial  B_{\frac{1}{\epsilon}}}&=\phi(x)\label{cbpar2} \\
\tilde{h}_\epsilon(x,0)&=\psi_\epsilon(x),\label{cipar2}
\end{align}
 where $W:\rr^2\to \rr$ is a  proper $\calC^2$ function,
 bounded below, with a finite number of critical points (denoted by $\{c_i\}_{i=1}^m$),   
%such that $W(v)\to \infty$ as  $|v|\to \infty$ and 
and such that the Hessian
of $W(u)$ is positive semidefinite for $|u|\geq K$, where $K$ is a fixed  real number. % for some real number $K$.
Then if $\tilde{h}_\epsilon(x,0)=\psi_\epsilon(x)$ is bounded 
%and  
 %$\nabla_v W(\psi_\epsilon(x,t))\to 0$ as $|x|\to \infty$ 
 there is a constant 
$C$ that depends only on $W$, $\phi$ and $\psi_\epsilon$ such that $|\tilde{h}_\epsilon(x,t)|\leq C$.
 \end{teo}
 \begin{proof}
Consider $l_\epsilon(x,t)=W(\tilde{h}_\epsilon)(x,t)$; then 
\begin{align*} (l_\epsilon)_t-\Delta l_\epsilon &=
\nabla_v W(\tilde{h}_\epsilon)\cdot (\tilde{h}_\epsilon)_t-\sum_i (\nabla_v W(\tilde{h}_\epsilon) \cdot (\tilde{h}_\epsilon)_{x_i})_{x_i}\\
&= \nabla_v W(\tilde{h}_\epsilon)\cdot (\tilde{h}_\epsilon)_t- 
(W''(\tilde{h}_\epsilon) \nabla \tilde{h}_\epsilon)\cdot \nabla \tilde{h}_\epsilon-\nabla_v W(\tilde{h}_\epsilon) \cdot \Delta \tilde{h}_\epsilon
\end{align*}
where $W''$ denotes the Hessian matrix of $W$.
Since $\tilde{h}_\epsilon$ satisfies \equ{ecpar2}, this becomes
\be (l_\epsilon)_t-\Delta l_\epsilon +\frac{|W'(\tilde{h}_\epsilon)|^2}{2} 
+(W''(\tilde{h}_\epsilon)\nabla u)\cdot\nabla \tilde{h}_\epsilon=0 \label{eqv} \ee

%First, as in the proof of lemma \ref{cw}, w
We are going to find bounds over 
$l_\epsilon$ at the boundary of $B_{\frac{1}{\epsilon}}$ and over its 
possible interior 
maxima in terms of $\max\phi$, $K$, $W(c_i)$ and $\max W(\psi(x))$. 

%If the maximum is attained at the boundary of $B_1$, s
Since for every $|x|=1$ holds $\tilde{h}_\epsilon(x,t)=\phi(x)$  and $\phi$ is uniformly bounded, we have that 
$$l_\epsilon(x)\leq \max W(\phi(x))\hbox{ for every } x\in \partial B_{\frac{1}{\epsilon}}.$$
%Since $W$ is proper this implies that there is a constant $K_1$ such that
%\be|\tilde{h}_\epsilon(x)|\leq K_1 \hbox{ for }x\in \partial B_{\frac{1}{\epsilon}}\label{cotbound}\ee

%Recall that by theorem \ref{v1.0} if
 %$\nabla_v W(\psi_\epsilon(x))\to 0$ as $|x|\to \infty$, 
%we have that $\nabla_v W(\tilde{h}_\epsilon(x,t))\to 0$ as $|x|\to \infty$. 
%We also know that if
 %$W$ has a finite number of critical points, then as $|x|\to \infty$ it holds  $W(\tilde{h}_\epsilon(x,t))\to W(c_i)$ for some critical point $c_i$. Therefore at infinity $v_\epsilon$ is bounded by $\max_i \{ W(c_i)\}$.

Suppose that $l_\epsilon$ has an interior maximum at $(x_0,t_0)$ and 
$|\tilde{h}_\epsilon(x_0,t_0)|\geq K$. 
Since $(x_0,t_0)$ is a maximum for $l_\epsilon$, it holds that
$(l_\epsilon)_t(x_0,t_0)\geq 0$ %(it is 0 if $t_0<\bar{t}$ and non-negative if $t_0=\bar{t}$) 
and $\Delta l_\epsilon(x_0,t_0)\leq 0$. We also have by hypothesis that
$W''(u)$ is positive semidefinite for $|u|\geq K$, hence
$$(l_\epsilon)_t-\Delta l_\epsilon +\frac{|\nabla _u W(\tilde{h}_\epsilon)|^2}{2} +(W''(\tilde{h}_\epsilon)\nabla \tilde{h}_\epsilon)\cdot
\nabla \tilde{h}_\epsilon \geq 0.$$
The inequality is strict (which contradicts \equ{eqv}) 
unless 
$\frac{|\nabla _u W(\tilde{h}_\epsilon)|^2}{\epsilon ^2} =(W''(\tilde{h}_\epsilon)\nabla \tilde{h}_\epsilon)\cdot
\nabla \tilde{h}_\epsilon =0$.
If $\nabla_v W (\tilde{h}_\epsilon(x_0,t_0))=0$, we would have  $\tilde{h}_\epsilon(x_0, t_0) =c_i$ for some $i$, therefore 
$W(\tilde{h}_\epsilon(x,t))\leq W(c_i)$. 
From this and the previous computations we conclude that $l_\epsilon$ is uniformly bounded. 

Since $W$ is a proper function, we have that there is a constant $C$ such that
$$|\tilde{h}_\epsilon(x)|\leq C \hbox{ for }x\in  \bar{B}_{\frac{1}{\epsilon}},$$
which finishes the proof

%\be |\tilde{h}_\epsilon|\leq \max\{K, K_1,a_i,\max_{x\in \rr^2} \psi(x)\},\label{cotavap}\ee  which finishes the proof
 \end{proof}
By observing that solutions to \equ{ecpar}-\equ{cbpar}-\equ{cipar} can be written as
$h_\epsilon(x,t)=\tilde{h}_\epsilon(x,t)-V_{\vec{q}}(x)$, where $\tilde{h}_\epsilon$ is a solution to \equ{ecpar2}-\equ{cbpar2}
-\equ{cipar2} we have

\begin{cor}
Let $h_{\vec{q}}(x,t):B_{\frac{1}{\epsilon}}\to \rr^2$ be a solution to 
\equ{ecpar}-\equ{cbpar}-\equ{cipar}, where $W:\rr^2\to \rr$ is a proper $\calC^2$ function,
 bounded below, with a finite number of critical points and   
%such that $W(v)\to \infty$ as  $|v|\to \infty$ and 
such that the Hessian
of $W(w)$ is positive semidefinite for $|w|\geq K$, where $K$ is a fixed  real number.
Then if $h_{\vec{q}}(x,0)=\psi_\epsilon(x)$ is bounded 
%and  
 %$\nabla_v W(\psi_\epsilon(x,t))\to 0$ as $|x|\to \infty$ 
 there is a constant 
$C$ that depends only on $W$, $\phi$, $U_{\vec{q}}$ and $\psi_\epsilon$ such that $|h_{\vec{q}}(x,t)|\leq C$.
\end{cor}

{\bf Proof of Lemma \ref{existprobpar}}
Let
$$\calC_{[\bar{t}_1,\bar{t}_2]}(B)=\{w:\bar{B}\times [\bar{t}_1,\bar{t}_2]\to \rr^2: w \hbox{
 is a uniformly bounded continuous function }
  \} $$
 with  the standard $\sup$ norm.
  %denoted by $|\cdot|_0$ 
%  and
 %$$B_{[\bar{t}_1,\bar{t}_2]}=\{u:\rr^2\times [\bar{t}_1,\bar{t}_2]\to \rr :W'(u)(x,t)\to 0 \hbox{ as } |x|\to \infty \hbox{ for all } t\in [\bar{t}_1,\bar{t}_2] \}.$$
% Let us first show that $B_{[\bar{t}_1,\bar{t}_2]}$ is a Banach space.
%\begin{lem}
%The space $$B_{[\bar{t}_1,\bar{t}_2]}=\{u:\rr^2\times [\bar{t}_1,\bar{t}_2]\to \rr :W'(u)(x,t)\to 0 \hbox{ as } |x|\to \infty \hbox{ for all } t\in [\bar{t}_1,\bar{t}_2] \}$$ with the standard $C^0$ norm  
%is a Banach space. 
%\end{lem}
%{\bf Proof:}

%Since $B_{[\bar{t}_1,\bar{t}_2]}\subset \calC_{[\bar{t}_1,\bar{t}_2]}$ and  $\calC_{[\bar{t}_1,\bar{t}_2]}$ is a Banach space, we only need to show that $B_{[\bar{t}_1,\bar{t}_2]}$ is closed in the topology of the $C^0$ norm.
%Let $u_n(x,t)\in B_{[\bar{t}_1,\bar{t}_2]}$ be Cauchy, so there exists $u\in  \calC_{[\bar{t}_1,\bar{t}_2]}$ such that
%$$\sup_{(x,t)\in\rr^2\times [\bar{t}_1,\bar{t}_2]}|u_n-u|(x,t)\to 0.$$
%We need to show that $u(x,t) \in B_{[\bar{t}_1,\bar{t}_2]}$.
%Since $u_n$ converges uniformly to $u$ there is a constant $C$ such that $|u|, |u_n|\leq C$ for every $n$. 
%Notice that
%$h_\epsilon-\phi_\epsilon$ satisfies

% \begin{align} P(h_\epsilon-\phi_\epsilon)&=-\frac{\nabla_v W(h_\epsilon)} {2\epsilon^2}+\Delta \phi_\epsilon\\
% h_\epsilon(x)-\phi_\epsilon(x)&=0 \hbox{ for every }x\in \partial\bar{ B_1}\\
% h_\epsilon(x,0)-\phi_\epsilon(x)&=\psi_\epsilon(x)-\phi_\epsilon(x) \end{align}

Consider some $\tau \geq 0$  
and  define $F^\tau_{\vec{q} }(\cdot,\psi^\tau_{\vec{q}}):\calC_{ [\tau ,\tau+\frac{2\beta }{M}]}(B_{\epsilon^{-1}})\to \calC_{ [\tau ,\tau+\frac{2\beta }{M}]}(B_{\epsilon^{-1}})$ by
%\newpage
$$ F^\tau_{\vec{q} }(w, \psi^\tau_{\vec{q}} )(x,t)=
\int_\tau ^t \int_{B_{\epsilon^{-1}}}\calH _{B_{\epsilon^{-1}}}(x,y,t-s) \left(\frac{-W'(w(y,s)+V_{\vec{q}})}{2  }+\Delta V_{\vec{q}}(y)\right) dy ds $$
\be+ \int_{B_{\epsilon^{-1}}}\calH _{B_{\epsilon^{-1}}}(x,y,t) \psi^\tau_{\vec{q}}(y)dy .
\label{F}\ee
%where $\calH(x,t)$ is the heat kernel in $\rr ^2$.
% and  
%the functions %
%$\psi^\tau_{\epsilon}(x)$ will be defined later.

%Lemma \ref{solempezandodesdetau} shows us that 

Notice that Duhamel's formula implies that fixed points of the function $ F^\tau_{\vec{q} }(\cdot,\psi^\tau_{\vec{q}})$ are solutions to \equ{ecpar} in $ [\tau ,\tau+\frac{2\beta }{M}]$.
%if $h_\epsilon$ is a solution to \equ{laeq} such that $h_{\vec{q}}(x,\tau)=\psi^\tau_{\vec{q}}(y)+V_{\vec{q}}$,  then $g_{\vec{q}}=h_{\vec{q}}-V_{\vec{q}}$ is a fixed point  of $F^\tau_{\vec{q} }(\cdot,\psi^\tau_{\vec{q}})$. 
Hence, in order to prove Lemma \ref{existprobpar} we will
use the following strategy: For every $\tau$, $\psi^\tau$ and appropriate constants $\beta, M$ we
find a fixed point of  $F^\tau_{\vec{q} }(\cdot,\psi^\tau_{\vec{q}})$ in some appropriate space; then we choose $\psi^\tau$ appropriately so the fixed points (that were found in the previous step) "glue"  together appropriately; we finish by showing that in fact \equ{ecpar} holds in the whole domain, as well as \equ{cbpar} and \equ{cipar}.% and  estimates for the norm of this fixed point. 

%Notice that for an appropriate  $\psi^\tau_{\vec{q}}$ solutions to \equ{ecpar} are not only fixed points of $ F^\tau_{\vec{q} }$, but also of
 %$F_{\vec{q} }$.
%if $h_\epsilon$ is a solution to \equ{laeq} such that $h_{\vec{q}}(x,\tau)=\psi^\tau_{\vec{q}}(y)+V_{\vec{q}}$,  then $g_{\vec{q}}=h_{\vec{q}}-V_{\vec{q}}$ is a fixed point  of $F^\tau_{\vec{q} }(\cdot,\psi^\tau_{\vec{q}})$. 
%Hence, in order to prove Lemma \ref{existprobpar} we will find a fixed point of  $F^\tau_{\vec{q} }(\cdot,\psi^\tau_{\vec{q}})$ in some appropriate space.% and  estimates for the norm of this fixed point. 
%We will need the following  lemmas.
%Consider th set of uniformly bounded continuous functions:
%$$ \calC_{ [\tau ,\tau+\frac{2\beta \epsilon ^2}{M}]}(B_1)=\{u:
%\bar{ B_1} \times [\tau ,\tau+\frac{2\beta \epsilon ^2}{M}]\to \rr^2:  \sup_{\bar{ B_1} \times [\tau ,\tau+\frac{2\beta \epsilon ^2}{M}]}|u|\leq C \hbox{ for some } C \}.$$

\begin{claim} \label{Fcont} If there is a 
  constant $M$ such that  $|W''|\leq M$ and  $\psi_{\vec{q}}^\tau$ is uniformly bounded, then $F^\tau_{\vec{q}}(\cdot,\psi^\tau_{\vec{q}}):
\calC_{ [\tau ,\tau+\frac{2\beta }{M}]}(B_{\epsilon^{-1}})\to \calC_{ [\tau ,\tau+\frac{2\beta }{M}]}(B_{\epsilon^{-1}})$ is well defined for each $\vec{q}\in \calQ$, where $\calQ$ is given by \equ{calq}. If additionally
 for 
any given $\tau$ and $\beta \in(0,1)$  we have that $\bar{t}$ 
satisfies $|\bar{t}-\tau|\leq \frac{2\beta }{M}$, then 
  $F^\tau_{\vec{q} }(\cdot,\psi^\tau_{\vec{q}})$ is a contraction  mapping
with constant $\beta$ in 
$\calC_{ [\tau ,\tau+\frac{2\beta }{M}]}(B_{\epsilon^{-1}})$.
\end{claim} 
%\begin{proof}
To prove that the function  $F^\tau_{\vec{q} }(\cdot,\psi^\tau_{\vec{q}}):
\calC_{ [\tau ,\tau+\frac{2\beta }{M}]}(B_{\epsilon^{-1}})\to \calC_{ [\tau ,\tau+\frac{2\beta }{M}]}(B_{\epsilon^{-1}})$ is well defined we need to show that $F^\tau_{\vec{q} }(\cdot,\psi^\tau_{\vec{q}})$ maps any uniformly bounded function into a uniformly bounded function, that is for any function $w$ that satisfies
$|w(x,t)|\leq C$ for all $(x,t)\in B_{\epsilon^{-1}} \times [\tau, \bar{t}]$ it holds that 
$|F^\tau_{\vec{q}} (w,\psi^\tau_{\vec{q}})(x,t)|\leq \bar{C}$ for all 
$(x,t)\in B_{\epsilon^{-1}}\times [\tau,\bar{t}]$.

By continuity of $W'$ we have that if  $\sup_{ B_{\epsilon^{-1}} \times [\tau,\bar{t}] }|w(x,t)| \leq C$  then there is a constant $C_1$ such that
$\sup_{(x,t)\in B_{\epsilon^{-1}}\times [\tau,\bar{t}] } |W'(w)(x,t)| \leq C_1$. 
Using the definition of $V_{\vec{q}}$, we can also find constants $C_2$ and $C_3$ that
$$|\Delta V_{\vec{q}}|
\leq 
C_2$$
and
$$|V_{\vec{q}}|\leq C_3$$
This implies 
\begin{align*}|F^\tau_{\vec{q}} & (w,\psi^\tau_{\vec{q}})|(x,t)\\ & \leq
 (C_1+C_2) \int_\tau^{\bar{t}} \int_{B_{\epsilon^{-1}}}\calH_{B_{\epsilon^{-1}}}(x,y,t-\tau-s)dy ds \\
& \quad + \sup_{x\in B_{\epsilon^{-1}}}|\psi^\tau _{\epsilon}(x) | 
 \int_{B_{\epsilon^{-1}}}\calH_{B_{\epsilon^{-1}}}(x,y,t-\tau)dy +C_3 \\
&\leq  (C_1+C_2)(\bar{t}-\tau) + \sup_{x\in B_{\epsilon^{-1}}}|\psi^\tau _{\epsilon} |(x) +C_3=\bar{C}< \infty, \end{align*}
for all $(x,t)$. Hence $F^\tau_{\vec{q} }(\cdot,\psi^\tau_{\vec{q}})$  is well defined for each $\vec{q} \in \calQ$  (where $\calQ$ is given by \equ{calq}).\\
 
Now we show that if $ |\bar{t}-\tau|\leq \frac{2\beta }{M}$, then
$F^\tau_{\vec{q} }(\cdot,\psi^\tau_{\vec{q}})$ is a contraction mapping.

Since $|W''|\leq M$ we have that
$$|W'(w_1)-W'(w_2)|\leq M|w_1-w_2|.$$
Then for every $x\in B_{\epsilon^{-1}}$ and $t \in [\tau, \bar{t}]$ it holds that

$$|F^\tau_{\vec{q}}(w_1,\psi^\tau_{\vec{q}})-F^\tau_{\vec{q}}(w_2,psi^\tau_{\vec{q}})|(x,t) =
\left| \int_\tau ^t \int_{B_{\epsilon^{-1}}}\calH _{B_{\epsilon^{-1}}}(x,y,t-s-\tau) \frac{-W'(w_1(y,s))+W'(w_2(y,s))}{2 } dy ds\right|$$
$$\leq \frac{M (\bar{t}-\tau)}{2 } \sup_{(x,t)\in B_{\epsilon^{-1}}\times [\tau,\bar{t}] } |w_1-w_2|(x,t). \hspace{1cm}$$
Then for $ |\bar{t}-\tau|\leq \frac{2\beta }{M}$ holds
%Choose $\bar{t}=\tau+ \frac{2\beta }{M}$ 
% the previous inequality implies
$$\sup_{(x,t)\in B_{\epsilon^{-1}}\times [\tau,\bar{t}] }|F^\tau_{\vec{q}}(w_1, \psi^\tau_{\vec{q}})-F^\tau_{\vec{q}}(w_2, \psi^\tau_{\vec{q}})|(x,t) \leq \beta\sup_{(x,t)\in B_{\epsilon^{-1}}\times [\tau,\bar{t}] } |w_1-w_2|(x,t)$$
%where $\bar{t}=\tau+\frac{2\beta}{M}$. 
and $F^\tau_{\vec{q} }(\cdot,\psi^\tau_{\vec{q}}):B_{\epsilon^{-1}} 
\times [\tau ,\tau+\frac{2\beta}{M}]\to B_{\epsilon^{-1}} \times [\tau ,\tau+\frac{2\beta }{M}]$ is a contraction with constant 
$\beta$.

We will assume that $|W''|\leq M$ and at the end of the proof we will point out the necessary modifications in the general case. Fix $\beta<1$ and let
%We prove this Theorem using  the following steps:
%\begin{itemize}
%\item 
%We start by proving existence for small time intervals by using the contraction mapping $F^\tau_\epsilon$ defined by \equ{F} . The length of these time intervals depends on $\epsilon$ and $\beta$, the contraction constant of  $F^\tau_\epsilon$ .

%\end{itemize}

%Fix $\bar{t}>0$ and 
 \be \tau_i=i \frac{2\beta }{M}\label{taui}\ee
\be \bar{t}_i=\tau_{i+1}, \label{ti}\ee
\be F_{\vec{q},i}=F^\tau_{\vec{q} }(\cdot,\psi^\tau_{\vec{q}})\ee
 with $i=0,\ldots ,I_\beta,$ where the constant $\beta, I_\beta\in \nn$ satisfy $ \frac{TM}{2 \beta }\leq I_\beta \leq 2 \frac{\bar{t}M}{2 \beta}$. By the definition of $\tau_i, \bar{t}_i$ we have that $\bar{t}_{I_\beta}\geq \bar{t}$. We will redefine  $\bar{t}_{I_\beta}=\bar{t}$.
 
By the previous claim $F_{\vec{q},i}$ is contraction, hence it has a unique fixed point, $h_{\vec{q}} ^i$. That is
\be F_{\vec{q},i}(h^i_{\vec{q}}(x,t))=h^i_{\vec{q}}(x,t). \label{fp}\ee
Moreover, since this this fixed point is bounded we have that $F^\tau_{\vec{q}}(h_{\vec{q}} ^i, \psi^\tau_{\vec{q}})\in C^{1,\frac{1}{2}}(B_{\epsilon^{-1}}\times (\tau_i, \tau_{i+1}])$. Recursively, $h_{\vec{q}}^i \in C^\infty$.
From \equ{fp} and
 Duhamel's formula  we can conclude that \equ{ecpar} and \equ{cbpar}  hold for $t\in [\tau_i,\bar{t}_i]$. We also have
%$h^i_{\vec{q}}=g_{\vec{q}} ^i+\phi_{\vec{q}}$  satisfies
%%%%%\be (h^i_{\epsilon })_t-
%\Delta h^i_{\epsilon } +\frac{1}{2 } W'(h^i_{\epsilon })=0 \label{eqinter}\ee
%\be h^i_\epsilon(x)=\phi_\epsilon(x)\hbox{ for every }x\in \partial\bar{ B_{\epsilon^{-1}}} \label{bcinter}\ee
\be h^i_{\vec{q} }(x,\tau_i)=\psi^{\tau_i}_{\vec{q}}(x)  \label{ciinter}\ee
 for $(x,t) \in B_{\epsilon^{-1}} \times (\tau_i, \bar{t}_i)$.

Now define recursively  $\psi^{\tau_i}_{\vec{q}}(x)$: 
\be \psi^{\tau_0}_{\vec{q}}(x)=\psi_{\vec{q}}(x) \label{psi0}\ee
\be \psi^{\tau_i}_{\vec{q}}(x)=h^{i-1}_{\vec{q} }(x,\tau_i). \label{psii}\ee
Then $h_{\vec{q}}(x,t)$ defined by
\be h_{\vec{q}}(x,t)=h^i_{\vec{q}}(x,t)
\hbox{ for } 
t\in [\tau_i,\bar{t}_i] \label{defueps}\ee
satisfies   \equ{ecpar} for $t\ne \tau_ i$. Moreover, 
by writing
$$h^{i+1}_{\vec{q}}(x,t)=\int_{\bar{t}_{i}} ^t \int_{B_{\epsilon^{-1}}}\calH _{B_{\epsilon^{-1}}}(x,y,t-\bar{t}_{i}-s) 
\frac{-W'(h^{i+1}_{\vec{q}}+V_{\vec{q}})(y,s)}{2 } dy ds $$
$$\hspace{2cm}+ 
 \int_{B_{\epsilon^{-1}}}\calH _{B_{\epsilon^{-1}}}(x,y,t-\bar{t}_{i})h_{\vec{q}}^i(y,\bar{t}_i)  dy,
$$
 standard computations show  that $h_{\vec{q}}$ satisfies \equ{ecpar}
for every $t$. Since
$h_{\vec{q}}$ also satisfies 
\equ{cbpar}-\equ{cipar} we have that $h_{\vec{q}}$ 
is the desired solution. In particular, this implies that $h_{\vec{q}}$ is the fixed point of $F_{\vec{q}}$. Uniqueness follows from the fact that fixed 
points of  contraction mappings are unique. 

In order to prove equation \equ{ineq0} we observe that since 
$h_{\vec{q}}$ is a fixed point of $F_{\vec{q}}^\tau$, standard 
computations imply for any function $w_{\vec{q}}$
\begin{align} |h_{\vec{q}}-w_{\vec{q}}| & \leq \frac{1}{1-\beta}\sup_{B_{\epsilon^{-1}}\times[\tau,\tau+\frac{2\beta}{M}]}  |F_{\vec{q}}^\tau(w_{\vec{q}})-w_{\vec{q}}|
\notag \\
 &\leq \frac{1}{1-\beta}\left(\sup_{B_{\epsilon^{-1}}\times[\tau,\tau+\frac{2\beta}{M}]}|F_{\vec{q}}^\tau(w_{\vec{q}})-F_{\vec{q}}(w_{\vec{q}})|+\sup_{B_{\epsilon^{-1}}\times[\tau,\tau+\frac{2\beta}{M}]}|F_{\vec{q}}(w_{\vec{q}})-w_{\vec{q}}|\right).\label{intereq}\end{align}
The definitions of $F_{\vec{q}}^\tau$ and $F_{\vec{q}}$ imply that 
$$P(F_{\vec{q}}^\tau(w_{\vec{q}})-F_{\vec{q}}(w_{\vec{q}}))=\frac{\nabla_vW(w_{\vec{q}})}{2}-\frac{\nabla_vW(w_{\vec{q}})}{2}=0,$$
and
$$F_{\vec{q}}^\tau(w_{\vec{q}})(x,\tau)-F_{\vec{q}}(w_{\vec{q}})(x,\tau)=
h_{\vec{q}}(x,\tau)-F_{\vec{q}}(w_{\vec{q}})(x,\tau).$$
Using Duhamel's formula we have
$$F_{\vec{q}}^\tau(w_{\vec{q}})-F_{\vec{q}}(w_{\vec{q}})
=\int_{B_{\epsilon^{-1}}}\calH_{B_{\epsilon^{-1}}}(x,y,t-\tau)
(h_{\vec{q}}(y,\tau)-F_{\vec{q}}(w_{\vec{q}})(y,\tau))dy.$$
Together with Lemma \ref{cotH}, this  implies
$$\sup_{B_{\epsilon^{-1}}\times[\tau,\tau+\frac{2\beta}{M}]}|F_{\vec{q}}^\tau(w_{\vec{q}})-F_{\vec{q}}(w_{\vec{q}})|\leq \sup_{B_{\epsilon^{-1}}}|
h_{\vec{q}}(x,\tau)-F_{\vec{q}}(w_{\vec{q}})|(x,\tau).$$
Using \equ{intereq} we conclude inequality \equ{ineq0}
$$|h_{\vec{q}}-w_{\vec{q}}|\leq \frac{1}{1-\beta}\left( 
2\sup_{B_{\epsilon^{-1}}\times[\tau,\tau+\frac{2\beta}{M}]}|F_{\vec{q}}
(w_{\vec{q}})-w_{\vec{q}}|+\sup_{B_{\epsilon^{-1}}}|h_{\vec{q}}-w_{\vec{q}}|(x,\tau)\right)$$

\medskip

For the general case (that is when there is no constant $M$ such that $|W''|\leq M$) we fix $K>0$ large enough. Then we replace $W$ for a function $\tilde{W}$ that satisfies:
\begin{itemize}
\item there is an $M$ such that $|\tilde{W}''|\leq M$,
\item $\tilde{W}(u)=W(u)$ for $u\leq \max\{2 C,K\}$, where $C$ is the constant given by Theorem \ref{vvc}.
 \item $\tilde{W}$ has the same critical points as $W$.
\end{itemize}
Then, our previous computations imply that there is a unique solution 
$h_{\vec{q}} $ to 
\begin{align} Ph_{\vec{q}}+\frac{\nabla_v \tilde{W}(h_{\vec{q}}+V_{\vec{q}})} {2}+\Delta V_{\vec{q}}&= 0\label{eqwtild}\\
 h_{\vec{q}}(x)&=0\hbox{ for every }x\in \partial B_{\epsilon^{-1}} \label{cbtilde}\\
 h_{\vec{q}}(x,0)&=\psi_\epsilon(x).\label{citilde}\end{align} 
Moreover for $w_{\vec{q}}$ as in the hypothesis holds
$$|h_{\vec{q}}-w_{\vec{q}}|\leq \frac{1}{1-\beta}\left( 2\sup_{B_{\epsilon^{-1}}\times[\tau,\tau+\frac{2\beta}{M}]}|\tilde{F}_{\vec{q}}(w_{\vec{q}})-w_{\vec{q}}| +\sup_{B_{\epsilon^{-1}}}|h_{\vec{q}}(x,\tau)-w_{\vec{q}}(x,\tau)|
,\right)$$
where $\tilde{F}_{\vec{q}}$ is analogous to $F_{\vec{q}}$ substituting $W$ for $\tilde{W}$.

However, following the proof Theorem \ref{vvc} we also have that $|h_{\vec{q}}|(x,t)\leq C$, where $C$ is the constant given by Theorem \ref{vvc}. This  fact and the construction of $\tilde{W}$ imply that 
 $h_\epsilon$ is not only a solution to \equ{eqwtild}-\equ{cbtilde}-\equ{citilde}, but
also to 
\equ{ecpar}-\equ{cbpar}-\equ{cipar} (since within this range $W=\tilde{W}$). 
Moreover, for $w_{\vec{q}}$ satisfying
$|w_{\vec{q}}|\leq K$ we will have
$\tilde{F}_{\vec{q}}(w_{\vec{q}})-w_{\vec{q}}=F_{\vec{q}}(w_{\vec{q}})-w_{\vec{q}}$, concluding that \equ{ineq0} holds and finishing the proof of the Theorem.
\end{proof}

\begin{teo}\label{cothder}
Let $h_{\vec{q}}$ be a solution to \equ{ecpar}-\equ{cbpar}-\equ{cipar}, 
then there is a constant $K$, independent of $\vec{q}$, 
such that for every $x\in B_{\frac{1}{\epsilon}}$
\be|Dh_{\vec{q}}|\leq K
.\ee
%Moreover there is a contant $K$ independent $C\leq\max \{$
\end{teo}

\begin{proof}
 Recall that $h_{\vec{q}}$ is vector-valued. We will denote the coordinate $i$-th of the vector $h_{\vec{q}}$ by $h_{\vec{q}}^i$ and, similarly, $\left(\nabla W(h_{\vec{q}})\right)^i$ is the the $i$th coordinate of 
$\nabla W(h_{\vec{q}})$.
We are going to prove separately that for each coordinate that there is a constant $C_i  $ such that $|\nabla  h_{\vec{q}}^i|\leq C_i$.

Let $f:\{(x,y):y\geq 0\}\to B_{\frac{1}{\epsilon}}$ be defined by
\be f(x,y)=\frac{1}{\epsilon}\left(\frac{x^2+y^2-1}{x^2+(y+1)^2},\frac{-2x}{x^2+(y+1)^2}\right).
\label{deftransf}\ee
In complex number notation, we can write for $z=x+iy$ 
$$f(z)=\frac{z-i}{z+i}.$$
Define
 \be s_{\vec{q}}^i(x,y,t)= h_{\vec{q}}^i(f(x,y),t)).\label{defsap}\ee
It satisfies
\begin{align*} \frac{8}{\epsilon(x^2+(y+1)^2)}\frac{d s_{\vec{q}}^i}{dt}-
\Delta s_{\vec{q}}^i&=-\frac{8}{\epsilon(x^2+(y+1)^2)}\left(\nabla W(h_{\vec{q}})\right)^i+\Delta v^i\hbox{ for }x\in \rr, y>0\\
s_{\vec{q}}^i(x,y,t)&=0 \hbox{ for } y=0 \hbox{ or }|(x,y)|\to \infty\\
s_{\vec{q}}(x,y,0)&=0.\end{align*}
Let $\tilde{P}$ be the operator defined by
$$\tilde{P}u=\frac{8}{\epsilon(x^2+(y+1)^2)}\frac{d u}{dt}-
\Delta u.$$

Theorem \ref{vvc}  and the definition of $s_{\vec{q}}^i$ implies that there is a constant $C$ independent of $\epsilon$ such that
$$ \tilde{P}  s_{\vec{q}}^i\leq \frac{C}{\epsilon}.$$
%A straightforward computations shows that there is a $C$ independent of $\epsilon$ such that for every $w\in B_{\frac{1}{\epsilon}}$ holds $|\nabla h_{\vec{q}}(w,0)|\leq C \beta \epsilon^2$. Therefore there is a  constant $C_2$  that satisfies $C_2\leq K  \beta \epsilon^2$ (where $K$ is a constant independent of $\epsilon$) and
Moreover,
\begin{align*} \frac{\partial s^i_{\vec{q}}}{\partial y}(x,0)=0.\end{align*}
Now define
$$w^i_{\vec{q}}(x,y,t)=s^i_{\vec{q}}(x,y,t)-\frac{C}{\epsilon}(y^2+y).$$
%Notice that, in particular, there is a constant $K$ independent of $\epsilon$ such that $C\leq K\max\{1,\epsilon^2 \beta\}$. 
Then
\begin{align*} \tilde{P}w^i_{\vec{q}}&=\tilde{P}s^i_{\vec{q}}-2\frac{C}{\epsilon}
\leq 0,\\
w^i_{\vec{q}}(x,0,t)&= 0 \hbox{ for every } x\in \rr^2 \hbox{ and }t>0,\\
w^i_{\vec{q}}(x,y,0)&<0\hbox{ for }|(x,y)|\to \infty \hbox{ and }y>0.
\end{align*}
Also,

\begin{align*} \frac{\partial w^i_{\vec{q}}}{\partial y}(x,y,0)=
-\frac{C}{\epsilon}(2y+1) \leq 0.\end{align*}
\begin{claim}
The maximum of $w^i_{\vec{q}}$ cannot be attained in the interior.
\end{claim}
%\begin{proof}
If the max is attained at some point in the interior, must hold that
$\Delta w^i_{\vec{q}} <0$
and
$\frac{d w^i_{\vec{q}}}{d t}\geq 0.$
Hence
$\tilde{P} w^i_{\vec{q}}\geq 0,$
which is a contradiction and finishes the proof of the claim.
\medskip

Since the maximum is attained on the boundary it must be attained at $y=0$. Therefore
$$ \frac{\partial w^i_{\vec{q}}}{\partial y}(x,y,t)\leq 0\hbox{ for every }t.$$
This implies that 
$$\frac{\partial s^i_{\vec{q}}}{\partial y}(x,y,t)\leq \frac{C}{\epsilon}(2y+1).$$
This procedure can be repeated for $-s^i_{\vec{q}}$, concluding that
 \be\left|\frac{\partial s^i_{\vec{q}}}{\partial y}(x,y,t)\right|\leq \frac{C}{\epsilon}(2y+1).\label{cotder1ap}\ee
Since the inverse function of $f$ is 
$$f^{-1}(w)=\frac{1+\epsilon w}{1-\epsilon w},$$
using \equ{defsap}, \equ{deftransf} and \equ{cotder1ap} we have (in complex number notation) for any $w\in B_{\frac{1}{\epsilon}}$ that
\be\left|\nabla h_{\vec{q}}^i(w,t)\cdot (1-\epsilon w)^2 \right|\leq 
2C\left(\frac{1-\epsilon^2|w|^2}{1+\epsilon^2|w|^2-\epsilon(w+\bar{w})} +\epsilon\right),\label{cotaderf} \ee
where $\bar{w}$ is the conjugate of $w$.

Similarly, if we define (by performing a rotation of $f$):
\be g(z)=\frac{i}{\epsilon}\frac{z-i}{z+i}\ee and
\be r(x,y,t)= h_{\vec{q}}(g(x,y),t),\ee
following the same method we obtain
\be\left|\nabla  h_{\vec{q}}^i(w,t)\cdot i(1+i\epsilon w)^2 \right|\leq 
2C\left(\frac{1-\epsilon^2|w|^2}{1+\epsilon^2|w|^2+i\epsilon(w-\bar{w})} +\epsilon\right).\label{cotaderg} \ee
Notice  for  $w$ away from $\frac{1}{\epsilon}$ and $\frac{i}{\epsilon}$ holds that $i(1+i\epsilon w)^2 $ and  $(1-\epsilon w)^2 $ are linearly independent as vectors in $\rr^2$.
Fixing some $\delta $ small enough and  considering $w$ such that
$|w-\frac{i}{\epsilon}|\geq \delta$ and $|w-\frac{1}{\epsilon}|\geq \delta$
we have that $\frac{1-\epsilon^2|w|^2}{1+\epsilon^2|w|^2-\epsilon(w+\bar{w})} +\epsilon$ and $\frac{1-\epsilon^2|w|^2}{1+\epsilon^2|w|^2+i\epsilon(w-\bar{w})} +\epsilon$ are bounded above and below independent of $\epsilon$. Hence
\be\left|\nabla h_{\vec{q}}^i(w,t)\right|\leq C\hbox{ for every }
\left|w-\frac{i}{\epsilon}\right|\geq \delta, \ \left|w-\frac{1}{\epsilon}\right|\geq \delta. \ee

Now considering rotation of $f$ of  $\pi$ and $\frac{3}{2}\pi$ radians (that is
 $\tilde{f}(z)=-\frac{1}{\epsilon}\frac{z-i}{z+i}$ and $\tilde{g}(z)=-\frac{i}{\epsilon}\frac{z-i}{z+i}$) and following the same procedure
we fund bounds for $\left|\nabla h_{\vec{q}}^i(w,t)\right|$ near $\frac{1}{\epsilon}$ and $\frac{i}{\epsilon}$, concluding the proof. \end{proof}
%\end{itemize}
%\cite{prueba}
Similarly it follows
\begin{cor}\label{cotkder}
Let $k_{\vec{q}}$ be defined by \equ{defk}. Then there is a constant $K$, independent of $\vec{q}$, 
such that for every $x\in B_1$
\be|Dk_{\vec{q}}|\leq \frac{K}{\epsilon}
.\ee
%Moreover there is a contant $K$ independent $C\leq\max \{$
\end{cor}

\def\cprime{$'$} \def\cprime{$'$} \def\cprime{$'$}

%\begin{thebibliography}{AAA}
%\bibliography{stationary191207}
%\end{thebibliography}

\end{document}